\documentclass[a4paper,11pt]{article}
\usepackage[top=25mm,left=22mm]{geometry}
\usepackage{amsmath,amssymb}
\usepackage{amsthm}
\usepackage{lineno,url} 
\usepackage{graphicx}
\usepackage{amsfonts}
\usepackage{footmisc}
\usepackage{multirow}
\usepackage{pstricks}

\newtheorem{theorem}{Theorem}[section]
\newtheorem{lemma}[theorem]{Lemma}
\newtheorem{corollary}[theorem]{Corollary}
\theoremstyle{definition}
\newtheorem{definition}[theorem]{Definition}

\newtheorem{remark}[theorem]{Remark}
\newtheorem{example}[theorem]{Example}

\def\bexa{\begin{example}}\def\eexa{\end{example}}
\def\brem{\begin{remark}}\def\erem{\end{remark}}
\def\bthm{\begin{theorem}}\def\ethm{\end{theorem}}
\def\blem{\begin{lemma}}\def\elem{\end{lemma}}
\def\bcor{\begin{corollary}}\def\ecor{\end{corollary}}
\def\bdefi{\begin{definition}}\def\edefi{\end{definition}}

\def\bmip{\begin{minipage}{\textwidth}}\def\emip{\end{minipage}}
\def\huga#1{\begin{gather} #1 \end{gather}}

\def\hual#1{\begin{align} #1 \end{align}}
\def\hualst#1{\begin{align*} #1 \end{align*}}

\newcommand{\R}{{\mathbb R}}
\newcommand{\C}{{\mathbb C}}\newcommand{\N}{{\mathbb N}}

\def\CO{{\cal O}}

\def\ga{\gamma}\def\om{\omega}\def\th{\theta}
\def\noi{\noindent}\def\ds{\displaystyle}
\def\pa{{\partial}}\def\lam{\lambda}
\newcommand{\bi}{\begin{itemize}}\newcommand{\ei}{\end{itemize}}
\newcommand{\ben}{\begin{enumerate}}\newcommand{\een}{\end{enumerate}}
\newcommand{\bce}{\begin{center}}\newcommand{\ece}{\end{center}}
\newcommand{\reff}[1]{(\ref{#1})}
\newcommand{\ov}[1]{{\overline {#1}}}

\newcommand{\spr}[1]{\left\langle #1 \right\rangle}
\newcommand{\vs}[1]{{\vspace{#1}}}

\def\aqui{\Leftrightarrow}
\def\ra{\rightarrow}

\newcommand{\barr}{\begin{array}}\newcommand{\earr}{\end{array}}
\newcommand{\bpm}{\begin{pmatrix}}\newcommand{\epm}{\end{pmatrix}}
\newcommand{\bsm}{\left(\begin{smallmatrix}}
\newcommand{\esm}{\end{smallmatrix}\right)}
\newcommand{\ba}{\begin{array}}\newcommand{\ea}{\end{array}}
\def\dd{\, {\rm d}}\def\ri{{\rm i}}

\def\er{{\rm e}}
\def\re{{\rm Re}}\def\im{{\rm Im}}
\def\om{\omega}\def\Om{\Omega}

\def\del{\delta}

\def\eex{\hfill\mbox{$\rfloor$}}

\def\Del{\Delta}

\def\al{\alpha}

\def\Ga{\Gamma}

\def\bd{\begin{displaymath}} \def\ed{\end{displaymath}}
\def\ba{\begin{array}} \def\ea{\end{array}}

\def\matlab{{\tt matlab}}

\def\dds{\frac{\rm d}{{\rm d}s}}\def\rds{{\rm ds}}
\def\ig{\includegraphics}\def\pdep{{\tt pde2path}}\def\auto{{\tt AUTO}}
\def\mlab{{\tt Matlab}}\def\ptool{{\tt pdetoolbox}}
\def\phome{{\tt www.staff.uni-oldenburg.de/hannes.uecker/pde2path}}
\def\bcs{BC}
\def\lamdtol{{\rm tol}_{\dot\lam}}
\def\vKar{von K\'arm\'an}
\def\utt{{\tt u}}\def\Ftt{{\tt F}}
\renewcommand{\arraystretch}{1}\renewcommand{\baselinestretch}{1.0}
\begin{document}
\mbox{}\vspace{0.1cm}\begin{center}\Large
pde2path -- a Matlab 
package for continuation and bifurcation in 2D elliptic systems \\[2mm] 
\normalsize
Hannes Uecker$^1$, Daniel Wetzel$^2$, Jens D.M.\ Rademacher$^3$  \\[2mm]
\footnotesize
$^1$ Institut f\"ur Mathematik, Universit\"at Oldenburg, D26111 Oldenburg, 
hannes.uecker@uni-oldenburg.de\\
$^2$  Institut f\"ur Mathematik, Universit\"at Oldenburg, D26111 Oldenburg, 
daniel.wetzel@uni-oldenburg.de \\
$^3$ Centrum Wiskunde \& Informatica, Science Park 123, 1098 XG Amsterdam, the Netherlands, rademacher@cwi.nl  \\[2mm]
\normalsize
\today

\end{center}
\begin{abstract}\noindent
\pdep\ is a free and easy to use \mlab\ 
continuation/bifurcation package for elliptic systems of PDEs 
with arbitrary many components, on general two dimensional domains, 
and with rather general boundary conditions. 
The package is based on the FEM of the \mlab\ \ptool, 
and is explained by a number of examples, including Bratu's problem, 
the Schnakenberg model, Rayleigh--B\'enard convection, and von Karman 
plate equations. These serve as templates to study new problems, 
for which the user has to provide, via \mlab\ function files, 
a description of the geometry, the boundary conditions, the coefficients 
 of the PDE, and a rough initial guess of a solution. 
The basic algorithm is a one parameter arclength--continuation with 
optional bifurcation detection and branch--switching. 
Stability calculations, error control and mesh-handling, and 
some elementary time--integration for the associated parabolic problem are also supported. 
The continuation, branch-switching, plotting etc are performed via 
\mlab\ command--line function calls guided by the \auto\ style. 
The software can be downloaded from \phome, where also an online 
documentation of the software is provided such that in this paper 
we focus more on the mathematics and the example systems. 
\end{abstract} 
\tableofcontents

\section{Introduction}\label{i-sec}
For algebraic systems, ordinary differential equations (ODEs), and 
partial differential equations (PDEs) in one spatial dimension 
there is a variety of software tools for the numerical 
continuation of families of equilibria and detection and 
following of bifurcations. These include, e.g., \auto\ \cite{auto}, 
{\tt XPPaut} \cite{xppaut} (which relies on \auto\ for the continuation part) and 
{\tt MatCont} \cite{matcont}, 
see also {\tt www.enm.bris.ac.uk/staff/hinke/dss/} for a 
comprehensive though somewhat dated list. 
Another interesting approach is the ``general continuation core'' {\tt coco}, \cite{coco}.

However,  for elliptic systems of PDEs with two spatial dimensions there appear to be few 
general continuation/bifurcation tools 
and hardly any that work out-of-the-box for non-expert 
users.\footnote{PLTMG \cite{pltmg} treats scalar equations, 
and there are many case studies using ad hoc codes, often based 
on \auto\ using suitable expansions for the second spatial direction; 
for 2D systems there also is 
{\tt ENTWIFE}, {\tt www.sercoassurance.com/entwife/introduction.html}, 
which however appears to be no longer maintained since 2001. 
For experts we also mention {\tt Loca}  \cite{loca}, which is 
designed for large scale problems, and {\tt oomph}  \cite{oomph}, 
another large package 
which also supports continuation/bifurcation, though this is not yet 
documented.}
Our software \pdep\ is intended to fill 
this gap.
Its main design goals and 
features are: 
\bi 
\item {\bf Flexibility and versality.} The software treats PDE systems 
\huga{\label{gform}
G(u,\lam):=-\nabla\cdot(c\otimes\nabla u)+a u-b\otimes\nabla u-f=0, 
}
where  $u=u(x)\in\R^N$, $x\in\Omega\subset\R^2$ some 
bounded domain, $\lam\in\R$ is a parameter, 
$c\in\R^{N\times N\times 2\times 2}$,
$b\in\R^{N\times N\times 2}$ (see \eqref{cten}, \eqref{bten} below), 
 $a\in\R^{N\times N}$ and $f\in\R^N$
 can  depend on $x,u,\nabla u$, and, of course, 
parameters.\footnote{The standard assumption is that $c,a,f,b$ depend on 
$u,\nabla u,\ldots$ locally, e.g., $f(x,u)=f(x,u(x))$; however, 
the dependence of $c,a,f,b$ on arguments {\em can} in fact be quite general,  
for instance involving global coupling, see \S\ref{acgc-sec}. 
In particular, we added the $-b\otimes\nabla u$ term to the \ptool--form 
for the effective evaluation of Jacobians, see below.}
The current version supports  ``generalized Neumann'' boundary conditions 
(\bcs) of the form
\begin{align}\label{e:gnbc}
{\bf n}\cdot (c \otimes\nabla u) + q u = g,
\end{align}
where ${\bf n}$ is the outer normal and again $q\in \R^{N\times N}$
and $g\in \R^N$ may depend on $x$, $u$, $\nabla u$ and
parameters. These boundary conditions include zero flux \bcs, 
and large prefactors in $q$, $g$ generate a ``stiff spring''
approximation of Dirichlet \bcs\ that we found to work
reasonably well. 

There are a number of predefined functions 
to specify domains $\Om$ and boundary conditions, or these can be 
exported from \matlab's pdetoolbox GUI, thus making it easy to deal 
with (almost) arbitrary geometry and boundary conditions. 

The software can also be used to time-integrate parabolic 
problems of the form 
\huga{\label{pprob} 
\pa_t u=-G(u,\lam), 
}
with $G$ as in \reff{gform}. This is mainly intended to easily find initial 
conditions for continuation. Finally,  
any number of eigenvalues 
of the Jacobian $G_u(u,\lam)$ can be computed, thus allowing stability 
inspection for stationary solutions of \reff{pprob}. 

\item {\bf Easy usage.} The user has to provide 
a description of the geometry, the boundary conditions, the coefficients 
 of the PDE, and a rough initial guess of a solution. There are a 
number of templates 
for each of these steps which cover some standard cases and should be easy 
to adapt. The software provides a number of \mlab\ functions 
which are called from the command line to perform continuation runs 
with bifurcation detection, branch switching, time integration, etc. 

\item {\bf Easy hackability and customization.} 
While \pdep\ works ``out--of--the--box'' for a significant 
number of examples, already for algebraic equations 
and 1D boundary value problems it is clear that there cannot be a general 
purpose ``solve--it--all'' tool for parametrized problems, see, e.g., 
\cite[Chapter 3]{seydel}. 
Thus, given a particular problem the user might want to customize \pdep. 
We tried to make the data structures and code as 
modular and transparent as possible. When dealing with a tradeoff 
between speed and readability we usually opted for the latter, and thus 
we believe that the software can be easily modified to 
add new features. In fact, we give some examples of ``customization'' 
below. Here, of course, having the powerful \mlab\ machinery 
at our disposal is a great 
advantage. 
\ei

\brem{\rm The $i^{{\rm th}}$ components of 
$\nabla\cdot(c\otimes\nabla u)$, $au$ and $b\otimes \nabla u$ in \reff{gform} 
are given by 
\hual{
&[\nabla\cdot(c\otimes\nabla u)]_i:=
\sum_{j=1}^N [\pa_xc_{ij11}\pa_x+\pa_x c_{ij12}\pa_y+\pa_yc_{ij21}\pa_x
+\pa_yc_{ij22}\pa_y]u_j\, ,\label{cten}\\
&[au]_i=\sum_{j=1}^N a_{ij}u_j,\qquad 
[b\otimes\nabla u]_i:=\sum_{j=1}^N [b_{ij1}\pa_x+b_{ij2}\pa_y]u_j, 
\label{bten}
}
and $f=(f_1,\ldots,f_N)$ should be seen as a column vector. 
If, for instance, we want to implement $-D\Delta u=-(d_1\Delta u_1,\ldots, 
d_N\Delta u_N)=-\nabla\cdot(D\nabla u)$ with $D$ a constant diagonal 
diffusion matrix, as it often occurs in applications, then 
\huga{\label{cdiag}c_{ii11}=c_{ii22}=d_i,\quad i=1,\ldots,N, 
\quad\text{and all other } c_{ijkl}=0, 
}
and there are special ways to encode this (and other symmetric situations 
for $c$ and $a$) in the {\tt pdetoolbox}. See the templates below, 
and \S\ref{coeffsec}. For $c,a,f$ see also the {\tt pdetoolbox} documentation, 
for instance {\tt assempde} in the 
\mlab\ help. 
\pdep\ also provides a simplified encoding for isotropic systems 
without mixed derivatives, see \S\ref{chem-sec}. 
Finally, the $i^{{\rm th}}$ component of 
${\bf n}\cdot (c \otimes\nabla u)$ is given by 
\huga{\label{bcdet}
[{\bf n}\cdot (c \otimes\nabla u)]_i=\sum_{j=1}^N 
[n_1(c_{ij11}\pa_x+c_{ij12}\pa_y)+n_2(c_{ij21}\pa_x+c_{ij22}\pa_y)]u_j, 
}
where ${\bf n}=(n_1,n_2)$, 
\eex}\erem 

\brem{\rm Clearly, the splitting between $a$ and $f$ ($b$ and $f$) 
in \reff{pprob} is not 
unique, e.g., for $G(u)=-\Delta u-\lam u+u^3$ we could use 
($a=-\lam$, $f=-u^3$) or ($a=0$, $f=\lam u-u^3$). 
Similarly, for, e.g., $G(u)=-\Delta u-\pa_x u$ we can use 
$b=(1,0)$ and $f=0$ or $b=(0,0)$ and $f=\pa_xu$. 
This flexibility of \reff{gform} has the advantage that in most cases 
the needed derivatives $G_u, G_\lam$ can be assembled efficiently from 
suitable coefficients $c,a,b$, and no numerical Jacobians are needed. 

Also note that \reff{gform} allows to treat 
equations in nondivergence form, too. For instance, we may write  
a scalar equation $-c(u)\Delta u-f(u)=0$ as 
$-\nabla\cdot(c(u)\nabla u)+(c'(u)\nabla u)\cdot\nabla u-f(u)=0$, and 
set $b_{111}(u)=-c'(u)\pa_x u$ and $b_{112}(u)=-c'(u)\pa_y u$, or 
add $-(c'(u)\nabla u)\cdot\nabla u$ to $f$. 
\eex}\erem 

Currently, the main drawbacks of \pdep\ are:  
\bi 
\item \pdep\ requires \mlab\ including 
the {\tt pdetoolbox}. Its usage explains the form \reff{gform}. 
One of its drawbacks is a somewhat slow performance, compared to, 
e.g., some Fortran implementations of the FEM.\footnote{Another drawback is 
a somewhat unhandy non--GUI description of geometry and boundary 
conditions, but for these we provide fixes. See also, e.g., 
\cite{pruefert}.} On the other 
hand, in addition to the \mlab--environment,  
the {\tt pdetoolbox} has a number of nice features: it also 
takes care of the geometry and mesh generation, 
it is well documented, it is fully based on sparse linear algebra techniques 
(which are vital for large scale problems), it 
exports (sparse) mass and stiffness matrices, and it provides a number of 
auxiliary functions such as 
adaptive mesh-refinement, or various plot options. 
\item Presently, only one parameter continuation is supported, and only bifurcations via 
simple eigenvalues are detected, located, and dealt with\footnote{In case symmetries cause multiple eigenvalues, artificial symmetry breaking sometimes is a viable ad hoc solution for the latter}. We plan to add new 
features as examples require them, and invite every user to do so as well. 
\ei 

In the following we first very briefly recall some basics of 
continuation and bifurcation. Then we explain design and usage of our 
software by a number of examples, 
mainly a modified Bratu problem as a standard scalar elliptic 
equation, some Allen--Cahn type equations, 
some pattern forming Reaction--Diffusion systems, including 
some animal coats intended for illustration of how to set up problems with 
complicated geometries. We give a rather detailed bifurcation diagram 
for the Schnakenberg system, and we consider three rather classical 
problems from physics: Rayleigh--B\'enard convection, 
some multi--component Bose--Einstein systems, and the \vKar\ plate equations. 
Thus, besides some mathematical aspects of 
continuation and the example systems, here we explain the syntax and usage of the software in a 
rather concise way. More comprehensive documentation of the 
data structures and functions is included in the software, or online at 
\cite{pde2phome}. 
\\[2mm]
\noi
{\bf Acknowledgements.} We thank Uwe Pr\"ufert for providing 
his extension of the \ptool\ and the documentation. 
Users familiar with \auto\ will recognize that \auto\ has been our guide in many 
respects, in particular concerning the design of the user interface. 
We owe a lot to that great software. We also thank Tomas Dohnal for  
testing early versions of \pdep\ and providing valuable hints for making 
\pdep\ and this manual more user friendly.  
JR acknowledges support by the NDNS+ cluster of the Dutch Science Fund NWO.

\brem 
Please report any bugs to {\tt pde2path@mathematik.uni-oldenburg.de}, 
as well desired additional features. We will appreciate any feedback, 
and will be happy to provide help with interesting applications. 
See also \cite{pde2phome} for an online 
documentation of the software, updates, FAQ, and general 
further information. 
\erem

\section{Some basics of continuation 
and bifurcation}
\subsection{Arclength continuation}\label{gtsec}
A standard method for numerical calculation of solution 
branches of $G(u,\lam)=0$, where $G:X\times\R\ra X$ is at least $C^1$, 
$X$ a Banach space, is (pseudo)arclength 
continuation.  For convenience and reference here we recall the basic ideas. 
Standard textbooks on continuation and bifurcation are 
\cite{seydel,govaerts,kuz,AllGeo90}, see also \cite{keller77,doedel07}, 
and the ``matrix-free'' approach \cite{Geo01}. 
Consider a branch $z(s):=(u(s),\lam(s))\in X\times\R$ parametrized by  
$s\in\R$ and the extended system 
\huga{\label{esys} 
H(u,\lam)=\bpm G(u,\lam)\\p(u,\lam,s)\epm=0\in X\times \R, 
}
where $p$ is used to make $s$ an approximation to arclength on the 
solution arc. 
Assuming that $X$ is a Hilbert space with inner product $\spr{\cdot,\cdot}$, 
the standard choice is as follows: given 
 $s_0$ and a point $(u_0,\lam_0):=(u(s_0),\lam(s_0))$, and 
additionally knowing a tangent vector 
$\tau_0:=(\dot u_0,\dot\lam_0):=\dds (u(s),\lam(s))|_{s=s_0}$ we use, 
for $s$ near $s_0$, 
\huga{
p(u,\lam,s):=\xi\spr{\dot u_0,u(s)-u_0}+(1-\xi)\dot\lam_0(\lam(s)-\lam_0)
-(s-s_0). 
}
Here $0<\xi<1$ is a weight, and $\tau_0$ is assumed to be normalized in 
the weighted norm 
$$
\|\tau\|_\xi:=\sqrt{\spr{\tau,\tau}_\xi}, \qquad 
\spr{\bpm u\\ \lam\epm,\bpm v\\ \mu\epm}_\xi:=\xi \spr{u,v}+(1-\xi)
\lam\mu. 
$$
For fixed $s$ and $\|\tau_0\|_\xi=1$, 
$p(u,\lam,s)=0$ thus defines a hyperplane perpendicular 
(in the inner product $\spr{\cdot,\cdot}_\xi$) 
to $\tau_0$ at distance $\rds:=s-s_0$ from $(u_0,\lam_0)$. We may then 
use a predictor 
 $(u^1,\lam^1)=(u_0,\lam_0)+\rds\,\tau_0$ for a solution \reff{esys} 
on that hyperplane, followed by a corrector using Newton's method in the form 
\huga{\label{newton} 
\bpm u^{l+1}\\
\lam^{l+1}\epm =\bpm u^{l}\\\lam^{l}\epm - A(u^l,\lam^l)^{-1}H(u^l,\lam^l),\quad 
\text{ where } A=\bpm G_u&G_\lam\\
\xi \dot u_0&(1-\xi)\dot\lam_0\epm.}
Since $\pa_s p=-1$, on a smooth solution arc we have 
\huga{
A(s)\bpm \dot u(s)\\\dot\lam(s)\epm=-\bpm 0\\\pa_s p\epm=\bpm 0\\ 1\epm. }
Thus, after convergence of \reff{newton} yields a new point $(u_1,\lam_1)$ 
with Jacobian $A^1$, the tangent direction $\tau_1$ at $(u_1,\lam_1)$ 
with conserved orientation, i.e., $\spr{\tau_0,\tau_1}=1$, 
can be computed from 
\huga{\label{tau1} A^1\tau_1=\bpm 0\\ 1\epm, \text{ with 
 normalization } \|\tau_1\|_\xi=1. 
} 
 Alternatively to \reff{newton} 
we may also use a chord method, where $A=A(u^1,\lam^1)$ is kept fixed 
during iteration, 
\huga{\label{chord} 
\bpm u^{l+1}\\
\lam^{l+1}\epm =\bpm u^{l}\\\lam^{l}\epm - A(u^1,\lam^1)^{-1}H(u^l,\lam^l).
}
 This avoids the costly evaluation of $G_u$ at the price 
of a usually modest increase of required iterations. 

\begin{figure}[ht!]
\bce 
\ig[width=50mm]{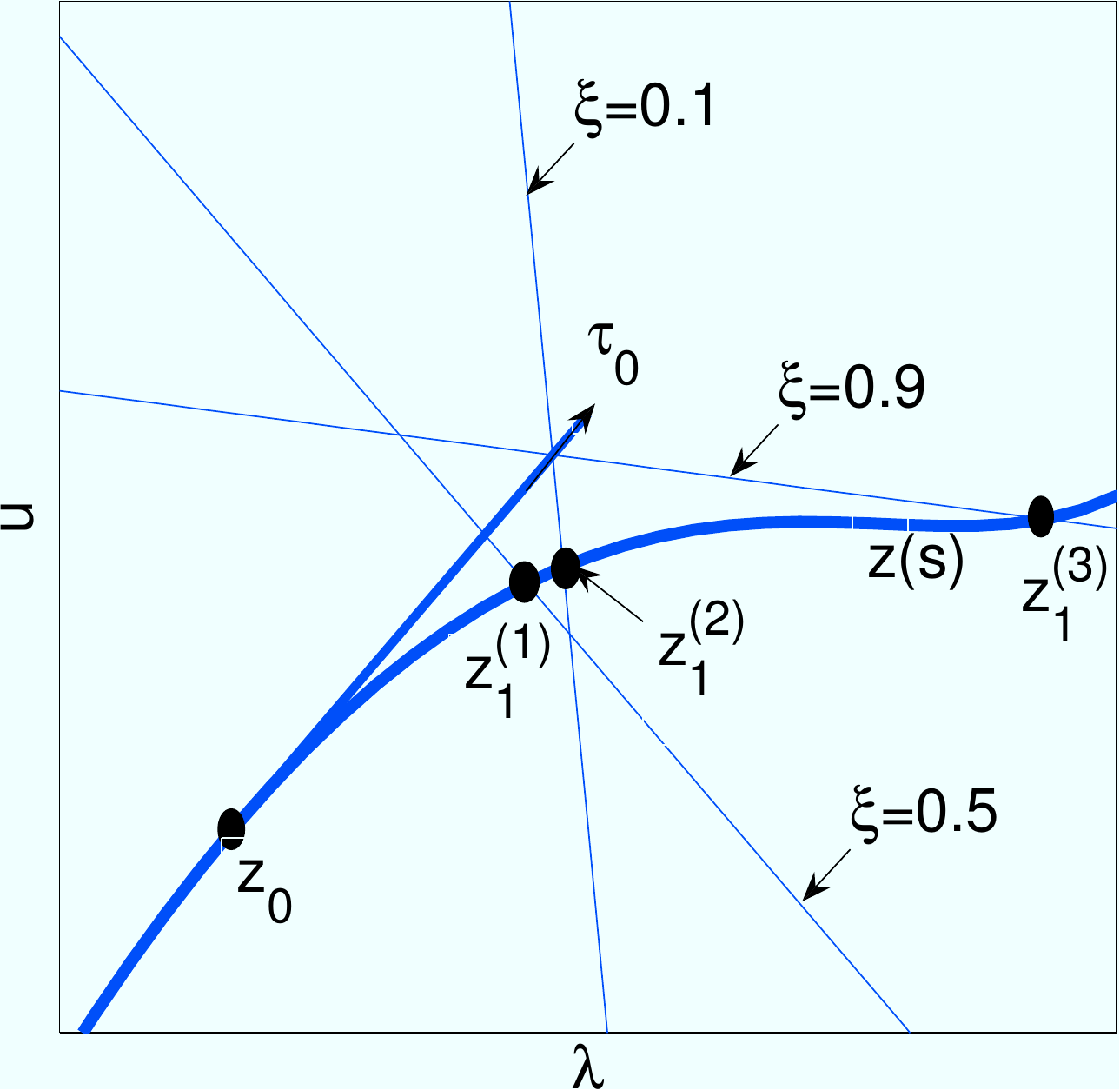}
\ece

\vs{-6mm}
\caption{{\small The role of $\xi$ in a one--dimensional ($u\in\R$) 
sketch with 
$\tau_0=(\dot u_0,\dot\lam_0)=(1,1)$ (unnormalized). Depending on $\xi$ we get 
different hyperplanes $\{(u,\lam)\in\R^2: \spr{\tau_0,(u,\lam)}_\xi={\rm ds}\}$ 
and consequently different ``next points'' $z_0^{(i)}$, $i=1,2,3$, on the 
solution curve $z(s)$. Small $\xi$ favors Newton search in $u$ direction 
(and thus orthogonal to a ``horizontal'' branch), while 
large $\xi$ favors the $\lam$ direction (parallel to a ``horizontal'' 
branch).  \label{xif}}} 
\end{figure}

The role of $\xi$ is twofold.\footnote{Here $\utt$ stands for the FEM 
approximation of $u$, and {\tt G(u,$\lam$)} for the FEM approximation 
of $G$. Below, the difference between the two  will be clear from the context, 
and thus we will mostly drop the different notations again.\label{fn:FEM}} 
 First, if ${\tt G}(\utt,\lam)=0$ comes from the 
discretization of a PDE $G(u,\lam)$ such as \reff{gform} 
over a domain $\Om$ with $n_p$ spatial points, 
then $\utt\in\R^p$ with large $p$, say $p=Nn_p$. 
Typically, we want to choose $\xi$ such that 
$\xi\|\utt\|_{\R^p}^2$ is an approximation of $\ds\frac 1{|\Om|}\|u\|_{L^2(\Om)}^2$. 
If (as usual), $u\equiv 1$ corresponds to $\utt_j=1$ for $j=1,\ldots,n_p$, 
then a rough estimate can be obtained by assuming that 
each component $u_i\equiv 1$, $i=1,\ldots,N$. Then 
\huga{\label{xiform}
\frac 1 {|\Om|}\|u\|_{(L^2)^N}^2=N\stackrel{!}{=}\xi \|\utt\|_{\R^p}^2=\xi N n_p, 
\text{ hence } \xi=1/n_p. 
}
This gives the basic formula for our choice of $\xi$. 
It is important that different $\xi$ may give different continuations: in the 
Newton loop, small $\xi$ favors changes in $u$, while larger $\xi$ 
favors $\lam$, see Fig.\ref{xif} for a sketch. 
Moreover, $\xi$ is also related to the scaling 
of the problem: if, e.g., we replace $\lam$ by, say, $\tilde{\lam}:=100\lam$, 
then $\xi$ should be adapted accordingly, i.e., $\tilde{\xi}=\xi/100$. 
In summary, $\xi$ should be considered as a parameter 
that can be used to tune the continuation, and that may also be changed 
during runs if appropriate.

Given a weight $\xi$, a starting point 
$(u_0,\lam_0,\tau_0)$,  and an intended step size $\rds$, the 
basic continuation algorithm thus reads as follows, already including 
some elementary stepsize control:\\[2mm]
\fbox{\parbox{160mm}{
{\bf Algorithm {\tt cont}}\\[-6mm]
\ben
\item {\bf Predictor.}  
Set $(u^1,\lam^1)=(u_0,\lam_0)+\rds\tau_0$. 
\item {\bf Newton--corrector.} Iterate \reff{newton} (or \reff{chord}) until convergence; 
decrease $\rds$ if \reff{newton} 
fails to converge and return to 1; increase $\rds$ for the next step 
if \reff{newton} converges quickly; 
\item {\bf New tangent.} Calculate $\tau_1$ from \reff{tau1}, set 
$(u_0,\lam_0,\tau_0)=(u_1,\lam_1,\tau_1)$ and return to step 1.
\een 
}}\\[2mm]
Theoretically, this does not work at possible ``bifurcation points'' 
where $A$ is singular.\footnote{Allthough generically continuation 
routines simply shoot past singular points.} 
More specifically, we {\em define:} \\[1mm]
{\bf B1.} A {\em 
simple bifurcation point} is a point $(u,\lam)$ where det$A$ changes 
sign. The implicit assumption is that this happens due to a simple 
eigenvalue of $A$ crossing zero.  

This clearly excludes folds (aka turning points), where a simple 
eigenvalue of $A$ reaches zero but det$A$ does not change sign, 
see \cite{keller77}.   
However, folds are no problem for the algorithm and 
can easily seen in the bifurcation diagram anyway. 
Therefore there 
is no special treatment of folds in the current 1-parameter version of \pdep. 
{\bf B1} also excludes 
bifurcations via even numbers of eigenvalues crossing, which are 
more complicated to deal with.  
\brem Numerically, for {\bf B1} we found it more robust to use $\xi=1/2$ 
 in the definition of $A$ for bifurcation purposes. For algorithmic 
reasons, we only use the first part of {\bf B1} for the detection 
of bifurcation points, i.e., 
the sign change of det$A$, which also occurs for an odd number 
of eigenvalues (counting multiplicities) crossing zero. Finally, 
to calculate sign(det\,$A$) we only calculate the set $\Sigma_0$ 
of eigenvalues $\mu_i, i=1,\ldots,n_{{\rm eig} }$ of $A$ closest 
to $0$ and then use
\huga{\label{detaform}
\text{sign(\,det\,$A$)=sign$(\ds \Pi_{i=1}^{n_{{\rm eig}}}\re\mu_i)$.}
}
Here the implicit assumption is that $n_{{\rm eig} }$ is sufficiently 
large such that all eigenvalues with negative real part are 
always contained in $\Sigma_0$. 
This is reasonable as $G_u+\gamma$ is a positive elliptic operator for 
sufficiently 
large $\gamma$.\footnote{Due to occasional numerical problems in the eigenvalue 
calculations, in the current standard setting 
of \pdep\ (controlled by switches, see below) we actually 
combine the sign-change of det$(A)$ with a consistency check with 
the eigenvalues of $G_u$. Moreover, in practice it is sometimes useful to 
turn off branch point 
detection and localization in the initial continuation (again by setting switches, see below),
but only monitor changes in the number of eigenvalues with negative real part 
in $G_u$. For each such change select a nearby starting point for a new 
continuation with branch detection. See also \S\ref{s:rbconv}. Finally, 
we also provide a routine {\tt findbif} which first scans a branch for a change of 
the number of unstable eigenvalues, and then uses {\bf B1} to locate 
a bifurcation. 

Our standard setting is $n_{{\rm eig}}=50$, but of course this is 
highly problem dependent and should be adapted by the user when needed. 
We give a warning if $|\mu_1|>|\mu_{{\rm neig}}|/2$ since 
then eigenvalues might wander out of $\Sigma_0$ to the left in the next steps.
}
\eex 
\erem 

After {\em detection} of a bifurcation between $s_k$ and $s_{k+1}$, 
the bifurcation is {\em located} by a bisection method. 
To switch branches we use ``Method I'' of \cite{keller77} (page 379). 
Let $(u_0,\lam_0)$ be a simple bifurcation point, $G_u=G_u(u_0,\lam_0)$, 
and $\tau_0=(\dot u_0,\dot\lam_0)$ be the tangent along the branch 
already computed. To obtain a tangent $\tau_1$ along the other 
branch we proceed as follows:\\[2mm]
\parbox{160mm}{
\fbox{\parbox{160mm}{
{\bf Algorithm {\tt swibra}}\\[-6mm] 
\ben
\item Calculate $\phi_1, \psi_1$ with $G_u^0\phi_1=0, {G_u^0}^T\psi_1=0$, 
$\|\phi_1\|=1$, $\spr{\psi_1,\phi_1}=1$. \\[-7mm]
\item Let $\al_0=\dot\lam_0$, $\al_1=\psi^T\dot u_0$, 
$\phi_0=\al_0^{-1}(\dot u_0-\al_1\phi_1)$. \\[-7mm]
\item Choose some small $\del>0$ and calculate the finite differences 
\hualst{a_1&=\frac{1}{\del}\psi_1^T
\bigl[G_u(u+\del\phi_1,\lam_0)-G_u^0\bigr]\phi_1, \\
b_1&=\frac{1}{\del}\psi_1^T
\biggl[\bigl[G_u(u_0+\del\phi_1,\lam_0)-G_u^0\bigr]\phi_0
+G_\lam(u_0+\del\phi_1,\lam_0)-G_\lam(u_0,\lam_0)\biggr]. 
}
Assuming $\al_0\ne 0$  (see \cite{keller77} if this is not the case), set 
$$\ds\ov{\al}_1=-\left(\frac{a_1\al_1}{\al_0}+2b_1\right) 
\text{ and }
\ds \tau_1=\bpm \ov{\al}_1\phi_1+a_1\phi_0\\a_1\epm. 
$$ 
Choose$^a$ a weight $\xi$ and a stepsize $\rds$, 
set $\tau_0=\tau_1/\|\tau_1\|_\xi$ and go to {\tt cont}, step 1.  \\[-7mm]
\een 
}} 
$^a$ {\footnotesize If branch-switching fails, i.e., 
if there is no convergence in {\tt cont} or if 
the solution falls back onto the known branch, then it may help  
to change {\tt ds} and/or $\xi$.}
} 

\subsection{Switching back and forth to the natural 
parametrization} 
If $\dot{\lam}:=\pa_s \lam$ does not change sign, then we know 
that a branch also has the 
``natural parametrization'' $(u(\lam),\lam)$, and,  
except at possible bifurcation points, 
$G_u(u,\lam)u'(\lam)=-G_\lam(u,\lam)$ has the unique solution 
$$
u'(\lam)=-G_u(u,\lam)^{-1}G_\lam(u,\lam). 
$$
Thus, given $(u_0,\lam_0)$ we may use the predictor 
 $(u^1,\lam_1)=(u_0,\lam_0)+\rds\,(u'(\lam_0),1)$ and then 
correct with fixed $\lam=\lam_1$. 
Algorithmically, however, we choose 
to keep the predictor $(u^1,\lam_1)=(u_0,\lam_0)+\rds\,\tau_0$, 
i.e., altogether, 
\huga{\label{trn}
(u^1,\lam_1)=(u_0,\lam_0)+\rds\,\tau_0, \qquad 
u^{l+1}=u^l-G_u(u^l,\lam_1)^{-1}G(u^l,\lam_1).  
}
After convergence to $(u_1,\lam_1)$ we 
calculate the new tangent via \reff{tau1}, and 
{\bf B1} can again be used as a check for bifurcation.
Moreover, with $\lamdtol>0$ a given tolerance, say $\lamdtol=0.5$, 
this gives a criterion when to switch back and forth between 
the algorithms, namely: 
\huga{\label{csw}
\text{If $|\dot\lam|>\lamdtol$, then use \reff{trn}, else use \reff{newton}.}
} 
Here again the weight $\xi$ is important: 
for fixed $\xi=1/2$ (say), $\dot\lam\ra 0$ as $n_p\ra \infty$ unless 
the branch is strictly horizontal, i.e., $\dot u=0$. 

\brem If applicable, \reff{trn} is usually slightly faster than 
\reff{newton}, as expected. On the other hand, 
we found that even for ``nearly horizontal'' branches, locating the bifurcation point 
typically works better with arclength continuation \reff{newton}. 
\eex\erem 

\section{Some scalar problems in \pdep}\label{s:scalar}
We now start the tutorial on \pdep\ by way of 
basic examples. The names in brackets refer to the sub-directory name 
of the directory {\tt demos}, which contains the given example. 

\subsection{Bratu's problem ({\tt bratu})}\label{br-sec}
Our first example is the scalar elliptic equation 
\huga{\label{brprob}
-\Delta u-f(u,\lam)=0, \quad f(u,\lam)=-10(u-\lam \er^{u}), \quad 
u=u(x)\in\R,
}
on the unit square with zero flux \bcs, i.e.,  
\huga{\label{brbc}
x\in \Om=(-1/2,1/2)^2,\quad \pa_n u|_{\pa\Om}=0.  
}
This problem has 
the advantage that a number of results can immediately be obtained analytically, 
that there are some nontrivial numerical questions (see below), and that we
 can compare with previous results, see, e.g., \cite{pltmg,bratu14,mittel86}

There is a primary homogeneous branch $u\equiv u_h(s)$, $\lam=\lam(s)$, 
``starting'' in $(0,0)$, on which $(u_h(s),\lam(s))$ satisfies 
$f(u)=-10(u-\lam \er^u)=0$. 
Bifurcation points $(u_k,\lam_k)$ are obtained from 
$G_uw=-\Delta w -f_uw\stackrel{!}{=}0$ which yields $10(u_k-1)=\mu_k$ where 
$\mu_k=(k_1^2+k_2^2)\pi^2$, $k\in\N_0^2$, are the eigenvalues 
of $-\Delta$ on $\Om$, see Table \ref{tab1}. 
\begin{table}[!ht]\bce
\begin{tabular}{|c|c|c|c|c|c|c|}
\hline
k&(0,0)&(1,0),(0,1),&(1,1)&(2,1),(1,2)&(2,2)&\ldots\\\hline
$u_k$&1 &$1+\pi^2/10$ & $1+\pi^2/5$ &$1+\pi^2/2$ & $1+4\pi^2/5$ &\ldots\\
$\lam_k=u_k\er^{-u_k}$&$1/e\approx 0.3679$ &$\approx 0.2724 $&$\approx 0.1520  $
&$\approx 0.0157$&$\approx0.0012$&\ldots\\
type& fold& double& simple& double& simple&\ldots\\
\hline
\end{tabular}
\caption{{Bifurcation from homogeneous branch in \reff{brprob}.}\label{tab1}}
\ece
\end{table}
From \S\ref{gtsec}, for arclength continuation 
the fold is nothing special and the two dimensional kernels 
$k=(k_1,k_2)$, $k_1\ne k_2$, will go undetected using B1. 
The simple bifurcation points should be detected, and branch switching 
can be tried. 
On bifurcating branches, further bifurcations may be expected. 

In \pdep, \reff{brprob},\reff{brbc} can be setup and run in a few steps 
explained  next, to quickly 
obtain the (basic) bifurcation diagram and solution plots in Fig.~\ref{bratubdf} 
on page \pageref{bratubdf}. 

\subsubsection{Installation and preparation}
The basic \pdep\ installation consists of a root directory, 
called {\tt pde2path}, with a subdirectory {\tt p2plib} 
containing the actual software, 
a subdirectory {\tt demos} with a further subdirectory for each problem, 
a subdirectory {\tt octcomp}, providing some basic octave compatibility 
(see \cite{pde2phome}), 
and one \mlab\ file {\tt setpde2path.m}, which is a utility function 
to set the \mlab\ path. 
Each of the demos comes with a file {\tt *cmds.m}, which contains the 
commands to run the example (and some comments), and which should be seen 
as a quarry for typical commands, and with a file 
{\tt *demo.m}, which produces more verbose output. 
To start we recommend to run {\tt setpde2path} (without arguments) 
in the root directory \pdep\ and then change into one of the 
demo--directories, e.g., type {\tt cd demos/bratu} in \mlab. 
Then inspect the file {\tt bratucmds.m} and copy paste the commands 
to the \mlab\ command line, or just execute {\tt bratucmds} or 
{\tt bratudemo}.

\subsubsection{General structure, initialization, 
and continuation runs }
In \pdep, a continuation and bifurcation problem is described by a 
structure, henceforth called {\tt p} (as in problem), which we now outline, see also 
Tables  
\ref{tab2} and \ref{tab2b}.\footnote{In addition to the fields/variables 
listed, there are quite a few more within {\tt p}. See 
{\tt stanparam.m} for these, and also for more comments on the ones which 
are listed. Some of the ``control fields'' are unlikely to be changed 
by the user, at least at the beginning, e.g., {\tt p.evopts.disp=0} (to suppress 
output during eigenvalue calculations), or 
{\tt p.pfig=1, p.brfig=2} (the figure numbers for plotting),  
and some of the additional 
fields/variables are only generated during computation, e.g., 
the residual {\tt res} and the error estimate {\tt err},  
which are put into {\tt p} for easy passing between subroutines and 
user access. Currently there are no global variables in \pdep, 
with the exceptions {\tt pj,lamj} which are set for 
numerical differentiation in {\tt resinj}, and possibly LU preconditioners 
for iterative linear system solvers, see \S\ref{lss-sec}. See also \S\ref{custsec}.}  Essentially, 
{\tt p} contains 
\bi 
\item function handles which describe the functions $c,a,b,f$ and the \bcs\ 
(and possibly the Jacobian) in \reff{gform};
\item fields which describe the geometry of the problem, 
including the FEM mesh;
\item fields which hold the current solution, i.e., $u,\lam$ and the tangent 
$\tau$;
\item a number of variables, switches and further functions (i.e., function handles) controlling the behaviour of the continuation and bifurcation algorithm, and filenames for file output. 
\ei 

\begin{table}[!ht]
\bce
{\footnotesize
\begin{tabular}{|p{23mm}|p{88mm}|p{42mm}|} 
\hline
func.\,handles&meaning&example (in {\tt initbratu})\\
\hline
f&{\tt [c,a,f,b]=f(p,u,lam)} PDE coefficients in \reff{gform}  &p.f=@bratuf\\
jac& {\tt [c,fu,flam,b]=jac(p,u,lam)} used to build 
$G_u, G_\lam$ from $f_u, f_\lam$ if desired 
(see jsw below), &p.jac=@bratujac\\
bcf&{\tt bc=bcf(p,u,lam)}, 
boundary conditions function&p.bcf=@(p,u,lam) bc$^{a,c}$\\
outfu&{\tt out=outfu(p,u,lam)}, defining output for bifurcation diagram, 
i.e.,~quantities saved on {\tt p.branch}. Default setting 
out=stanbra(p,u,lam) gives out=[$\|u_1\|_\infty; \|u_1\|_{L^2}$];
&p.outfu=@stanbra$^b$\\
ufu&{\tt cstop=ufu(p,brout,ds)}, 
user function, called after each calculation of a point, e.g.~for 
printing information and checking stopping criteria
&p.ufu=@stanufu$^b$\\
headfu&{\tt headfu(p)} (no return arguments); 
print headline of screen output&p.headfu=@stanheadfu$^b$\\
blss, lss& (bordered) linear system solver, see \S\ref{lss-sec}&
p.blss=@blss$^b$, p.lss=@lss$^b$.\\
\hline\hline
various var.&meaning&example (in {\tt initbratu})\\\hline 
geo&geometry (and also BC in call to recnbc1)  
&[p.geo,bc]=recnbc1(0.5,0.5);$^{a,b,c}$\\
points,edges,tria&mesh&p=stanmesh(p,30,30)\\
bpoints,..,btria&{\em b}ackground mesh (used for mesh-adaption) 
&p=setbmesh(p)$^b$\\
neq, np, nt& number of equations, mesh points, triangles&automatically from mesh\\
u,lam,tau,ds&essential continuation data&p.u=0.2*ones(p.np); \ldots;\\
branch, bifvals& lists generated via p.outfu, used to 
plot bif.\,diagrams&p.branch=[]; p.bifvals=[];\\
\hline\hline
file names&meaning& standard (in {\tt setfn})\\\hline 
pre& name of subdir for files; by default automatically set to the struct name used in call to {\tt cont} or {\tt *init} & \\
\mbox{pname,bpname}&(base)filenames for output of points, 
bifurcation points; actual filenames  augmented by counter&
p.pname='p' for {\em p}oint, p.bpname='bp' for {\em b}if.{\em p}oint\\
\hline
\end{tabular} }
\ece
\footnotesize{
$^a$ typical setup with {\tt bc} independent of $\lam, u$ and hence defined in advance\\
$^b$ ``standard'' choices provided by \pdep\\
$^c$ see  {\tt bratuinit.m} and documentation of \ptool\ for explanation
}
\caption{{Basic variables in structure {\tt p} of a problem, 
and their initialization in {\tt bratuinit}.
}\label{tab2}}
\end{table}
\begin{table}[!ht]
\bce{\footnotesize
\begin{tabular}{|p{44mm}|p{113mm}|} 
\hline
main functions&purpose\\
\hline
p=cont(p)&main continuation routine\\
p=swibra(pre,file,npre)&branch-switching at bifurcation point 
from previous run, prefix set to npre \\
p=meshadac(p,varargin)&adapt mesh in p, see \S\ref{numremsec} for 
varargin\\
p=findbif(p,ichange)&locate bifurcation point based on $G_u$; cont itself uses {\tt bifdetec} (bifchecksw$>0$)\\ 
p=loadp(pre,file,npre)&load solution struct p from file and reset prefix to npre\\
plotbra(p,wnr,cmp,aux)&plot component cmp of branch over $\lam$, 
in figure number wnr \\ 
plotbraf(pre,file,wnr,cmp,aux) &as plotbra but 
from file (saved previously, leave out the '.mat'),\\
plotsol(p,wnr,cmp,pstyle)
&plot solution, use plotsolf(pre,file,wnr,\ldots) to plot from file\\
\hline
\end{tabular}}
\ece

\vs{-2mm}
\caption{{Main ``user'' functions in \pdep.}\label{tab3}}
\end{table}
\begin{table}[!ht]
\bce
{\footnotesize
\begin{tabular}{|p{35mm}|p{122mm}|}
\hline
Newton\&Cont&meaning \\\hline 
imax,normsw&Newton controls: max number of steps and selection of a norm (`norm-switch')\\
tol&stop-crit.~for Newton, typically should be around 1e-10; \\
nsw&0 for Newton, 1 for chord\\
jsw&switch for derivatives ($G_u, G_\lam$). \mbox{ 0: ($G_u, G_\lam$) by 
$(c,f_u,b,f_\lam)$,}\\
& \mbox{ 1: $G_u$ by $c,f_u,b$,\quad $G_\lam$ by FD,\qquad}
2: $G_u$ by FD, $G_\lam$ by
$f_\lam$, \quad 3: both by FD.\\
dsmin,dsmax,dlammax&min/max stepsize in s, max stepsize in $\lam$\\
lammin,lammax&min/max $\lam$, preset to $\mp 10^6$, reset to use as stopping criteria\\
nsteps&number of steps to take\\
parasw,lamdtol&parametrization switch and tolerance: if parasw=0 
resp.~parasw=2 then always use \reff{trn} resp.~\reff{newton}. If 
parasw=1 then use \reff{csw} with $\lamdtol=$lamdtol\\
amod, maxt, ngen& controls for mesh adaption: adapt every amod-th step, 
aim at maxt triangles, in at most 
ngen refinement steps\\
errchecksw, errtol&switch and tol for a posteriori error estimate and handling, 
see \S\ref{numremsec}. \\
\hline 
Bif.,Plot\&User-control&meaning\\[1mm]\hline
bifchecksw&0 for no checks, 1 for check via B1 with consistency 
with eigenvalues of $G_u$, \mbox{2 for B1 alone}\\
neig&number of eigenvalues to be calculated in {\tt spcalc} and {\tt bifdetec}, 
default 50\\
eigsstart&0 to start eigs randomly, 1 to start with $(1,\ldots,1)$ (default) \\
bisecmax&max number of bisections to locate bifurcation; turned off 
by bisecmax=0\\
spcalcsw & 0/1 for eigenvalue calculation off/on\\
pmod/pstyle&plot each pmod-th step in style pstyle (1 mesh, 2 pcolor, 3 rendered 3D, \ldots)\\
pcmp, bpcmp&component to plot, component of branches to plot\\
smod&save every smod-th step\\
isw,vsw&interaction/verbosity switch: 0=none, 1=some, 2=much; \\
pfig,brfig,ifig&figure-numbers for u-plot, branch-plot and info-plot 
during cont.; in ifig we plot add.~information, e.g., 
after mesh adaption, or the new tangent after {\tt swibra}\\
timesw&if $>0$, print timing info after {\tt cont}. See 
{\tt stanparam.m} for the timers in {\tt p}. \\
nbp&number of user--components of branch to be printed on screen 
(p.ufu=@stanufu) \\
\hline
resfac, mst, pmimax& {\tt pmcont} only, see \S\ref{pmcont-sec}\\
\hline 
\end{tabular}}
\ece
\caption{{Main switches and controls in a structure {\tt p} used in {\tt cont}, 
see {\tt stanparam.m} for typical values and a number of additional 
switches with detailed comments. 
See also \S\ref{pmcont-sec} for additional parameters controlling 
{\tt pmcont}, the {\tt p}arallel{\tt m}ulti--continuation version.}\label{tab2b}}
\end{table}

\noi Studying a continuation and bifurcation problem using \pdep\ thus 
consist of:  
\bi
\item Setting up a file defining the coefficient functions 
$c$, $a$, $b$ and $f$ (and usually a function for the Jacobians) 
in \reff{gform}, e.g. {\tt bratuf.m} (and {\tt bratujac.m}). 
(Here we assume that the \bcs\ function p.bcf is defined inline as in {\tt bratuinit.m} 
and many of the further examples) 
\item Setting up an initialization function 
file, e.g., {\tt bratuinit} filling {\tt p}. The main steps are (1) define 
p.f and p.jac, (2) define geometry and mesh, and usually the \bcs\ by an inline function, 
(3) set the parameters and provide 
a starting point. In many cases most parameters and switches can be set to 
``standard values''. For this we provide the function {\tt p=stanparam(p)}, 
which should be called first, and afterwards individual parameters 
can be reset as needed. For the mesh generation we provide an elementary 
function {\tt p=stanmesh(p,hmax)}, where {\tt hmax} is the maximal 
triangle side--length. This is based on {\tt initmesh} from the {\tt pdetoolbox}, i.e., 
a Delaunay algorithm. 
For rectangular domains the syntax  
{\tt p=stanmesh(p,nx,ny)} is also allowed, which is based on 
{\tt poimesh} from the {\tt pdetoolbox}, with obvious meaning. 
\footnote{If applicable, {\tt poimesh} obviously is faster and gives 
more regular meshes, if used with care, i.e., choose {\tt nx/ny} according 
to $L_x/L_y$ where $L_x,L_y$ are the sidelengths of the rectangle. 
However, we found 
that in some cases the Delaunay mesh gives more robust numerics, 
in particular after mesh refinement (see \S\ref{numremsec}). Thus we 
recommend to experiment with both.}
\item Calling a number of \pdep\ functions. The basic call is 
{\tt p=cont(p);} (a continuation run), which can be followed, e.g.,  
by a repeated call to extend the branch. 
Or, in case a branch point has been found, a call to branch switching 
and subsequent continuation by {\tt q=swibra('p','bp1','q'); q=cont(q);} 
where, e.g., {\tt ./p/bp1.mat} is data written at a branch point during the previous run. 
Inbetween runs, {\tt p} can be modified from the command line, 
e.g., type {\tt p.imax=5} (say) to (re)set the maximal number of 
Newton-iterations, or call {\tt p=meshref(p)} to refine the mesh 
before a subsequent run; afterwards, call {\tt p=cont(p)} again. 
According to the settings, data is plotted and written to files in a 
sub-directory with name {\tt p.pre}. There are also
functions for further postprocessing, e.g., plotting of solutions 
and bifurcation diagrams, whose documentation is mainly provided within 
the corresponding matlab files, and by the calls in the example directories. 
\ei 

\begin{figure}[ht!]
\bce 
\begin{tabular}{llll}
(a)&(b)&(c)&(d)\\
\ig[width=36mm]{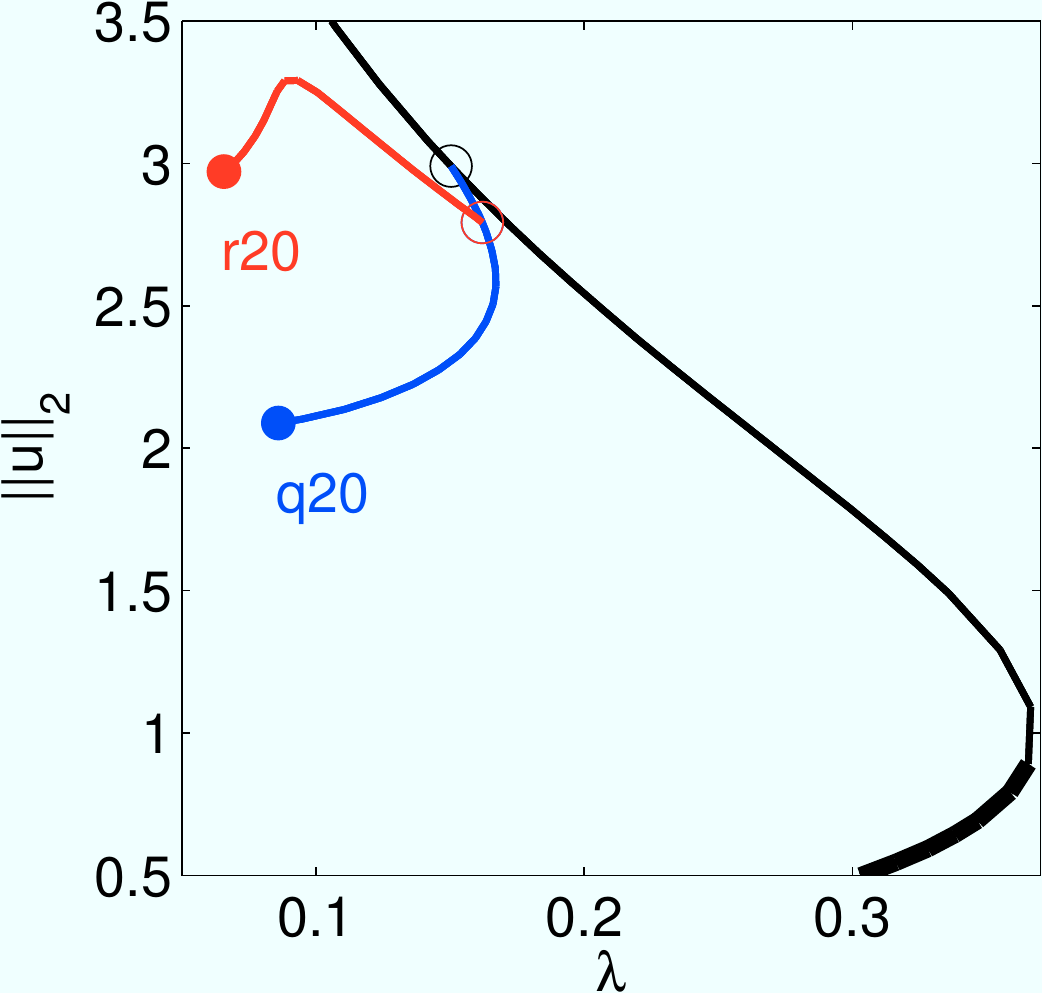}
&\ig[width=38mm]{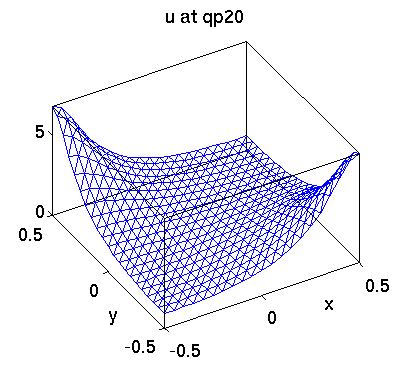}
&\ig[width=38mm]{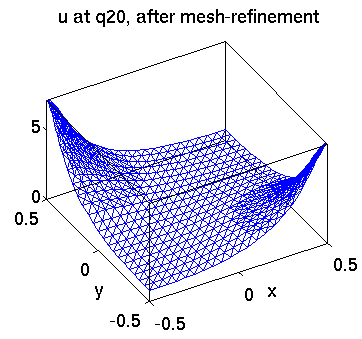}
&\ig[width=38mm]{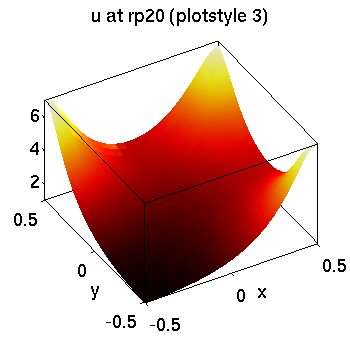}
\end{tabular}
\ece
\caption{{\small (a) Elementary bifurcation diagram for 
\reff{brprob} ($\|u\|_{L^2}$ over $\lam$) generated by \pdep, 
over a uniform mesh with 800 triangles. Thick lines 
indicate stable (parts of) branches, thin lines unstable branches, 
$\circ$ bifurcation points. (b), (d) Some solution plots. (c) Preview 
of mesh refinement. 
 See \S\ref{numremsec} 
for the quality of the mesh at the ``ends'' of branches $q,r$ 
and the due mesh--refinement. 
\label{bratubdf}}}
\end{figure}

The fundamental user provided functions thus are the coefficient 
function {\tt p.f} and the additionally recommended jacobian function 
{\tt p.jac}, see \S\ref{coeffsec}. 
To study a new problem, we  recommend to edit copies of the files {\tt *init.m},  
 {\tt *f.m} and {\tt *jac.m}  
of a suitable example (e.g.~{\tt *=bratu} for a scalar problem) in an empty directory, and start with calling 
{\tt p=[];p=newinit(p);} {\tt p=cont(p)}. 
The bifurcation diagram in Fig.~\ref{bratubdf} and the 
solution plots are generated from 
the commands in Table \ref{mtab2}, 
either given from the command line, 
or put into a \mlab\ script 
\footnote{The figures generated from the commands in {\tt bratucmds.m} (and similarly 
for the remaining demos to come) may differ slightly from the ones in this manual; 
this is due to 
minimal postprocessing via \mlab's ``Edit Figure Properties'' to set fontsizes, 
axis labels, etc.}. 

\begin{table}[!ht]
\bce{\footnotesize
\begin{tabular}{|p{160mm}|}
\hline\\\vs{-10mm}
\begin{verbatim}
function p=bratuinit(p) % init-routine, see bratuinit.m for more comments
p=stanparam(p); p.neq=1; p.f=@bratuf; p.jac=@bratujac; [p.geo,bc]=recnbc1(0.5,0.5); 
p.bcf=@(p,u,lam) bc; % typical inline definition of the BC function 
nx=20;p=stanmesh(p,nx,nx); p=setbmesh(p); % mesh and "background" mesh 
pre=sprintf('%s',inputname(1)); p=setfn(p,pre); % set filename (prefix)
p.xi=1/p.np; p.dlammax=0.02; p.lammin=0.02;p.tau=1; % set p.tau to something
p.lam=0.2; p.u=0.1*ones(p.np,1); p.ds=0.05; % "trivial" branch 
\end{verbatim}
\vs{-8mm}
\\\hline\hline\\\vs{-10mm}
\begin{verbatim}
function [c,a,f,b]=bratuf(p,u,lam) %% coeff for Bratu 
u=pdeintrp(p.points,p.tria,u); c=1; a=0; f=-10*(u-lam*exp(u)); b=0; 
\end{verbatim}\vs{-8mm}\\\hline\hline  \vs{-5mm}
\begin{verbatim}
function [c,a,f,b]=bratujac(p,u,lam) %% Jacobian for Bratu 
u=pdeintrp(p.points,p.tria,u); c=1; fu=-10*(1-lam*exp(u)); flam=10*exp(u); b=0;
\end{verbatim}
\vs{-8mm}\\\hline\hline 
\vs{-4mm}
\begin{verbatim}
% commands to run bratu in pde2path (selection, see also script file bratucmds.m)
p=[];p=bratuinit(p); p=cont(p); 
q=swibra('p','bp1','q'); q.lammin=0.1;q.nsteps=20;q=cont(q); 
plotbra(p,3,2,'ms',12,'lw',5,'fs',16,'cl','k');plotsolf('q','p20',4,1,1);
\end{verbatim}
\vs{-8mm}\\\hline
\end{tabular}}
\ece

\vs{-2mm}
\caption{{\small The basic init--routine {\tt bratuinit.m}, the definitions 
of  PDE coefficients and Jacobian (see \S\ref{coeffsec}), 
and some selected commands (see {\tt bratucmds.m}) to run \pdep.  
See also the files for detailed comments. 
}\label{mtab2}}\end{table}

\subsubsection{The PDE--coefficients and Jacobians}
\label{coeffsec}
The coefficient function {\tt p.f} is fundamental, and the jacobian function 
{\tt p.jac} is recommended. 
The input argument {\tt u} of both is the vector of nodal values, 
with $u_1(\cdot)=${\tt u(1:p.np)}, $u_2(\cdot)=${\tt u(p.np+1:2*p.np)}, \ldots, 
 $u_N(\cdot)=${\tt u((p.neq-1)*p.np+1:p.neq*p.np).} 
For the outputs $c,a,f,b$ (of {\tt p.f}) we allow two 
forms, i.e., (arrays of) constants, or (arrays of) 
values {\em on the triangle midpoints} of the FEM--mesh, 
essentially as explained in \cite{pdetdoc}. 
There are two major ways to generate $c,a,f,b$ from $u$. We first focus on $f$: 
\bi
\item[a)] Use  {\tt u=pdeintrp(p.points,p.tria,u)} to first interpolate $u$ to the 
triangle values (again called 
$u$), which yields a matrix 
$$u=\bpm u_{11}&u_{12}&\ldots&u_{1,n_t}\\ 
\ldots&\ldots&\ldots&\ldots\\
u_{N1}&u_{N2}&\ldots&u_{Nn_t}\epm, 
$$
where $n_t$ is the number of triangles in the mesh. 
Then write, e.g., $f(u)$ in a standard 
\mlab\ way, i.e., {\tt f=-10*(u-lam*exp(u));} see {\tt bratuf} in Table 
\ref{mtab2}. 

\item[b)] First express, e.g., $c,f$  as 
``Matlab text expression in x,y,u,ux,uy'' (from \cite{pdetdoc}), 
with obvious meaning. For this, a parameter {\tt lam} must be converted 
to a string.  Afterwards, {\tt pdetxpd} is called to 
evaluate the text expression on the triangle midpoints. 
See file {\tt bratuft.m} for an example. 
\ei 
Option a) has the advantage that it is more ``natural'', more flexible, and, 
at least for 
simple expressions, slightly shorter. If, however, $f$ depends on 
$x,u_x,\ldots$, then option b) might be shorter. 
To some extent it is a matter of taste which way to generate 
$f$ is preferred (and similarly $c$,\ldots), therefore for \reff{brprob} 
we provide both as templates. However, b) only allows local 
dependence $f(u,\ldots)=f(u(x),\ldots)$, while in a) we have more flexibility 
and, for instance, can also call external functions in a simple way. 
\footnote{See, e.g., \S\ref{acgc-sec}.}

For {\tt c}, which in principle is the $N\times N\times 2\times 2=1{\times}1{\times}2{\times}2$ tensor 
$c_{11kl}=\bpm 1&0\\ 0&1\epm$ we simply write {\tt c=1} which corresponds 
to the simplest symmetric case. \footnote{There are various 
special coding schemes 
for diagonal or symmetric cases, see \cite{pdetdoc}. For convenience 
here we just note that in the general case {\tt c} is a $4N^2\times {\tt p.nt}$ 
matrix (or in case of constant coefficients a $4N^2$ column vector) 
with $c_{ijkl}$ in row $4N(j-1)+4i+2l+k-6$, $i,j=1,\ldots,N$, $k,l=1,2$. In this 
scheme we would that have {\tt c=[1;0;0;1]} to encode the Laplacian. 
Similarly, there are special storage schemes for symmetric $a$, 
but in general $a=a_{ij}$, $i,j=1,\ldots,N$, is stored as a 
$N^2\times {\tt p.nt}$ matrix (resp.~a $N^2$ column vector in case of 
constant coefficients) with $a_{ij}$ in row $N(j-1)+i$. Note the somewhat 
non--lexicographical order.}

The tensor $b=b_{ijk}$, $i,j=1,\ldots,N$, $k=1,2$ in \reff{gform} 
is not part of the \ptool. 
Its coding and storage mimics that of {\tt c}. In detail, 
{\tt b} is an $2N^2\times m$ array, 
where $m=1$ (constant case) or $m=${\tt p.nt}, in the order
\huga{\label{bform}{\small
b=\bigl[b_{111}; b_{112}; b_{211}; b_{212}; \ldots b_{N11}; b_{N12};\quad  
b_{121}; b_{122}; \ldots; b_{N21}; b_{N22}; \ \ldots\ 
b_{1N1};\ldots;b_{NN2}\bigr],}
}
i.e., $b_{ijk}$ is in row $2N(j-1)+2i+k-2$. 
Unlike the \ptool\ setup for $c,a$ there are currently no 
schemes to encode special situations such as, e.g., 
$b\otimes \nabla u=\al\pa_x u$ which corresponds to advection 
into direction $x$. Thus, in this case, for $N=2$, {\tt b} reads 
$b=[\al;0;0;0;\,\,0;0;\al;0]$, while, e.g., $\beta\pa_y u$ 
yields $b=[0;\beta;0;0;\,\,0;0;0;\beta]$. 
Again we remark that advective terms can also be put into $f$ such that  
{\tt p.f} may simply return {\tt b=0}. 
However, {\tt b} is needed if there are advective terms and we 
want explicit Jacobians as explained next.

If {\tt p.jsw}$<3$ then the user must also set {\tt p.jac} to a function 
handle with outputs ${\tt c,fu,flam,b}$ again defined on triangle midpoints, 
where in fact ${\tt c}$ is the same as in {\tt p.f}.  
Compare Table \ref{mtab2}, and Remark \ref{jacrem} -- also
concerning the meaning of {\tt fu} and ${\tt b}$. Providing the same 
{\tt c} (and also often the same ${\tt b}$) twice (in {\tt p.f} and {\tt p.jac}) 
is redundant, 
and there are situations where only either {\tt c,fu,b} or {\tt flam} 
are needed, namely {\tt jsw=1} resp.~{\tt jsw=2}. However, 
we found the small overhead of recomputing {\tt c,b}, resp.\,unnessessarily 
computing ${\tt fu,b}$ resp.\,${\tt flam}$ acceptable to have a clear code. 
On the other hand, splitting the calculation of PDE coefficients and 
Jacobian coefficients into two routines is reasonable since often 
at least for testing it is convenient to use {\tt jsw=3} where 
analytical jacobians are never needed.

\brem\label{jacrem} 
Given coefficients $c,a,f,b=c(u_0),a(u_0),f(u_0),b(u_0)$, 
the FEM transforms \reff{gform} 
into the algebraic system ${\tt K(u_0)u-F(u_0)}\stackrel{!}=0$ 
where {\tt K} is called the stiffness matrix, assembled from $c,a$ and $b$, and 
$\Ftt$ is the FEM representation of $f$. 
Thus, {\tt u} solves the FEM discretization of \reff{gform} if the residual  
${\tt r=resi(p,u,lam):=K(u)u-F(u)\stackrel{!}=0}$.  
The basic \ptool\ routine to assemble {\tt K=K}$_0$ (in case $b=0$) and {\tt F} 
is {\tt [K,F]=assempde(\ldots,c,a,f)}, where ``\ldots'' stands for boundary 
conditions and mesh-data. To assemble 
the advection matrix $B$ we additionally provide 
{\tt B=assemadv(p,t,b)}. The full system--matrix then 
is $K=K_0-B$ . 

Thus, if $a=b=0$ and 
$c$ does not depend on $(u,\lam)$,  then with 
the local derivatives ${\tt fu}=f_u$ and ${\tt flam}=f_\lam$ 
returned from {\tt p.jac}, the 
Jacobian $G_u$ and the derivative $G_\lam$ can be obtained 
from 
\huga{\label{jacass}{\tt [Gu,Glam]=assempde(\ldots,c,-fu,-flam)}, }
and this is done for {\tt jsw=0}. 
If $a\ne 0$ is independent of $u$ (and still $b=0$), then {\tt p.jac} 
must return 
$$
{\tt fu}=f_u-a.
$$
If $b=0$ but, e.g., $a$ depends on $u$, then the formulas must be adapted 
accordingly. In other words, 
{\em in calculating {\tt fu} assume that $a=0$ in \reff{gform}, i.e., 
if $a\ne 0$ in {\tt p.f}, 
then define $\tilde{f}(u)=f(u)-au$ and set ${\tt fu}=\tilde{f}_u$ in 
{\tt p.jac}}. If $c=c(u)$ or $b=b(u)\ne 0$, 
then $G_u$ can still be assembled using $c(u)$ and suitable 
$\tt b$ and ${\tt fu}$, see \S\ref{acqlsec} and \S\ref{chem-sec} for examples. 
Similarly, {\tt -flam} must always be understood in a ``generalized'' sense, 
i.e., it must assemble to $G_\lam$, even if this involves high 
derivatives of $u$ which originally were implemented via $c$. 
In any case, remember {\em that the notations {\tt fu}, {\tt flam} and 
${\tt b}$ 
in {\tt p.jac} are only conventions}, motivated by the fact that 
the case that only $f$ depends on $(u,\lam)$ is the most common. 

For {\tt jsw=1}, {\tt Gu} is still assembled but {\tt Glam} is 
calculated by finite differences.\footnote{Since 
approximating $G_\lam$ by finite differences only takes one additional call 
of {\tt resi}, speed is not an issue for choosing between {\tt jsw=0} and 
{\tt jsw=1}. 
In the latter case, simply set {\tt flam=0} in {\tt p.jac}. 
See \S\ref{para-sw-sec} for an example.} 
For ${\tt jsw}=2$ we use {\tt numjac} to calculate {\tt Gu}; 
to be efficient 
this requires a sparsity structure {\tt S}, and here we assume that 
$\Ftt_i$ depends only on (all components of) $u$ on the i-th node and all 
neighboring nodes, which corresponds to the sparsity structure obtained 
by {\tt [Gu,Glam]=assempde(\ldots,0,a,0)} with 
$a_{ij}=1$ for $i,j=1,\ldots,N$. Our experience is that 
numerical Jacobians are fast enough for moderate size problems, 
i.e., for up to a few thousand degrees of freedom. 
Of course this also depends on the structure 
of the problem: diagonal diffusion or not, weak or strong 
coupling of the different components of $u$. 
Still, assembling Jacobians is usually much faster. 
For ${\tt jsw}=3$, both {\tt Gu} and {\tt Glam} are approximated by finite differences. 
In any case, for both (${\tt jsw}\le 1$ and ${\tt jsw}\ge 2$) 
we assume local dependence of $f$ on $u$. 
See, however, \S\ref{acgc-sec} for some modifications for the case of global coupling. 
\eex 
\erem 

\brem\label{bcrem} 
 The boundary conditions, see \S\ref{geosec}, are 
updated from {\tt bc=p.bcf(p,u,lam)} before assembling. 
In the (frequent) case that the \bcs\  
do not change during continuation we set may 
{\tt p.bcf=@(p,u,lam) bc} in the init-routine (after generating {\tt bc}). 
See, however, \S\ref{s:acbc} for examples 
with $\lam$--dependent \bcs.
\eex\erem 

As mentioned, when applicable, assembling $G_u$ via $c,a,f_u$ and $b$ ({\tt p.jsw=0,1}) gives a matrix $G_{u,a}$ and is faster 
(by orders of magnitude for large $Nn_p$)
than numerical differentiation by {\tt numjac} ({\tt p.jsw=2,3}), which 
gives a matrix $G_{u,n}$ which is in general close to but not equal to $G_{u,a}$. 
Intuitively we might also expect $G_{u,a}$ to be ``more 
accurate'' than $G_{u,n}$. However, we need some caution: in fact, 
$G_{u,n}$ often gives better convergence of the Newton loop 
for the algebraic system $r(u)=K(u)u-F(u)\stackrel{!}{=}0$.  The reason 
is that $G_{u,a}$ involves 
interpolation of nodal values to triangle values in $c,a,b$ and $f_u$, while for 
$G_{u,n}$ this is done on $c,a,f,b$, consistent with the definition 
of $r(u)$. 
This effect becomes prominent on poor (underresolved) 
meshes, where the relative 
error $e=\|G_{u,n}-G_{u,a}\|/\|G_{u,n}\|$ can be of order $0.05$ or larger.  
However, $e\ra 0$ for mesh spacing $h\ra 0$. 
For convenience we provide the function {\tt  [Gua,Gun]=jaccheck(p)} 
which returns $G_{u,a}, G_{u,n}$, produces spy--plots of these 
matrices, and prints the timing and some diagnostics. 

Thus, for small $n_p$ it might appear that $G_{u,n}$ is favorable. 
On the other hand, $G_{u,a}$ is obtained  in 
acceptable time on much finer meshes, where the FEM solution $u_h$ 
should be much closer to a PDE solution $u$. In fact, using 
{\tt jsw=2,3} can even be dangerous in the sense that 
it may mask the fact that a FEM solution $u_h$ is not close to a PDE solution $u$. See \S\ref{numremsec}.

\subsubsection{The geometry and the boundary conditions}\label{geosec}
The domain $\Om$ is typically described as a polygon. As the \ptool\ syntax 
is somewhat unhandy, and the rectangular case is quite common we provide 
the function {\tt geo=rec(lx,ly)} which yields $\Om=[-l_x,l_x]\times [-l_y,l_y]$. 
An extension is the function {\tt polygong}, see \cite{pruefert} for 
its syntax. 
Setting up ``arbitrary complicated'' geometries $\Om$ 
is most convenient if there is a drawing {\tt img.jpg} of $\Om$ 
in the current directory. 
Type {\tt im=imread('img.jpg'); figure(1);image(im); [x,y]=ginput;} 
which yields a crosshair. Click (counterclockwise) on $\pa\Om$, stop 
with {\tt return}. The obtained vectors {\tt x,y} can be saved as a {\tt *.mat} 
or {\tt *.txt} file, and piped through {\tt geo=polygong(x,y)}. 
Finally, geometries can also be exported from the \ptool\ GUI. 

The \ptool\ syntax for the boundary conditions \reff{e:gnbc} is also somewhat 
unhandy. For scalar equations the most common \bcs\ are homogeneous Dirichlet 
or Neumann \bcs.  The routine {\tt [geo,bc]=recdbc1(lx,ly,qs)} approximates Dirichlet \bcs\ over rectangles via Robin \bcs\ of the form ${\bf n}\cdot(c\otimes \nabla u)+q_s hu=0$ with a large $q_s={\tt qs}$. 

For {\tt c} of order 1 typically 
$q_s=\CO(10^2)$ or $q_s=\CO(10^3)$ works well. 
For homogeneous Neumann \bcs\ we provide {\tt [geo,bc]=recnbc1(lx,ly)}, 
with the extension {\tt [geo,bc]=recnbc2(lx,ly)} to 2--component 
systems. 

For the genuine systems case (or the case of non--rectangular domains) 
we provide the routines {\tt bc=gnbc(pneq,varargin)} 
and {\tt bc=gnbcs(pneq,varargin)}.  For a system with {\tt neq}
components and a domain with {\tt nedges} edges, 
{\tt bc=gnbc(neq,nedges,q,g)} creates ``generalized
Neumann \bcs'' \eqref{e:gnbc} that are given as numerical
data. Different boundary conditions $q_j,g_j$ at the edge with index
$j$ are generated by a call of the form {\tt
  bc=gnbc(neq,q$_1$,g$_1$,\ldots,
q$_{\tt nedges}$,g$_{\tt
    nedges}$)}. For $g$, $q$ given in terms of an explicit formula,
e.g., involving $x, u, \nabla u$, the function {\tt gnbcs} accepts a
string variable encoding of {\tt g} and {\tt q} or {\tt g$_j$}, {\tt
  q$_j$} and otherwise works in the same way. 
\subsubsection{Error estimates and mesh adaption}
\label{numremsec} 

\begin{figure}[ht!]
\bce 
\begin{tabular}{lll}
(a)&(b)&(c)\\
\ig[width=38mm]{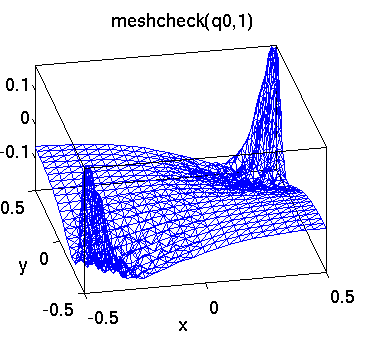}
&\ig[width=38mm, height=37mm]{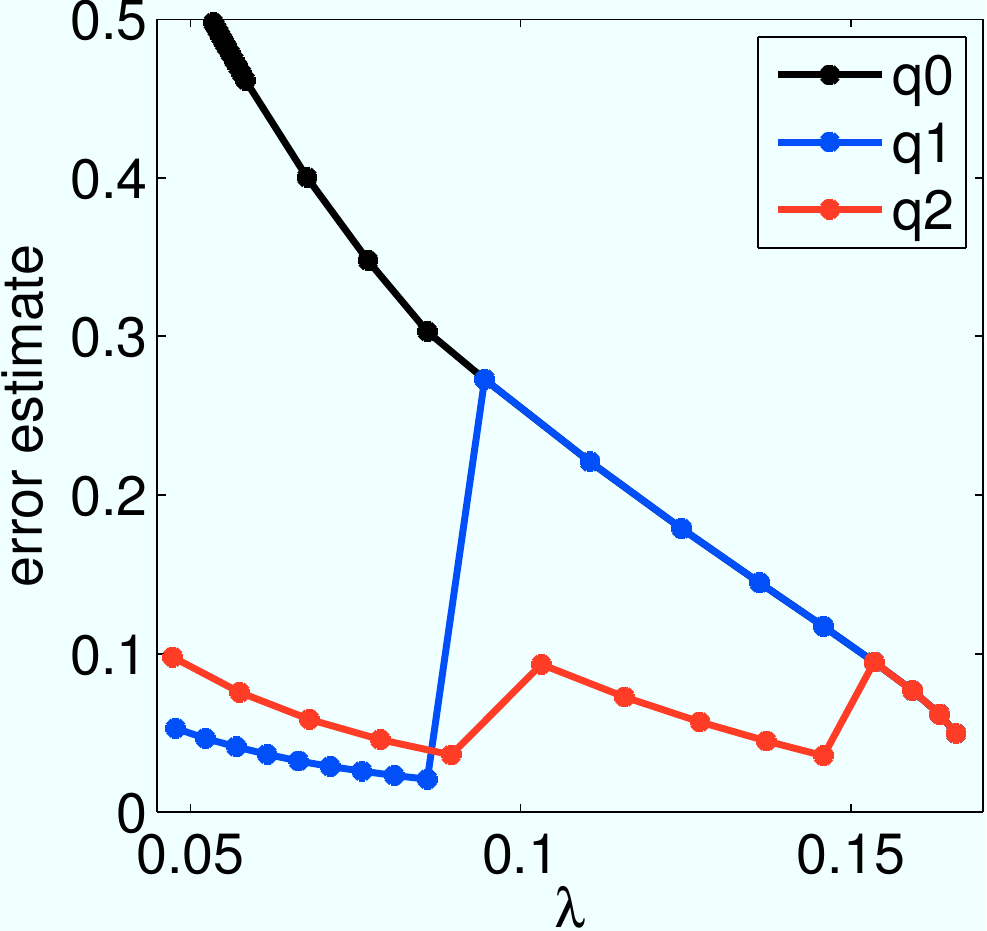}
&\ig[width=38mm]{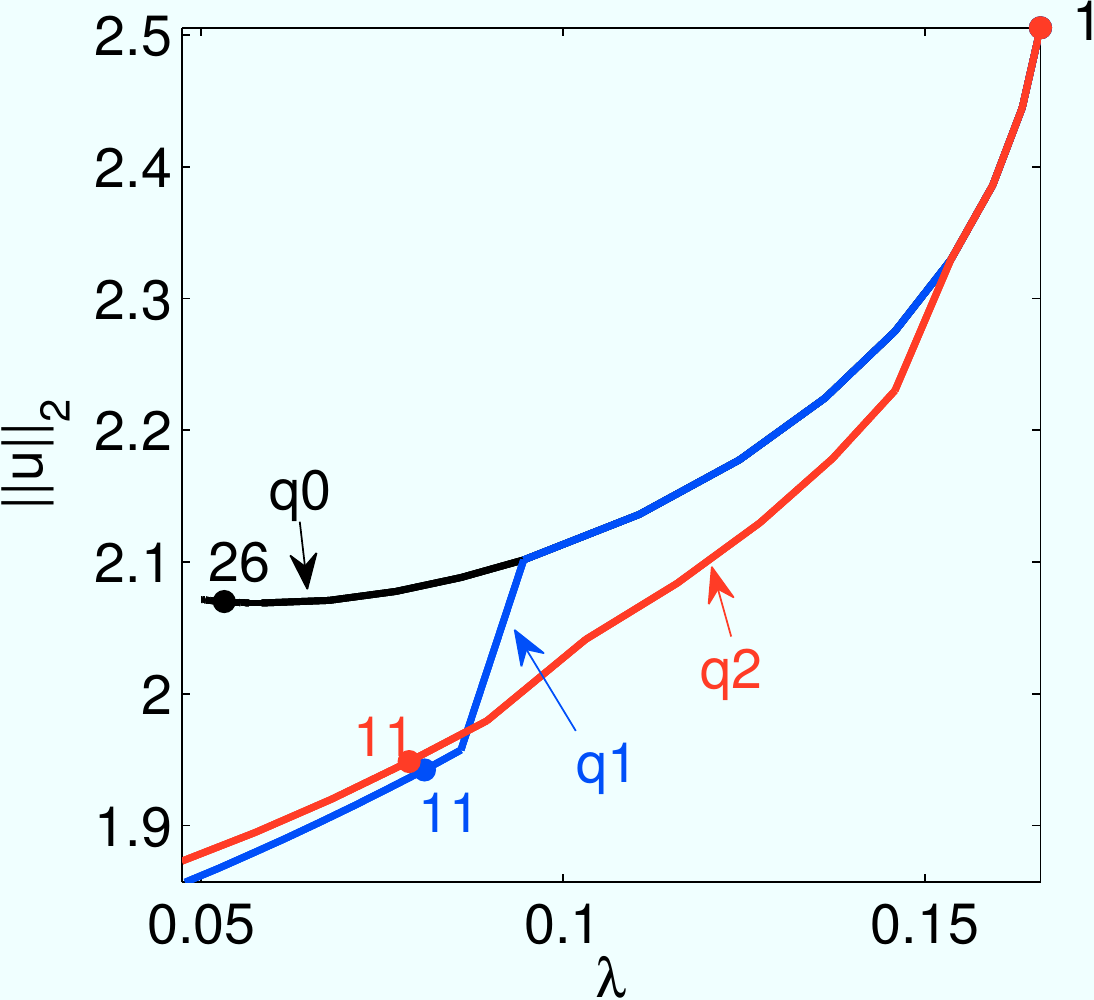}
\end{tabular}
\ece
\caption{{\small (a) Error in u at {\tt q0}=point 20 on the q branch from Fig.~\ref{bratubdf} 
obtained from calling  {\tt meshcheck(q0,1)}. (b) 
{\tt p.err} over $\lam$ for continuing from point 10 on the q branch, without 
mesh-adaption (labeled {\tt q0}), with strategy (i) ({\tt amod=10}, labeled 
{\tt q1}) and with 
strategy (ii) ({\tt p.errbound=0.1}, labeled {\tt q2}).  
(c) The bifurcation diagrams belonging to (b). 
In {\tt q1} a rather large jump appears from 
refinement after the 10$^{th}$ step.
\label{bratunf}}}
\end{figure}

As an ad hoc way to check whether a FEM solution $u_h=${\tt p.u} approximates a 
PDE solution $u$ we provide {\tt [q,ud]=meshcheck(p,cmp)}. This (adaptively) 
refines the FEM mesh 
in {\tt p} to roughly the double number of triangles and calculates 
a new FEM solution $u_{h,{\rm new}}$ from the old solution $u_{h,{\rm old}}$. Then 
$u_{h,{\rm old}}$ is interpolated to $\tilde{u}_{h,{\rm new}}$ on the new mesh, 
$u_{{\rm diff}}=u_{h,{\rm new}}-\tilde{u}_{h,{\rm new}}$ is formed, 
and  $\|u_{{\rm diff}}\|_{\infty}$ 
and the relative error $\|u_{{\rm diff}}\|_{\infty}/
\|u_{h,{\rm new}}\|_{\infty}$ are printed. The new solution structure 
{\tt q} and 
the difference {\tt ud}$=u_{{\rm diff}}$ are returned, 
and for  {\tt cmp}$>0$ we additionally generate a plot 
of the {\tt cmp}$^{th}$ component of $u_{{\rm diff}}$. 
For instance, in Fig.~\ref{bratunf}a) we check the mesh of point 20 on the q 
branch in Fig.~\ref{bratubdf} which indicates 
that the mesh in 2 corners is clearly too poor, 
and before continuing for smaller $\lam$ we should refine the 
mesh. \footnote{See {\tt bratucmds.m} for the calling sequence for Fig.~\ref{bratunf}a). 
It is also often useful to repeat calls 
to {\tt meshcheck}, i.e., here call next {\tt q=meshcheck(q,1)}, 
until the relative error becomes suffiently small, indicating 
that $u_{h}$ is a good approximation.} 

Thus, some (automatic) mesh adaption may be 
vital for reliable continuation. The \ptool\ comes with 
mesh refinement based on an a posteriori error estimator as follows. 
For the scalar Poisson problem $-\Delta u=f$, $u|_{\pa\Om}=0$, 
let $u_h$ be 
the FEM solution and $u$ the PDE solution. Then, with $\al,\beta>0$ 
some constants independent of the mesh, 
\huga{
\|\nabla(u-u_h)\|_{L^2}\le \al \|hf\|+\beta D_h(u_h), 
}
where $h=h(x)$ is the local mesh size, and 
$ D_h(v)=\bigl(\sum_{\tau\in E_i}
h_\tau^2(\pa_{n_\tau} v)^2\bigr)^{1/2}$. Here $\pa_{n_\tau} v$ is the 
jump in normal derivative of $v$ over the edge $\tau$, $h_\tau$ the 
length of the edge, 
and $E_i$ the set of all interior edges. For equations 
$-\nabla(c\otimes\nabla u)+au=f$ this suggests the error indicator function 
\huga{\label{trierr}
E(K)=\al\|h(f-au)\|_K+\beta\left(\frac 1 2 \sum_{\tau\in\pa K} 
h_\tau^2(n_\tau\cdot c\nabla u_h)^2\right)^{1/2}
}
for each triangle, which is calculated by the \mlab\ routine {\tt pdejmps}. 
For convenience we provide the interface 
routine {\tt err=errcheck(p)}.\footnote{\label{bfoot} 
With {\tt a,f,b} returned from {\tt p.f}, this also re--calculates $f$ (via 
{\tt fb=bgradu2f(p,f,b,u)}) to include $b\otimes\nabla u$ into {\tt fb}, as {\tt pdejmps} 
does not take  $b\otimes \nabla u$ into account directly. Similarly, 
the mesh refinement below always recalculates $f$ to include $b\otimes \nabla u$.  
Moreover {\tt errcheck} also contains our 
settings for some tunable parameters of {\tt pdejmps}.} Calling, e.g., 
{\tt err=errcheck(p)} yields {\tt err=0.273} which somewhat overestimates 
the error in Fig.~\ref{bratunf}(a). In {\tt cont}, for {\tt errchecksw>0} 
we call {\tt errcheck} after each successful step and store the result 
in {\tt p.err}. 
 
The (basic) mesh refinement 
strategy then is to introduce new triangles where $E(K)$ 
is large.\footnote{It is not a priori clear if this is 
also suitable for systems. 
Nevertheless we found {\tt pdejmps} to work well also for the problems 
considered below. 
}
This is done by the \ptool\ routine {\tt refinemesh}, but we provide the 
interface routine {\tt p=meshref(p, varargin)}.\footnote{where 
{\tt varargin} takes pairs {\tt 'maxt',maxt}=number of triangles aimed at, 
or {\tt 'ngen',ngen}=number of refinement steps, or {\tt 'eb',eb}=error bound.  
Calling for instance {\tt q0r=meshref(q0,'eb',0.0025,'maxt', 50000, 'ngen',20)} 
shows that 
to achieve an estimated error$\le 0.0025$ we need about 30.000 triangles. }
Since {\tt meshref} also interpolates the tangent $\tau$ to the new mesh, 
we can continue with {\tt cont} immediately after mesh 
refinement. 
However, instead of mesh refinement, 
which means introduction of new points into the mesh, we rather 
need mesh adaption, which means refinement where necessary, but coarsening 
where possible, to limit storage requirements. 
In \pdep, mesh--adaption is implemented in an ad hoc way in the function 
{\tt p=meshadac(p,varargin)} 
by first interpolating a given solution to a (typically somewhat coarse) 
``base--mesh'' or ``background-mesh'' and then refining. 
\footnote{base-mesh given by {\tt  p.bpoints, p.bedges, p.btria}, {\tt varargin} 
as above. The implementation of some true adaption routine is on our to-do-list.} 

During continuation runs there are basically two strategies for this, 
which can also be mixed: 
 
(i) call {\tt meshadac} every {\tt p.amod}$^\text{th}$ step, 
for {\tt p.amod>0}, or 

(ii) call {\tt meshadac} 
whenever {\tt p.err}$>${\tt p.errbound} (choose {\tt p.errchecksw>1} for 
this)\footnote{where we refine until {\tt err<p.errbound/2} in order 
to allow some margin for the next steps}. 

 See Figures \ref{bratunf}(b),(c) for a comparison of the different approaches. 

In summary, mesh adaption strategies and error bounds are 
highly problem dependent, and moreover, may not be rigorously justified 
for the system case or general \bcs. 
Thus, although it works well for a number 
of examples we considered, mesh-adaption should be applied 
with care.  Also note that mesh-adaption may lead to spurious jumps in 
det$A$. 

\brem{\rm For bifurcation from trivial branches, another good strategy 
is to prepare a {\em finer} base mesh than the starting mesh,  for 
instance if the trivial branch consists of spatially homogeneous solutions,  
but the bifurcating solutions develop sharp gradients. For convenience we 
also provide the functions {\tt p=newmesh(p)}, which interpolates the 
current $(u,\tau)$ to a new mesh generated after user input 
in the form {\tt hmax} or {\tt nx,ny}, 
and {\tt p=setbmesh(p)} which sets the 
base mesh to the current mesh. This should be called if it is expected 
that the current mesh is a good base for adaption in the steps to come. 
}
\eex\erem

\subsubsection{The linear system solvers}\label{lss-sec}
Recall that after discretization with $n_p$ points 
we have nodal values $u\in \R^p$ with $p=Nn_p$ large, and 
\huga{\label{aagain} 
G_u\in \R^{p\times p} \quad\text{and}\quad 
A=\bpm G_u&G_\lam\\
\xi\dot{u}&(1-\xi)\dot{\lam}\epm \in\R^{(p+1)\times (p+1)}
} 
are large, but sparse (block) matrices. The question is how to best 
solve $G_u v=r$ and the bordered systems such as
$A\tau=z$, respectively. 
In all the examples that we considered, our experience is that 
the highly optimized matlab solver {\tt z=A}$\backslash${\tt b} of 
$Az=b$ 
works remarkably well, but for easy customization of the code we never call 
$\backslash$ directly but use two interface routines: 
\ben
\item {\tt v=p.lss(M,r,p,lam)} to solve $M v=r$ with $M=G_u\in\R^{p\times p}$; 
\item {\tt z=p.blss(A,b,p,lam)} to solve $Az=b$ with 
$A\in\R^{{p+1}\times{p+1}}$. 
\een 
Here {\tt blss} and {\tt lss} stand for ({\tt b}ordered){\tt l}inear 
{\tt s}ystem {\tt s}olver. 

The default solvers {\tt lss} and {\tt blss} just contain one 
command, namely {\tt v=M$\backslash$r} resp.~{\tt z=A$\backslash$b}. 
Nevertheless, for large systems or for some special classes of problems iterative solvers might 
work better, and as templates we provide the two routines 
{\tt ilss} and {\tt iblss}, using {\tt gmres} with (incomplete) 
LU factorization {\tt luinc} 
as preconditioners. These should, of course, be reused 
as long as gmres converges quickly, and here (and in {\tt resinj.m}) we thus 
introduce some {\em global} variables, namely {\tt global L U;} 
resp.\,{\tt global bL bU}. Thus, when using, e.g., {\tt ilss} the user 
must also issue {\tt global L U; L=[]; U=[];} from the command line. 
\footnote{The reason for this construction is that we do not want to 
make $L,U$ a part of {\tt p} since this needs a lot of disk space 
when saving {\tt p}: typically, we get fill--in factors for $L,U$ of 
10 and larger. \label{lssfoot}} 

It turns out that for 
scalar problems these outperform the direct solvers {\tt lss} and 
{\tt blss} for large $n_p$, $n_p>10^5$, say. On the other hand, 
for systems, {\tt luinc} becomes exceedingly slow such that $\backslash$ 
beats {\tt gmres} with $LU$ preconditioning even for very large $n_p$. 
In summary, the iterative solvers {\tt ilss} and {\tt iblss} should 
only be regarded as template files to create problem specific 
iterative solvers when needed. See also \S\ref{acgc-sec} for 
adaptions of {\tt lss} and {\tt blss} to some special situation. 

Finally, various approaches have been proposed 
for the solutions of the bordered systems $Az=b$, see, e.g., 
\cite{keller77, govaerts}. As alternatives to {\tt blss} we provide {\tt bellss} 
(``bordered elimination'') and {\tt belpolss} (``bordered elimination plus one'').  
To use these, simply set, e.g., {\tt p.blss=@belpolss}. 
In our tests the performance of {\tt bellss} and {\tt belpolss} is roughly the same as {\tt blss}. 

\subsubsection{Screen output, plotting, convergence failure, auxiliary functions}
The screen output during 
runs is controlled by the two functions {\tt p.headfu} (headline) 
and the function {\tt p.ufu}. 
These are preset in {\tt stanparam} as {\tt p.headfu=@stanheadfu, p.ufu=@stanufu} 
to first print a headline and then, 
after each step, some useful information. To print 
some other information the user should adapt {\tt stanheadfu} and 
{\tt stanufu} to a local copy, say {\tt myhead.m}, and set 
{\tt p.headfu=@myhead}, and similar for {\tt stanufu} and {\tt p.ufu}.\footnote{For 
{\tt p.timesw>0} we also plot some timing information at the end 
of {\tt cont}.}
The bifurcation diagram and solution plots are also generated 
during continuation runs, but in general it is more convenient to 
postprocess via {\tt plotbra, plotsolf} etc.

The files {\tt p*.mat} and {\tt bp*.mat} contain the complete data of 
the respective point on a branch.\footnote{Obviously this is often quite 
redundant, but it is necessary if, e.g., some mesh refinement 
occured during continuation. To save disk space, however, we deliberately 
chose to {\em not} make $G_u$ a part of {\tt p}. See also 
footnote \ref{lssfoot}.} Thus, a run which is no longer 
in memory can be simply reloaded by, e.g., 
{\tt q=loadp(pre,pname,'q')}, where {\tt pre, pname} is the 
name data of a previously saved point, and the third argument is used 
to set the directory name for the newly created struct. The loaded point will often be either the last 
one or the first; 
in the latter case, to change direction of the branch, use, e.g., 
{\tt q=loadp('p','p1','q'); q.ds=-q.ds; q=cont(q);} 

If the Newton--loop does not converge even after reducing 
{\tt ds} to {\tt dsmin} then {\tt cfail.m} is called. The standard 
option is to simply abort {\tt cont}, but we offer 
a number of alternatives, e.g., to change some parameters like 
{\tt dsmin} or {\tt imax}, or to try, e.g., some mesh refinement or 
adaption. Clearly, the choice here is strongly problem dependent, 
and thus we recommend to adapt {\tt cfail.m} if needed; see 
\S\ref{custsec} for remarks on such ``customization without function 
handles''. 

Besides those already mentioned we provide further auxiliary functions,  
see \cite[m2html]{pde2phome} for a complete documented list.

\subsection{The Allen--Cahn equation with Dirichlet boundary 
conditions ({\tt ac})}\label{ac-sec}
\begin{figure}[ht!]
\bce 
{\small 
\begin{tabular}[t]{p{30mm}p{35mm}p{35mm}p{35mm}}
(a) 
\ig[width=30mm,height=25mm]{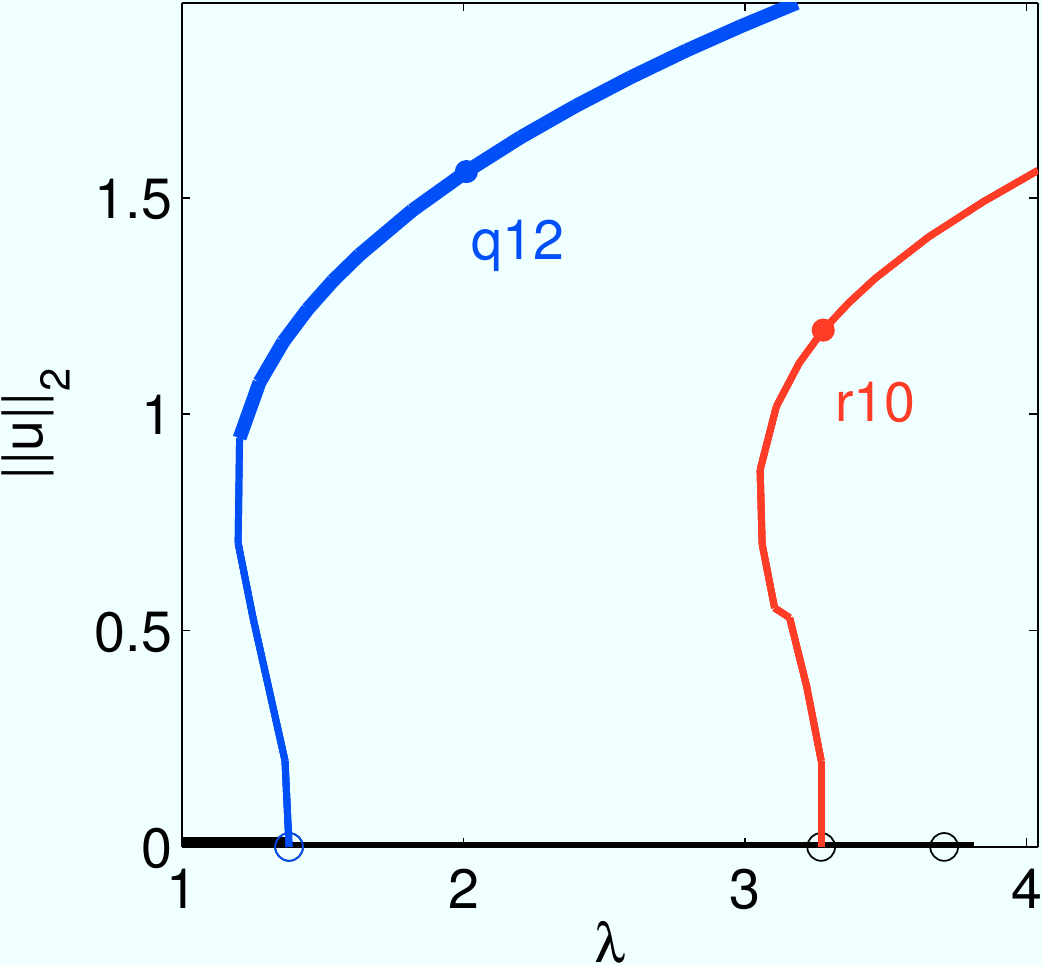}
&(b) \ig[width=35mm,height=25mm]{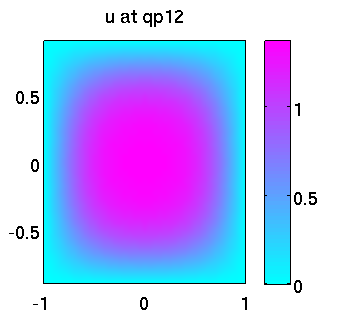}&
(c) \ig[width=35mm,height=25mm]{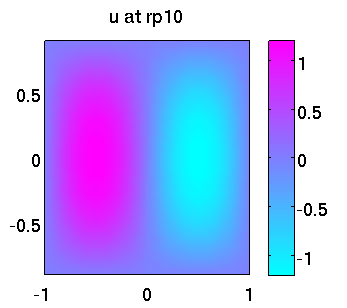}&
(d) \ig[width=35mm,height=25mm]{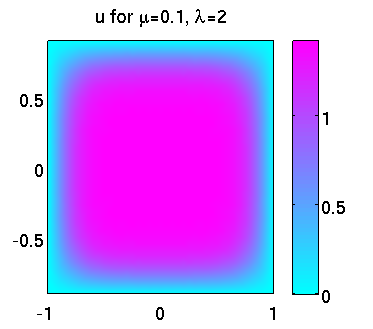}
\end{tabular} 
}
\ece

\vs{-5mm}
\caption{{\small  (a) Elementary bifurcation diagram for \reff{ace} with 
$\mu=0.25$. 
 No secondary bifurcations occur, 
and the mode--structure on each branch is completely determined at bifurcation.
(b),(c) some points on branches as indicated. 
(d) Solution after continuation in $\mu$ from (b) to $\mu=0.1$.
\label{acbdf}}}
\end{figure}

In our second example we use Dirichlet boundary conditions (DBC), 
and explain some ad hoc parameter switching, and time integration. 
We consider a cubic--quintic (to have folds) Allen--Cahn equation 
\huga{-\mu\Delta u-\lam u-u^3+u^5=0\quad\text{on}\quad \Om=[-L_x,L_x]
\times[-L_y,L_y], \quad  u|_{\pa\Om}=0, 
\label{ace}
}
with two parameters $\mu>0$ and $\lam\in\R$. We use 
{\tt [p.geo,bc]=recdbc1(lx,ly,1e3)} to approximate 
the DBC, and set $L_x=1$ and $L_y=0.9$ to 
break the square symmetry present in {\tt bratu} 
in order to have only simple bifurcations, 
namely at $\lam_{kl}=\mu \pi^2((k/L_x)^2+(l/L_y)^2)$. First we 
fix $\mu=0.25$  which yields 
$\lam_{11}=1.3784, \lam_{21}=3.2289, \lam_{12}=3.6630, \ldots$, 
and continue in $\lam$, which yields Fig.~\ref{acbdf}(a)--(c). 
After branch switching we turn on mesh--adaption after each 5 steps. 
See {\tt acdemo.m} or {\tt accmds.m} for more details, which also 
contain an example of perturbing a solution and subsequent time 
integration.

\subsubsection{Parameter switching}\label{para-sw-sec}
Unlike \auto, \pdep\ (currently) has no switches or presets for multi--parameter 
continuation. However, 
switching to a new parameter for continuation can be achieved in a 
simple and flexible way by  modifying the structure {\tt p} 
from the command line. As an example we want to continue in 
$\mu$ from point 10 on q to $\mu=0.05$. This is achieved by 
the commands in Table \ref{pstab}, and yields Fig.~\ref{acbdf}(d). 
The basic idea is to copy {\tt q} to {\tt w} (this is not strictly necessary) 
and then reset {\tt w.f} and {\tt w.jac}. We make our life 
simple by setting {\tt w.jsw=1} such that we do not need $G_\lam$, 
and set $\xi=10^{-6}$ since the dependence on $\mu$ is quite sensitive. 

We also introduce a new parameter {\tt w.up1} 
( '{\tt u}ser {\tt p}arameter 1', but any name will be fine) which is used 
to pass the current $\lam$ to {\tt acfmu},  see Table~\ref{pstab}. 

\begin{table}[!ht]
\bce{\small
\begin{tabular}{|p{160mm}|}
\hline\\\vs{-10mm}
\begin{verbatim}
w=q;w.up1=w.lam;w.lam=0.25;w.lammin=0.05;w=setfn(w); w.ds=-0.01;
w.f=@acfmu;w.f=@acjacmu;w.jsw=1;w.parasw=0;w.xi=1e-6;w.restart=1;w=cont(w);
 \end{verbatim}
\vs{-12mm}\\\hline
\vs{-5mm}
\begin{verbatim}
function [c,a,f,b]=acfmu(p,u,lam)  % AC for cont in mu
u=pdeintrp(p.points,p.tria,u); c=lam; a=0; f=p.up1*u+u.^3-u.^5; b=0;
 \end{verbatim}
\vs{-12mm}\\\hline
\vs{-5mm}
\begin{verbatim}
function [c,fu,flam,b]=acjacmu(p,u,lam) % Jacobian for cont of AC in mu; run with jsw=1
u=pdeintrp(p.points,p.tria,u); c=lam; flam=0; b=0; fu=p.up1+3*u.^3-5*u.^4; 
 \end{verbatim}
\vs{-12mm}\\\hline
\end{tabular}}
\ece

\vs{-2mm}
\caption{{\small Switching to continuation in $\mu$, commands, and modified 
coefficient and Jacobian functions. 
}\label{pstab}}\end{table}

\subsubsection{Time integration}\label{tint-sec}
For time integration of \reff{pprob} using the struct {\tt p} we 
provide a simple semi-implicit Euler method. Writing $u^{(n)}$ for 
$u(t_n,\cdot)$, choosing a time--step $h$, approximating 
$\pa_t u(t_n)\approx \frac 1 h(u^{(n+1)}-u^{(n)})$ where 
$t_n=t_0+nh$, and evaluating, e.g., 
$\nabla\cdot(c\otimes\nabla u)$ 
as $\nabla\cdot(c(u^{(n)})\otimes\nabla u^{(n+1)})$ 
we obtain, on the FEM level, 
\hualst{
&\frac 1 h M(u^{(n{+}1)}{-}u^{(n)})=-K(u^{(n)})u^{(n{+}1)}{+}F(u^{(n)})\\ 
\aqui\ & u^{(n+1)}=(M{+}hK(u^{(n)}))^{-1}(Mu^{(n)}{+}hF(u^{(n)})). 
}
Here $M$ is the mass matrix and $K$ is the stiffness 
matrix on time--slice $n$. This is implemented in 
{\tt tint(p,h,nstep,pmod)}, 
where {\tt nstep} is the number of time steps, and a plot (of component {\tt p.pcmp}) is 
generated each {\tt pmod}'th step. Thus 
we may call, e.g., {\tt p.u=p.u+0.1*rand(p.neq*p.np,1); 
p=tint(p,0.1,50,4,10)} to first perturb a given solution and 
then time-step. See {\tt accmds.m} for an example, where 
we perturb a solution on the unstable part of the {\tt q} branch 
into both directions of the unstable manifold; in the subsequent 
time integration the solution converges to the stable trivial solution or 
the stable {\tt q} branch, respectively, as expected. 
However, the main purpose of {\tt tint} is 
to generate (stable) initial data for continuation, i.e., after {\tt tint} call 
{\tt cont}. See also 
\S\ref{s:rbconv} for an example where {\tt tint} is used in this spirit. 

\brem\label{tintrem} In {\tt tint} we assemble $K(u^{(n)})$ and 
solve $(M+hK(u^{(n)}))u^{(n+1)}=g^{(n)}$ in each 
step by {\tt lss}. 
Clearly, for special cases this can be optimized: for instance, 
if $c,a,b$ do not depend on $u$, then the textbook approach would be 
to assemble $K$ at the start, 
followed by some incomplete LU--decomposition of $M+hK$ 
combined with some iterative solver. 
However, similar remarks as in \S\ref{lss-sec} apply, and 
thus we use the very 
elementary form above, 
but stress again that {\tt tint} 
in its present form is not intended for heavy time-integration. 
\eex\erem 

\subsection{The Allen--Cahn equation with mixed $\lambda$-dependent boundary conditions ({\tt achex})}
\label{s:acbc}
\begin{figure}[ht!]
\bce 
{\small 
\begin{tabular}[t]{p{35mm}p{40mm}p{35mm}p{35mm}} 
(a) \ig[width=35mm,height=35mm]{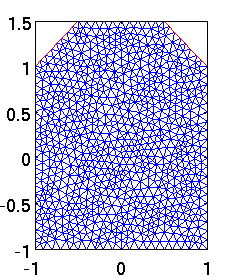}&
(b) \ig[width=40mm,height=35mm]{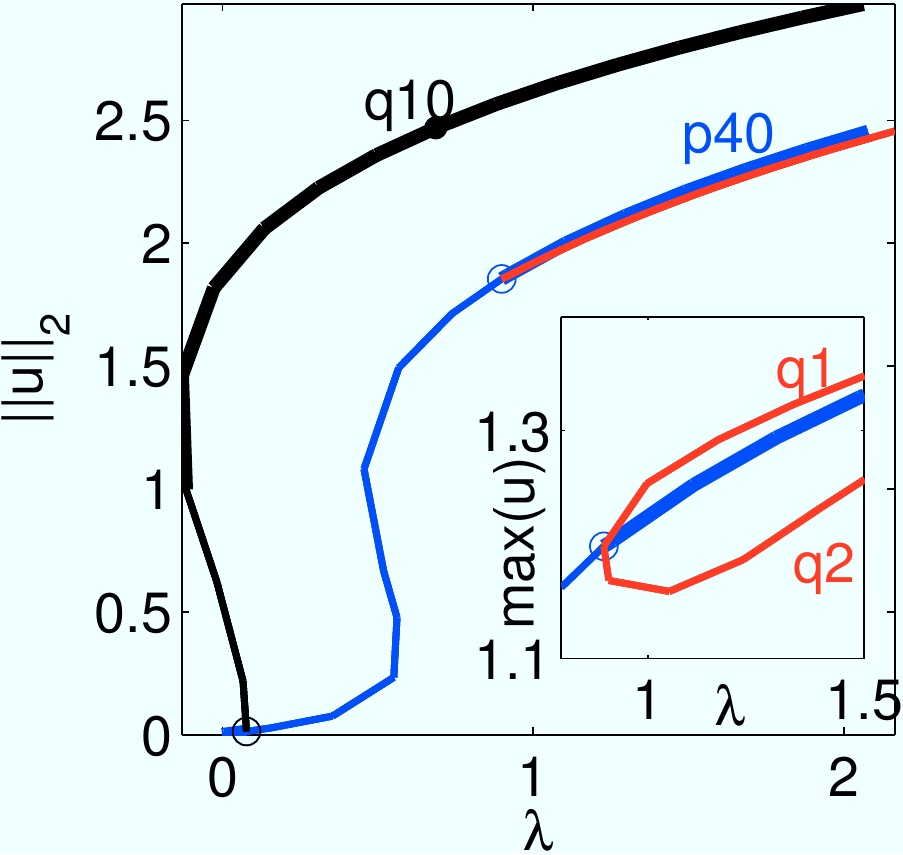}& 
(c) \ig[width=35mm,height=35mm]{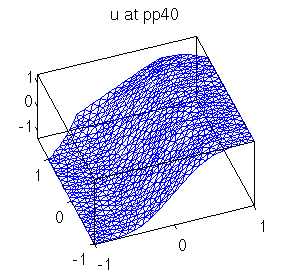}&
(d) \ig[width=35mm,height=35mm]{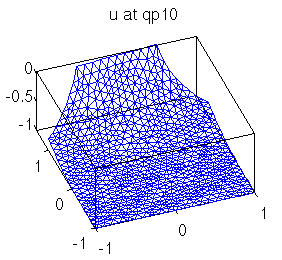}
\end{tabular}
}
\ece

\vs{-0mm}
\caption{{\small (a) A hexagonal domain, already with a basic mesh, created by 
{\tt p=stanmesh(p,0.1)}. (b) Bifurcation diagram for \reff{ace} with \bcs\ given 
by \reff{acb1}. (c),(d) selected solution plots. 
\label{acbcf}}}
\end{figure}

We illustrate a few more possibilities with \pdep\ by modifying the
Allen-Cahn example \eqref{ace} from above. We consider again 
\reff{ace}, i.e., 
$-0.25\Delta u-\lam u-u^3+u^5=0$, but instead of 
homogeneous Dirichlet \bcs\ on a rectangle, 
we consider hexagonal domains $\Om$ and parameter dependent mixed 
Dirichlet-Neumann \bcs. Figure \ref{acbcf}(a) shows an example for 
$\Om$, which basically consists of a square, 
with the top boundary shifted by $\del_y=0.5$ between $[-l_x,l_x]$, 
$l_x=0.5$. Denote this part of $\pa\Om$ by $\Ga_D$, and set 
$\Ga_N=\pa\Om\setminus\Ga_D$. 
To define the domain we could,  e.g., 
 use the {\tt pdetool} GUI to draw a polygon composed of six edges, one for
each segment and export the geometry. 
However, usually the function {\tt polygong} \cite{pruefert} is 
much more convenient.\footnote{see {\tt geo=hexgeo(lx,dely)}, which 
also contains a slightly edited output of the GUI for comparison.}  
Second, we want to define the boundary conditions 
\huga{\label{acb1}
n\cdot\nabla u=0 \mbox{ on }\Ga_N, \quad 
 u= \lambda x \mbox{ on } \Ga_D. 
}
To implement this we use a stiff spring approximation on $\Ga_D$ in 
via {\tt gnbcs}, i.e., 
\huga{\label{e:acbc6ximp}
{\small  \begin{array}{l}
    \mbox{{\tt qd=mat2str(10$\hat{\;}$4);gd=[mat2str(10$\hat{\;}$4*lam)~'*x']; qn='0'; gn='0';}}\\
    \mbox{{\tt bc=gnbcs(1,qn,gn,qn,gn,qn,gn,qn,gn,qd,gd,qn,gn);}}
  \end{array}
}
}
With \pdep\ we perform a continuation starting from the trivial zero
solution and obtain the bifurcation diagram plotted in
Fig.~\ref{acbcf}; see {\tt ac6cmds}. Bifurcation detection and branch 
switching work without problems, and the error estimate is always well 
below 0.01. To generate both parts of the {\tt r} branch we first 
call {\tt r1=swibra('p','bp2','r1',-0.1);r1=cont(r1)} and then 
{\tt r2=loadp('r1','p1','r2'); r2.ds=-r1.ds; r2=cont(r2)} to proceed 
in the other direction. At the end of {\tt ac6cmds} we also run an 
example with $u=\lam$ on $\Ga_D$ implemented via {\tt gnbc}. 

\subsection{A quasilinear Allen--Cahn equation ({\tt acql})}
\label{acqlsec}
To give an example of a more complicated Jacobian we modify \reff{ace} to the 
quasilinear Allen--Cahn equation
\huga{-\nabla\cdot[(0.25+\del u+\ga u^2)\nabla u] -f(u,\lam)=0\quad\text{on}\quad \Om=[-L_x,L_x]
\times[-L_y,L_y], \quad  u|_{\pa\Om}=0, 
\label{acql}
}
with $f(u,\lam)=\lam u+u^3-u^5$ and $L_x=1, L_y=0.9$ as before. 
See {\tt acqlf.m}. 
The linearization around $u$ gives the linear 
operator 
$$
G_u(u,\lam)v=-\nabla\cdot[(0.25+\del u+\ga u^2)\nabla v]+[-f_u(u,\lam)-\del\Delta u
-2\ga(\nabla u\cdot\nabla u+u\Delta u)]v-[(\del+2\ga u)\nabla u]\cdot \nabla v. 
$$
Hence, in {\tt acqljac.m} we now have 
${\tt fu}=f_u+\del\Delta u+2\ga(\nabla u\cdot\nabla u+u\Delta u)$, 
and ${\tt b}_{111}=(\del+2\ga u)u_x$ and 
${\tt b}_{112}=(\del+2\ga u)u_y$, cf.~Remark \ref{jacrem}. 
To generate $(u_x,u_y)$ and $\Delta u$ as coefficients in {\tt acqljac.m} 
we use {\tt pdegrad} resp.~{\tt pdegrad, pdeprtni} and {\tt pdegrad} 
again, see \cite{pdetdoc}. 

\begin{figure}[ht!]
\bce 
\begin{tabular}{llll}
\ig[width=38mm]{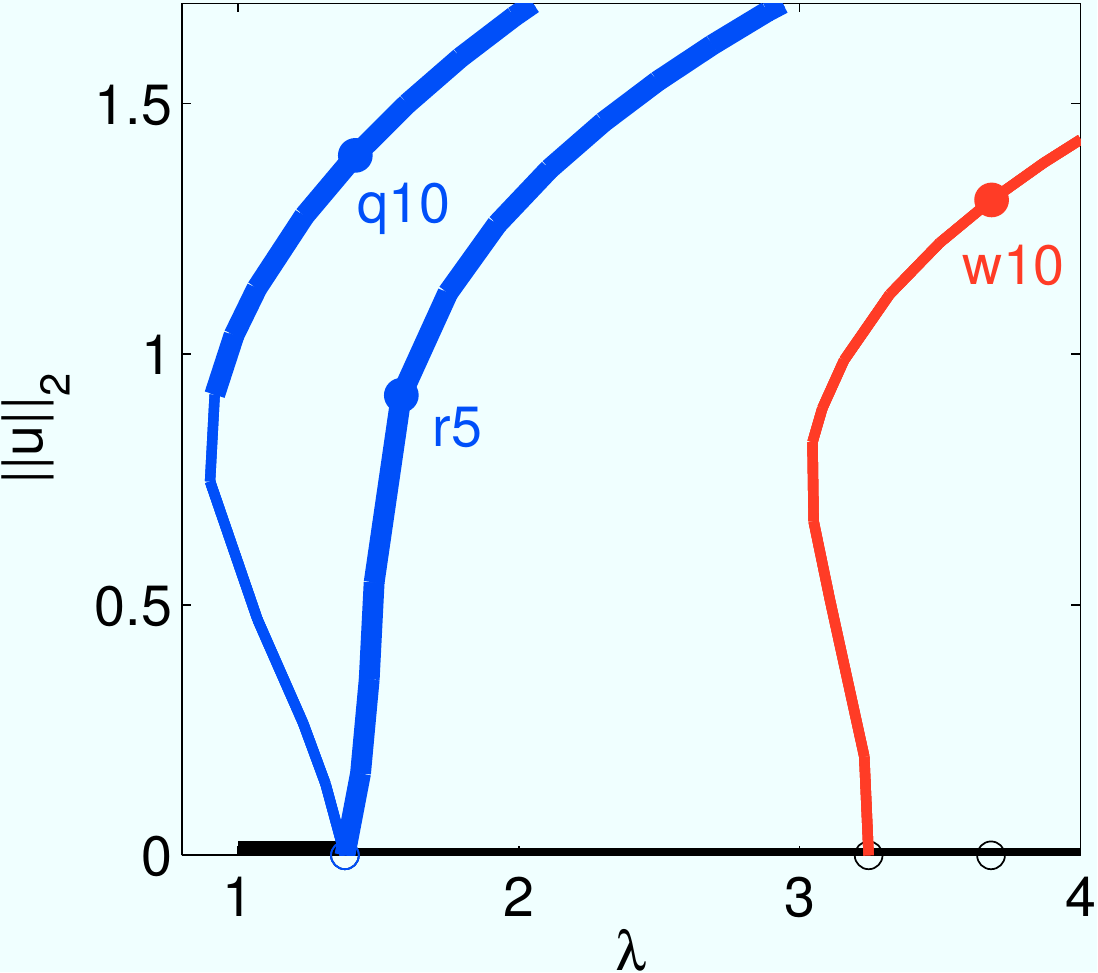}
&\ig[width=36mm]{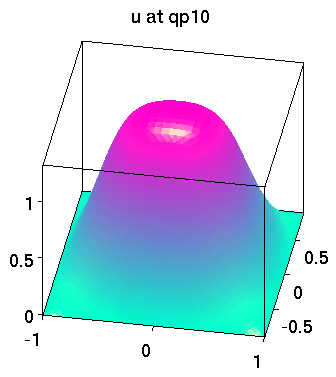}
&\ig[width=36mm]{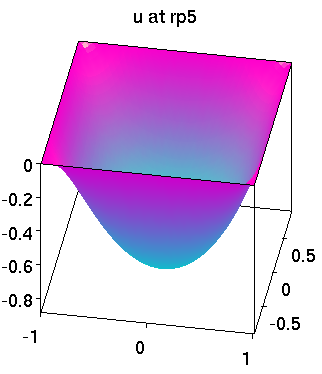}
&\ig[width=36mm]{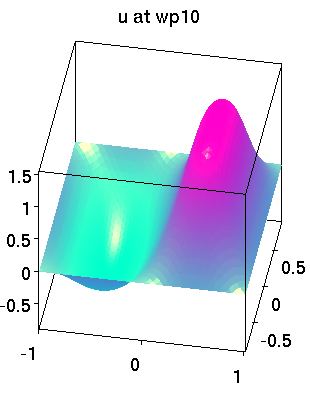}
\end{tabular}
\ece
\caption{{\small Elementary bifurcation diagram for 
\reff{acql} with $\del=-0.2$ and $\ga=0.05$, and 
some solution plots. The two blue branches are in fact one 
branch, and the corner at the transcritical bifurcation from the 
trivial branch is due to the choice of vertical axis. 
The symmetry $u\mapsto -u$ no longer holds, and ``up'' 
humps are steeper than ``down'' humps due to $\delta<0$. 
The {\tt w}--branch is still double due to the $x\mapsto -x$ 
symmetry. }\label{acqlbdf}}
\end{figure}

The term $\del u$ in $c$ changes 
the $u\mapsto -u$ symmetry of the Allen-Cahn equation \reff{ace}. 
The bifurcation points from the trivial branch 
$u\equiv 0$ in \reff{acql} are as in \reff{ace}, but the bifurcations change, 
see Fig.~\ref{acqlbdf}. 
In particular the first bifurcation changes from pitchfork to 
transcritical. 

\subsection{An Allen--Cahn equation with global coupling ({\tt acgc})}\label{acgc-sec}
As an example of a ``non--standard'' elliptic equation 
we treat an Allen--Cahn equation 
with a global coupling. 
We fix $\mu=0.1$ and $\lam=1$ in \reff{ace}, introduce a 
new parameter (again called $\lam$) and consider 
\huga{
G(u,\lam):=
-0.1\Delta u-u-u^3+u^5-\lam\spr{u}=0\quad\text{on}\quad 
\Om=[-\pi/2,\pi/2]^2, \quad  u|_{\pa\Om}=0,
\label{acgc}
}
where $\spr{u}=\int_\Om u\dd x$. The term $\lam\spr{u}$ is called a global coupling 
or global feedback, positive for $\lam>0$ resp.~negative for $\lam<0$. 
Problems with global coupling occur, e.g., 
in surface catalysis, where global coupling 
arises through the gas phase \cite{rbmi96}, in semi-conductors and 
gas-discharges \cite{woesler,stoll}, and as ``shadow systems'' 
in pattern formation when there is a very fast inhibutor diffusion 
\cite{ward2000}. 

The global feedback does not fit into the framework of 
\reff{gform} if $f$ is assumed to be local. For the definition of $G(u)$ 
this is not yet a problem as we may simply define $f$ as, e.g,
$$ 
\verb# f=u+u.^3-u.^5+lam*triint(u,p.points,p.tria);#
$$ 
where {\tt triint(g,points,tria)} 
is the Riemann sum of $\int g(x)\dd x$ over the given mesh. However, 
for continuation we make extensive use of Jacobians, and $G_u(u)$ 
is now given by 
$$
[G_u(u)v](x)=-0.1\Delta v(x)-(1+3u(x)^2-5u^4(x))v(x)-\lam\spr{v}. 
$$
As yet we cannot deal with last term, cf.~Remark \ref{jacrem}. 
The first try would be to simply ignore it in continuation, 
but this in general only works for small $|\lam|$ while for larger $|\lam|$ we 
loose convergence in the (false) Newton loop. 
We can express $\spr{v}$ on the FEM level 
via a matrix $M$ such that $G_u(u)v=(K-\lam M)v$. Essentially, 
for ``natural parametrization'' we need to solve 
\huga{\label{acls} 
G_u(u)v=r, \text{ where } G_u(u)=(K-\lam \nu \eta^T) \text{ with }
\nu,\eta\in\R^{n_p}. 
} 
Here $\eta=(a T)^T$ where $(a_1,\ldots,a_{n_t})$ contains the 
triangle areas, $T\in \R^{n_t\times n_p}$ 
interpolates $u\in\R^{n_p}$ from nodal values to triangle values
(such that $\spr{u}=${\tt triint(g,points,tria)=eta*u}), 
and $\nu_i=\int_\Om 1\phi_i\dd x$ corresponds 
to adding $\spr{v}$ to all nodes with the correct weight. 
However, $(K-\lam\nu\eta^T)$ is a full matrix and should never even be formed. 
Instead we customize {\tt lss} to use a Sherman--Morrison formula which gives 
(for \reff{acls}) 
\huga{\label{smf} 
v=K^{-1}r+\al(K^{-1}\nu)(\eta^TK^{-1})r, \quad \al=\frac \lam 
{1-\lam\eta^TK^{-1}\nu}.
}
In {\tt acgcjac.m} we then just ignore the term $\lam\spr{u}$. 
Similar remarks apply to the bordered systems solved by {\tt blss}.

\begin{figure}[htbp]
\bce 
{\small 
\begin{tabular}[t]{p{27mm}p{29mm}p{30mm}p{30mm}p{30mm}} 
(a) \ig[width=27mm,height=30mm]{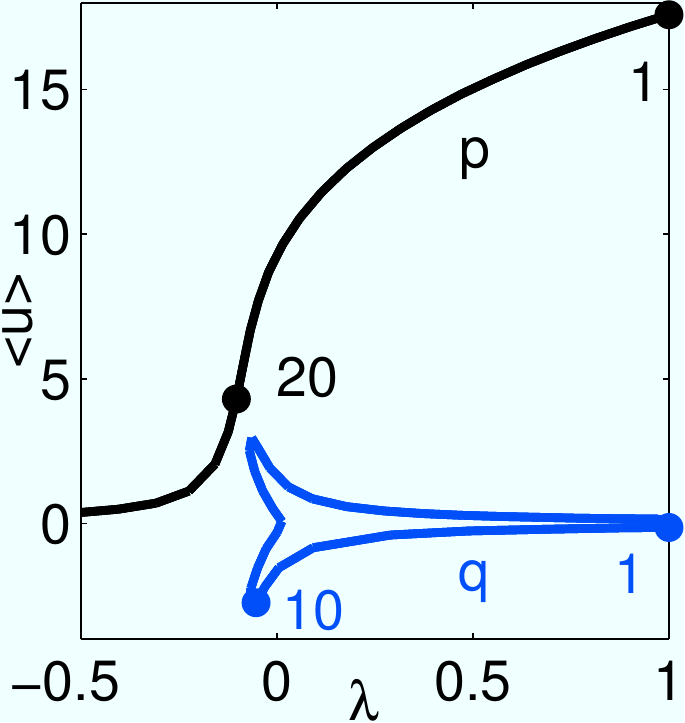}&
(b) \ig[width=29mm,height=35mm]{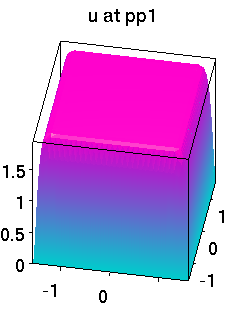}& 
(c) \ig[width=30mm,height=35mm]{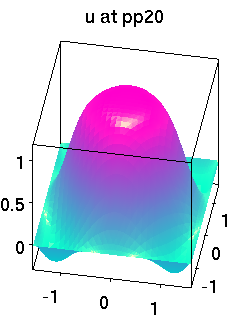}&
(d) \ig[width=30mm,height=35mm]{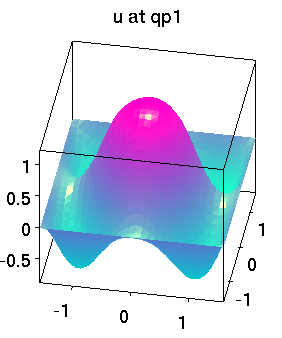}&
(e) \ig[width=30mm,height=35mm]{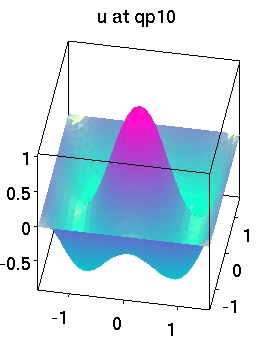}
\end{tabular}
}
\ece

\vs{-0mm}
\caption{{\small (a) Two solution branches for \reff{acgc},  and 
four selected solutions. By positive global feedback, the plateau in (b) 
($u$ around $1.93$) is substancially above the zero $(1+\sqrt{5})/2\approx 1.62$ of 
$f(u)=u+u^3-u^5$. Here, some mesh refinement near the boundary is 
also crucial. Decreasing $\lam$ to slightly negative values $u$ gets pushed 
below $0$ near the boundary (c). On the other branch some somewhat localized 
solutions are found (d), (e). Note 
that \reff{acgc} is symmetric w.r.t.~$(u,\lam)\mapsto (-u,\lam)$. 
\label{acgcf}}}
\end{figure}

In the actual implementation we introduce {\em global variables} {\tt nu,eta}. 
The idea is that it is sufficient to calculate $\nu$,$\eta$ once for a given 
mesh, for instance in {\tt acgcf.m}, as this is always called before the 
Jacobian {\tt acgcjac} or the linear system solvers {\tt gclss} or {\tt gcblss}.  
If we set aside mesh--refinement then we could calculate $\nu,\eta$ 
at startup and store them e.g.~as {\tt p.nu,p.eta} but with mesh refinement 
global variables are more convenient. See  Table \ref{acgctab} for the full code for 
{\tt acgcf.m, lss.m} for this example, 
and Fig.~\ref{acgcf} for the result of the basic continuation runs 
contained in {\tt acgccmds.m}. We switch off spectral 
calculations and bifurcation checks by setting {\tt spcalcsw=0; bifchecksw=0;} 
since out of the box these would be based on the (wrong) local Jacobian, 
and the two branches were generated by using two different starting points.

\begin{table}[!ht]
\bce{\small
\begin{tabular}{|p{160mm}|}
\hline\\\vs{-10mm}
\begin{verbatim}
function [c,a,f,b]=acgcf(p,u,lam)  % AC global coupling 
global eta nu; try se=size(eta,2); catch; eta=[]; se=0; end
if(se~=size(u,1))  % eta not yet set, or mesh is refined 
  C=n2triamat(p.points,p.tria); ta=triar(p.points,p.tria); eta=ta*C;
  [M,nu]=assempde(p.bc,p.points,p.edges,p.tria,0,0,1); end
um=eta*u; u=pdeintrp(p.points,p.tria,u); c=0.1; a=0; b=0; f=u+u.^3-u.^5+lam*um; 
\end{verbatim}\vs{-8mm}\\\hline\hline  \vs{-5mm}
\begin{verbatim}
function x=gclss(A,b,p,lam) % lss for AC with global coupling, Sherman-Morrison
global eta nu; y=A\b; z=A\nu; al0=lam*eta*z; al=lam*eta*y/(1-al0); x=y+al*z; 
\end{verbatim}
\vs{-8mm}\\\hline
\end{tabular}}
\ece

\vs{-2mm}
\caption{{\small Definition of $f$ and customized {\tt gclss.m}; 
see also {\tt acgcjac.m} and {\tt gcblss.m}. 
}\label{acgctab}}\end{table}

\subsection{First summary, and some remarks on 
customization}\label{custsec}
We end this introductory section based on scalar examples with a first 
summary and some implementation remarks. 

The {\tt p.f=@..} syntax has the advantage that multiple version of {\tt f} 
can be maintained and switching can be done by only changing {\tt p.f=@..}. 
On the other hand, we do not want to overwhelm the user with 
such options, and thus we restricted the ``user--definable'' functions 
to {\tt p.f,\ldots,p.headfu} from Table \ref{tab2}, 
where in fact in most cases the user only needs to 
set up {\tt p.f}, and {\tt p.jac} for {\tt p.jsw}$\le 2$. Nevertheless, 
as outlined above 
any function of \pdep\ can be customized for a given problem  by just 
 copying it  from {\tt ../p2plib/} 
to the current directory (where \mlab\ searches first) and then modifying it. Main 
candidates for customization are, e.g., 
{\tt plotbra.m, plotsol.m} if additional features/options are desired in plotting 
the bifurcation diagram or/and the solutions. See, e.g., \S\ref{gpsec} and 
\S\ref{s:rbconv}. 

 Most functions of \pdep\ only require a few input/output arguments. 
An important exception is {\tt plotbra(p,wnr,cmp,varargin)} where {\tt varargin} is 
a possibly long list 
of argument/value pairs. See {\tt plotbra.m} for a detailed description, 
and also the various {\tt plotbra} example calls in the demos. 

As mentioned, by default there are no global variables in \pdep, 
with the exceptions {\tt pj,lamj} which are set for 
numerical differentiation in {\tt resinj}, and possibly LU preconditioners 
for iterative linear system solvers, see \S\ref{lss-sec}. On the other hand, 
e.g., {\tt p.f}, {\tt p.jac}, {\tt p.lss} etc.~do not return {\tt p}. 
This is to have a somewhat clean distinction between functions for specific 
calculations and function like {\tt p=cont(p), p=swibra(...), p=meshref(p)} 
which modify the structure {\tt p}, including the mesh. 
As a result, the user might want to introduce some global variables 
to streamline calculations, see \S\ref{acgc-sec} for an example. These 
should then be declared before initialization of {\tt p}. 

For convenience, in Table \ref{cfig} we summarize the typical 
steps in the usage of the software.  

\begin{table}[!ht]
{\small 
\bce
\begin{tabular}{|p{160mm}|} 
\hline
{\bf Initialization.} Declare (user defined) global variables (if any). 
Initialize structure {\tt p}, typically by first calling {\tt p=stanparam(p)}, 
followed by problem dependent calls to define (function handles for the) 
PDE coefficients, \bcs\ and Jacobian, and the geometry, mesh, and starting 
point.  \\ \hline 
{\bf {\tt function p=cont(p)}}\\[-6mm]
\ben 
\item If restart=1 then {\tt inistep}: generate first two points on branch and (secant) $\tau_0$
\item Predictor $(u^1,\lam^1)=(u_0,\lam_0)+\rds\tau_0$ with stepsize ds
\item Corrector: depending on {\tt parasw} and $\dot{\lam}_0$ 
use {\tt nlooppde} for \reff{trn} or {\tt nloopext} for \reff{newton} or 
\reff{chord}. This uses {\tt getder, getGu} resp.\,{\tt getGlam} to obtain 
derivatives, and {\tt lss} resp.\,{\tt blss} as linear systems solver. 
\item Call {\tt sscontrol} to assess convergence ({\tt res,iter} returned 
from  {\tt nlooppde} resp.~{\tt nloopext}): 
If {\tt res}$\le${\tt p.tol} accept step, i.e., goto 5, 
(and increase {\tt ds} if {\tt iter<imax/2}). 
If {\tt res>p.tol} and {\tt ds> dsmin} then decrease {\tt ds} 
and goto 2.  
If {\tt res>p.tol} and {\tt ds=dsmin} then no convergence, hence 
call {\tt cfail}.
\item Postprocessing: calculate new tangent $\tau_1$ by 
\reff{tau1}, call {\tt spcalc} (if spcalcsw=1), {\tt bifdetec} (if {\tt bifchecksw=1}). 
Check for error and mesh adaption. 
Update {\tt p}, i.e., put $u_0=u_1$, $\lam_0=\lam_1$, 
$\tau_0=\tau_1$ into {\tt p}, call {\tt out=outfu(p,u,lam)}, 
plot and save to disk. Call {\tt p.ufu} for printout and further user-defined actions.
\item If stopping criteria met ({\tt p.ufu} returned 1 or 
stepcounter$>${\tt nsteps}) 
then stop, else next step, i.e., goto 2.
\een 
\\ \hline
{\bf Post--processing.} Plot bifurcation diagrams via {\tt plotbra} 
({\tt plotbraf}) 
and solutions  
via {\tt plotsol} ({\tt plotsolf}). 
If bifurcations have 
been found, use {\tt swibra} (and {\tt cont} to follow some of these). 
\\ \hline 
\end{tabular}
\ece
}
\caption{Typical software usage, including pseudo--code of 
{\tt p=cont(p)}, with main function calls. \label{cfig}}
\end{table}

\section{Some prototype Reaction--Diffusion Systems}
\label{rd-sec}
Pattern--formation in  Reaction--Diffusion Systems (RDS), in particular 
from mathematical biology \cite{Mur}, 
is one of the main applications of path-following and bifurcation 
software. Here we first consider 
a quasilinear  two-component system with ``cross diffusion'' from chemotaxis 
\cite{mm91} to explain the setup of $c$ in this rather general case, 
and the setup of general domains in \pdep. 
We essentially recover the bifurcation diagrams from 
\cite{mm91} without special tricks or customization. 

Our second example is the Schnakenberg model \cite{schnak}, which is semilinear 
with a diagonal constant diffusion matrix, and thus in principle simpler 
than the first example. However,  here we are interested in a more complete 
bifurcation picture, and the Schnakenberg model shows many 
bifurcations already on small domains. Therefore 
we need some adaptions of the basic {\tt cont} algorithm to 
{\tt pmcont} ({\tt p}arallel {\tt m}ulti
{\tt cont}inuation), and we introduce {\tt findbif} to locate some 
first bifurcations from the homogeneous branch.

\subsection{Chemotaxis} 
\label{chem-sec}
An interesting system from chemotaxis has been analyzed in \cite{mm91}, 
including some numerical path-following and bifurcations using ENTWIFE. 
The (stationary) problem reads 
\huga{\label{chem1} 
0=G(u,\lam):=-\bpm D\Delta u_1-\lam \nabla\cdot(u_1\nabla u_2)
\\\Delta u_2\epm 
-\bpm r u_1(1-u_1)\\ \frac{u_1}{1+u_1}-u_2\epm. 
}
where $\lam\in\R$ is called the chemotaxis coefficient and 
$D>0$ and $r\in\R$ are additional parameters. In \reff{chem1} 
we  have $b\equiv 0$, may set $a=0$, and 
identify the second term with $f(u)$. The 
linearization reads 
\huga{\label{chem1j}
G_u(u,\lam)\bpm v_1\\ v_2\epm =
\left[{-}\bpm D\Delta&{-}\lam\nabla\cdot (u_1\nabla \cdot)\\
0&D\Delta\epm
+\bpm r(2u_1{-}1)+\lam \Delta u_2&0\\
{-}(1+u_1)^{{-}2}&1\epm\right]
\bpm v_1\\ v_2\epm +\lam\bpm \nabla u_2\cdot \nabla v_1 \\0 \epm, 
}
and in the notation from Remark \ref{jacrem}, 
the first matrix in \reff{chem1j} relates to {\tt c}, 
the second  to $-{\tt fu}$, and the last term gives $-b\otimes \nabla v$ 
with $b_{111}=-\lam \pa_x u_2$, $b_{112}=-\lam\pa_y u_2$, 
and $b_{ijk}=0$ else. 

\subsubsection{Bifurcation diagram over rectangles ({\tt chemtax})}
 Following \cite{mm91} we first study \reff{chem1} on a rectangular 
domain $\Om=[-L_x/2,L_x/2]\times [-L_y/2,L_y/2]$ 
with homogeneous Neumann \bcs. Again a number of results can then be obtained 
analytically. There are two trivial stationary branches, namely $u=(0,0)$ 
which is always unstable, and $u=u^*=(1,1/2)$. From the \bcs, the eigenvalue 
problem $M v=\mu v$ for the linearization around $u^*$ has 
solutions of the form $\mu=\mu(m,l,\lam)$, $v=v(m,l,\lam;x)
=\phi e_{m,l}(x,y)$ with $\phi\in\R^2$ 
and $\ds e_{m,l}(x,y)=\cos\left(\frac {m\pi}{L_x}(x+\frac{L_x}{2})\right)
\cos\left(\frac {m\pi}{L_y}(y+\frac{L_y}2)\right)$, 
$(m,l)=(1,0),(0,1),(2,0),(1,1), \ldots$.  
To study bifurcations from $u^*$ we solve $\mu(m,l,\lam)\stackrel{!}{=}0$ 
for $\lam$ which yields 
$$ \lam_{m,l}:=4(Dk^2+r)(k^2+1)/k^2, \text{ where }
k^2:=\pi^2\left(\frac {m^2}{L_x^2}+ \frac {l^2}{L_y^2}\right). $$ 
As in \cite[Fig.3]{mm91} we choose 
 $D=1/4$, $r=1.52$ and the ``$1\times 4$'' 
domain $L_x=1$, $L_y=4$, which yields Table \ref{cbiftab1}. 
\begin{table}[!ht]\bce
\begin{tabular}{|c|c|c|c|c|c|c|c|}
\hline
$(m,l)$&(0,2)&(0,3)&(0,1),&(1,0),(0,4)&(1,1)&(1,2)&\ldots\\\hline
$\lam_{ml}$&12.01 &13.73&17.55&17.57&18.15&19.91&\ldots\\
\hline
\end{tabular}
\caption{{\small Bifurcation from $u^*=(1,1/2)$ in \reff{chem1}, 
$D=1/4$, $r=1.52$.}\label{cbiftab1}}
\ece
\end{table}

\begin{figure}[htbp]
\bce 
{\small 
\begin{tabular}[t]{p{60mm}p{20mm}p{20mm}}
(a)&(b)&(c)\\
 \ig[width=50mm,height=40mm]{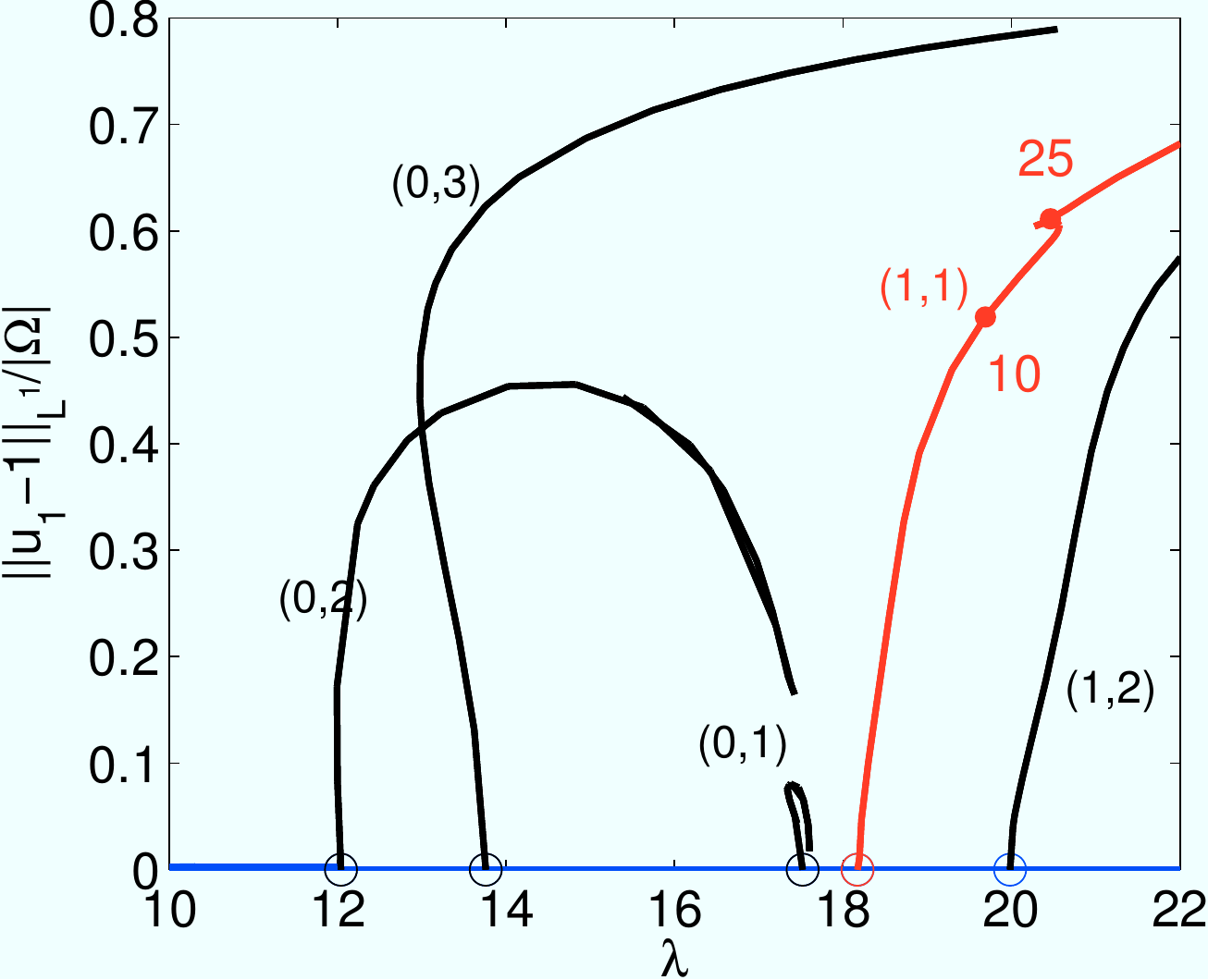}&
\ig[width=17mm,height=40mm]{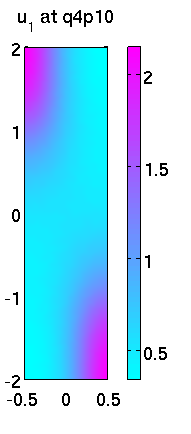}&
 \ig[width=17mm,height=40mm]{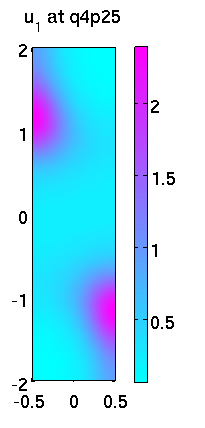}
\end{tabular}}
\ece

\vs{-4mm}
\caption{{\small (a) Bifurcation diagram for \reff{chem1} with {\tt jsw=3}, 
i.e., numerical Jacobians, and {\tt nt=2376}.  Of the bifurcating branches 
only the $(0,2)$--branch is stable in a certain $\lam$ range, 
and a number of secondary bifurcations occur on each branch. 
(b),(c) The shape of solutions before and after the loop on the $(1,1)$ branch. 
For {\tt jsw=1} we need finer meshes (${\tt nt}\approx 10^4$ to $2\cdot 10^4$ 
and 
adaptive refinement), which, 
while giving smaller error-estimates, also 
destroy the speed advantage of assembled Jacobians. 
\label{ctbdf}}}
\end{figure}

To encode \reff{chem1} we note that 
$c_{1111}=c_{1122}=D, c_{1211}=c_{1222}=-\lam u_1, c_{2211}=c_{2222}=1$, 
and all other entries of $c$ are zero. In particular, $c$ is isotropic 
and thus we may use {\tt isoc.m} to encode it, see Table \ref{chemtab1}
for {\tt chemf.m}. For convenience and illustration, 
here by default we first use {\tt p.jsw=3} in {\tt cheminit.m} such that 
{\tt p.jac} need not be set. 
With {\tt p=stanmesh(p,0.075)} leading to {\tt p.nt=2376} this still gives 
quick results, which moreover essentially do not change under mesh refinement. 
Also, in {\tt cheminit} we introduce 
{\tt p.vol}=$|\Om|$; we want to use this quantity in {\tt chembra.m} 
since the bifurcation diagrams in \cite{mm91} plot 
$\|u_1-1\|_{L^1}/|\Om|$ over $\lam$. Again this is a simple example 
that the user can augment the structure {\tt p} with whatever is useful. 
\begin{table}[!ht]
\bce{\small
\begin{tabular}{|p{160mm}|}
\hline\\\vs{-10mm}
\begin{verbatim}
function [c,a,f,b]=chemf(p,u,lam) % chemotaxis system with isoc 
u=pdeintrp(p.points,p.tria,u);a=0;b=0; v1=ones(1,p.nt); 
f1=r*u(1,:).*(1-u(1,:));f2=u(1,:)./(1+u(1,:))-u(2,:); f=[f1;f2];
D=0.25;r=1.52; c=isoc([[D*v1 -lam*u(1,:)];[0*v1 v1]],p.neq,p.nt);
\end{verbatim}
\vs{-8mm}\\\hline
\end{tabular}}
\ece

\vs{-2mm}
\caption{{{\tt chemf.m} as a prototype for definition of PDE coefficients 
in case of a (nonsymmetric) $c$ depending on $u$ and $\lam$. See 
{\tt isoc} and the {\tt assempde} documentation for the order of $c_{ijkl}$ in {\tt c}.} 
\label{chemtab1}}\end{table}
The commands in {\tt chemcmds.m} yield the bifurcation diagram in 
Fig.\ref{ctbdf}, where the bifurcation values $\lam_{m,l}$ (except 
for $\lam_{10}=\lam_{04})$ are found with resonable accuracy, and 
which agrees well with \cite[Fig.3(a)]{mm91}, with one exception:  
on the $(1,1)$ branch there is a loop near $\lam=20.5$ 
with two bifurcations, during which the solution structure changes 
as detailed in (b),(c). Presumably, this loop was just missed in 
\cite{mm91} due to a larger stepsize. 

Alternatively, to run \reff{chem1} with {\tt jsw=1} 
we also provide {\tt chemjac.m} which 
encodes \reff{chem1j}. For this, however, we need considerably 
finer meshes, mainly since the calculation 
of the coefficient $\Delta u_2$ (needed for {\tt jsw<2}) 
via {\tt pdegrad} and {\tt pdeprtni} does 
not go together well with Neumann boundary conditions, since 
the averaging involved in {\tt pdeprtni} produces some error 
at the boundaries. Therefore, we also replace {\tt p=cheminit(p)} 
by {\tt p=cheminitj(p)}, which resets a number of switches to 
(re)run \reff{chem1} with {\tt jsw=1}. See the end 
of {\tt chemcmds.m} resp.~{\tt chemdemo.m}.

\subsubsection{Drawing general domains ({\tt animalchem})}
\begin{figure}[ht!]
\bce 
{\small 
\begin{tabular}[t]{p{35mm}p{35mm}p{35mm}p{35mm}}
(a) \ig[width=35mm,height=35mm]{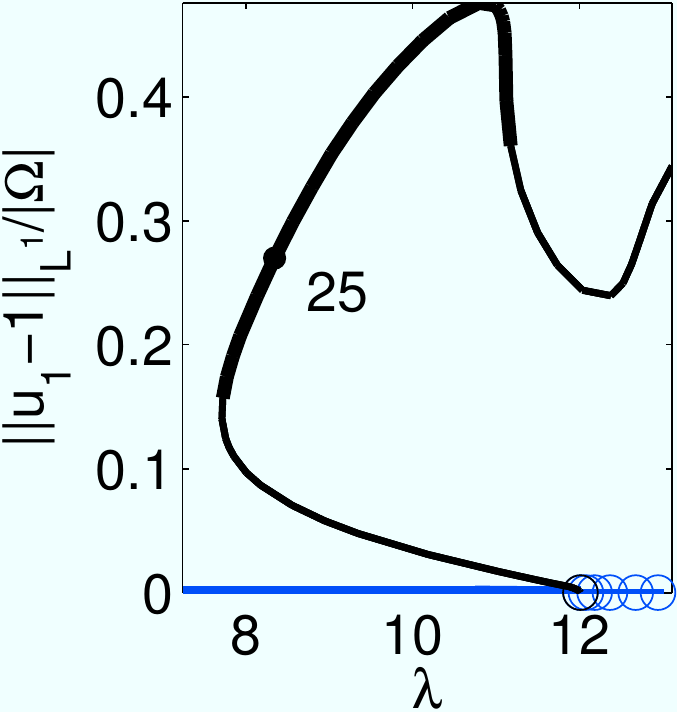}&
(b) \ig[width=35mm,height=35mm]{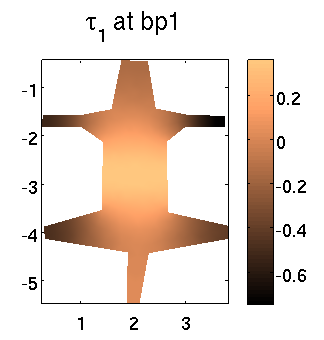}& 
(c) \ig[width=35mm,height=35mm]{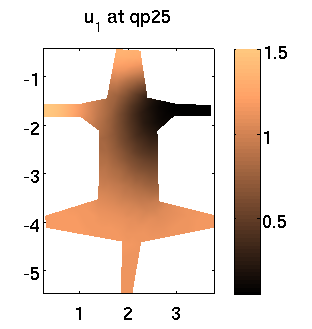}&
(d) \ig[width=35mm,height=35mm]{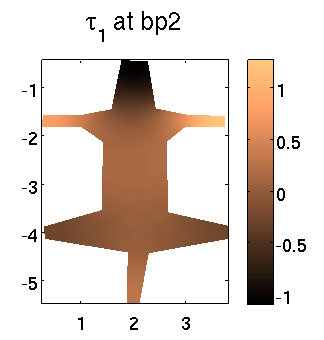}
\end{tabular}
}
\ece

\vs{-4mm}
\caption{{\small Bifurcation diagram for \reff{chem1} 
($\|u_1-1\|_{L^1}/|\Om|$ over $\lam$), first 
bifurcation direction from the trivial branch, 
stable solution on the bifurcating branch, and the second bifurcation 
from trivial branch. The bifurcation points on the trivial branch are quite 
close. Thus, on the trivial branch we need to run {\tt cont} with 
small {\tt dlammax}. \label{abdf}}}
\end{figure}

We now consider \reff{chem1} on the animal--shaped 
domain in Fig.\ref{abdf}, taken from \cite{Mur}, with 
Neumann \bcs. To set up $\Om$ we proceed graphically as explained in 
\S\ref{geosec}, see {\tt animalgeo.m}, also for the setup of the \bcs. 
The plots in Fig.\ref{abdf} are generated from the commands 
in {\tt animalcmds.m}. For problems of this type, the bifurcation directions 
from a trivial branch are often most interesting.

\subsection{The Schnakenberg model ({\tt schnakenberg})}\label{schnak-sec}
\begin{figure}[htbf]

\begin{minipage}{0.16\textwidth}
\includegraphics[width=\textwidth]{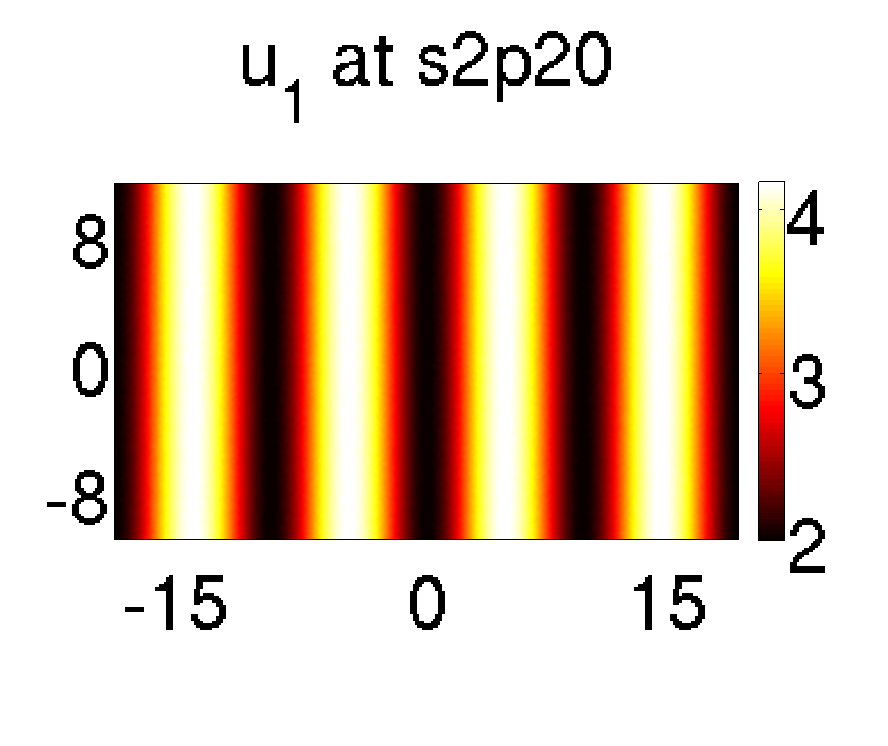}\\
\includegraphics[width=\textwidth]{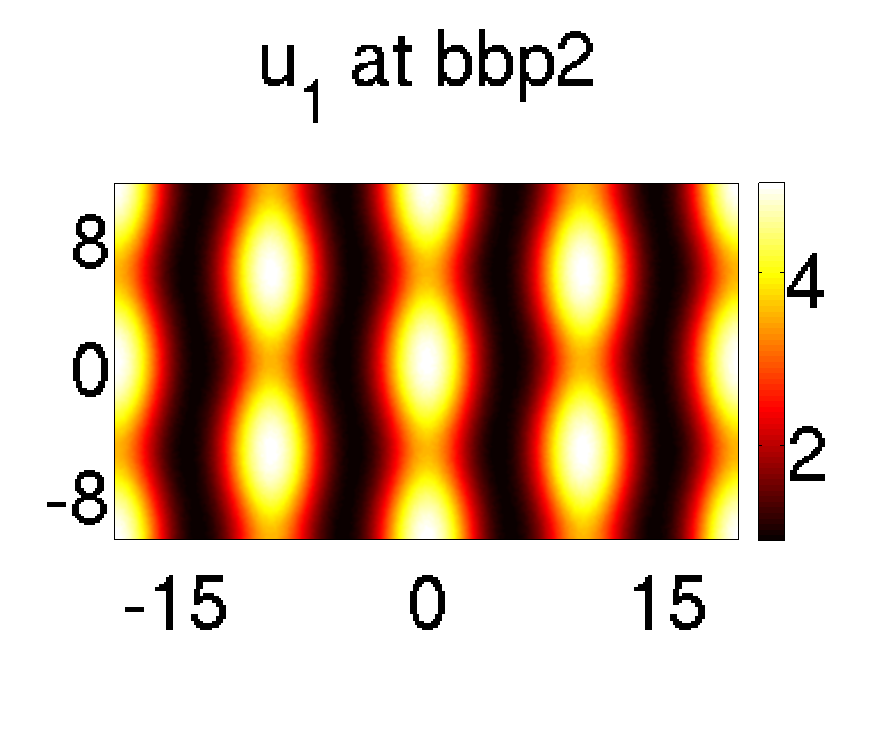}
\end{minipage}
\begin{minipage}{0.16\textwidth}
\includegraphics[width=\textwidth]{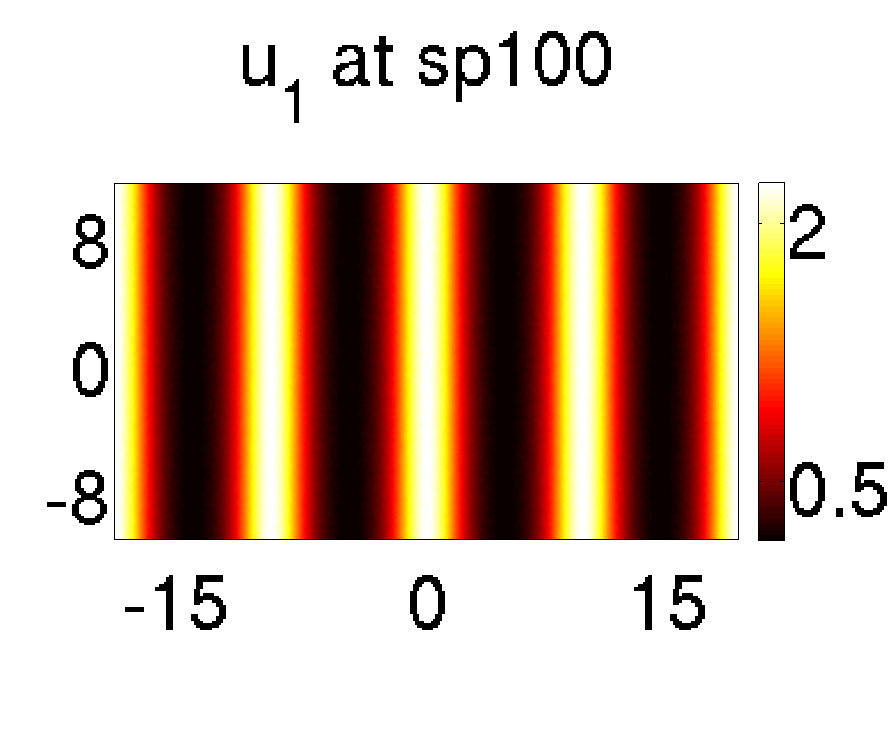}\\
\includegraphics[width=\textwidth]{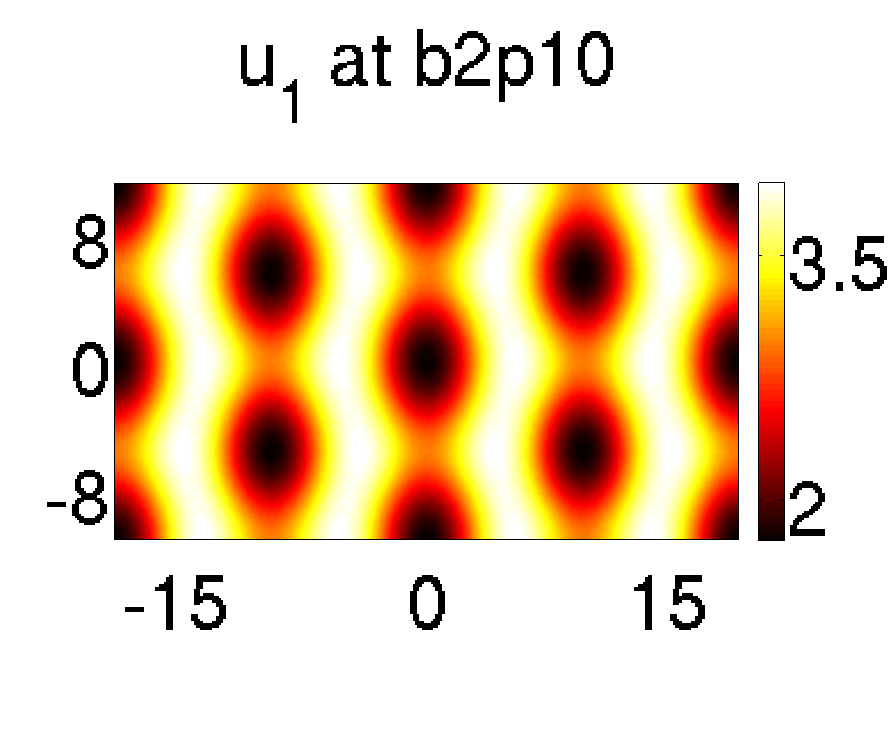}
\end{minipage}
\begin{minipage}{0.5\textwidth}
\includegraphics[width=\textwidth]{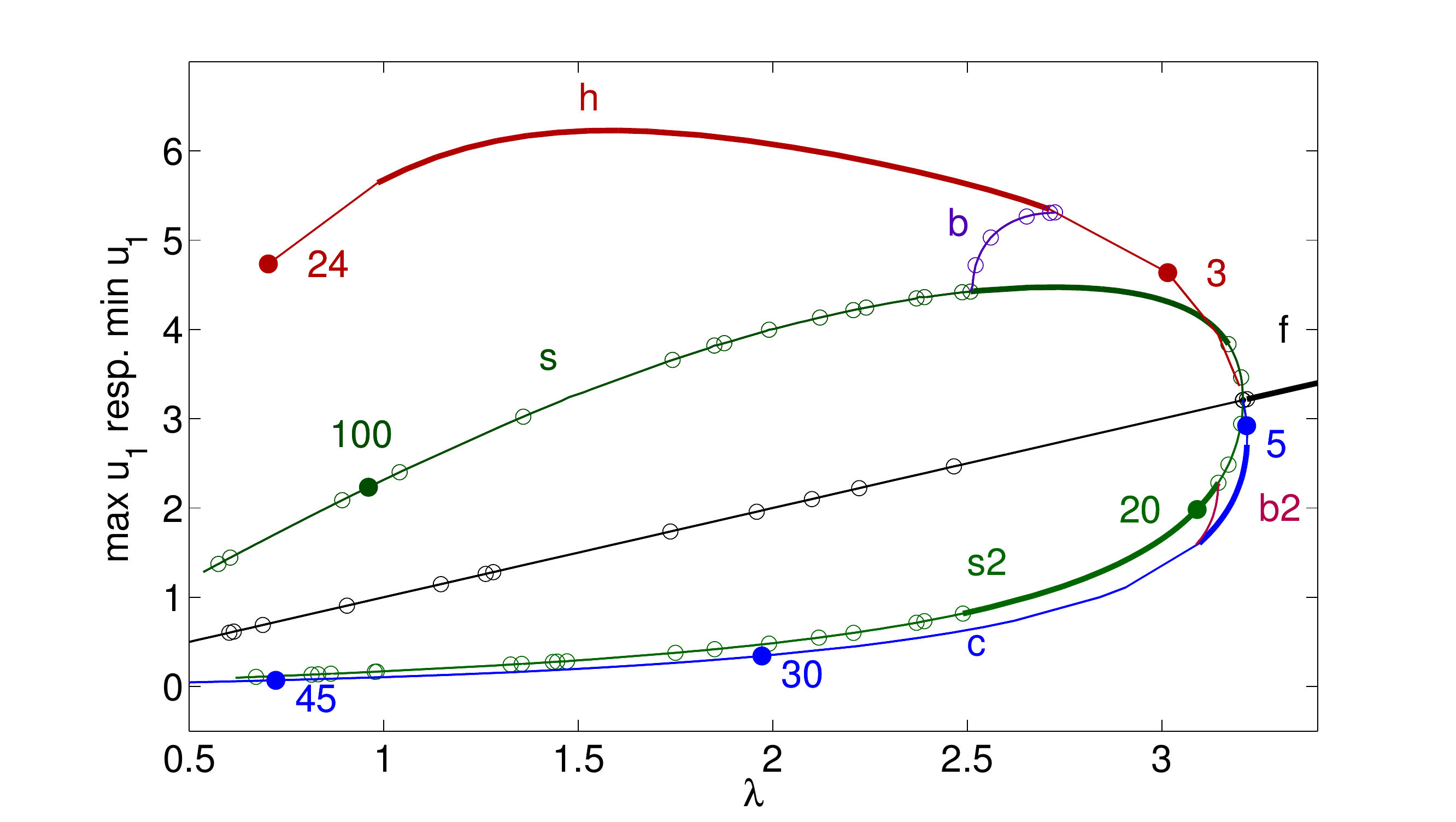}
\end{minipage}
\begin{minipage}{0.16\textwidth}
\includegraphics[width=\textwidth]{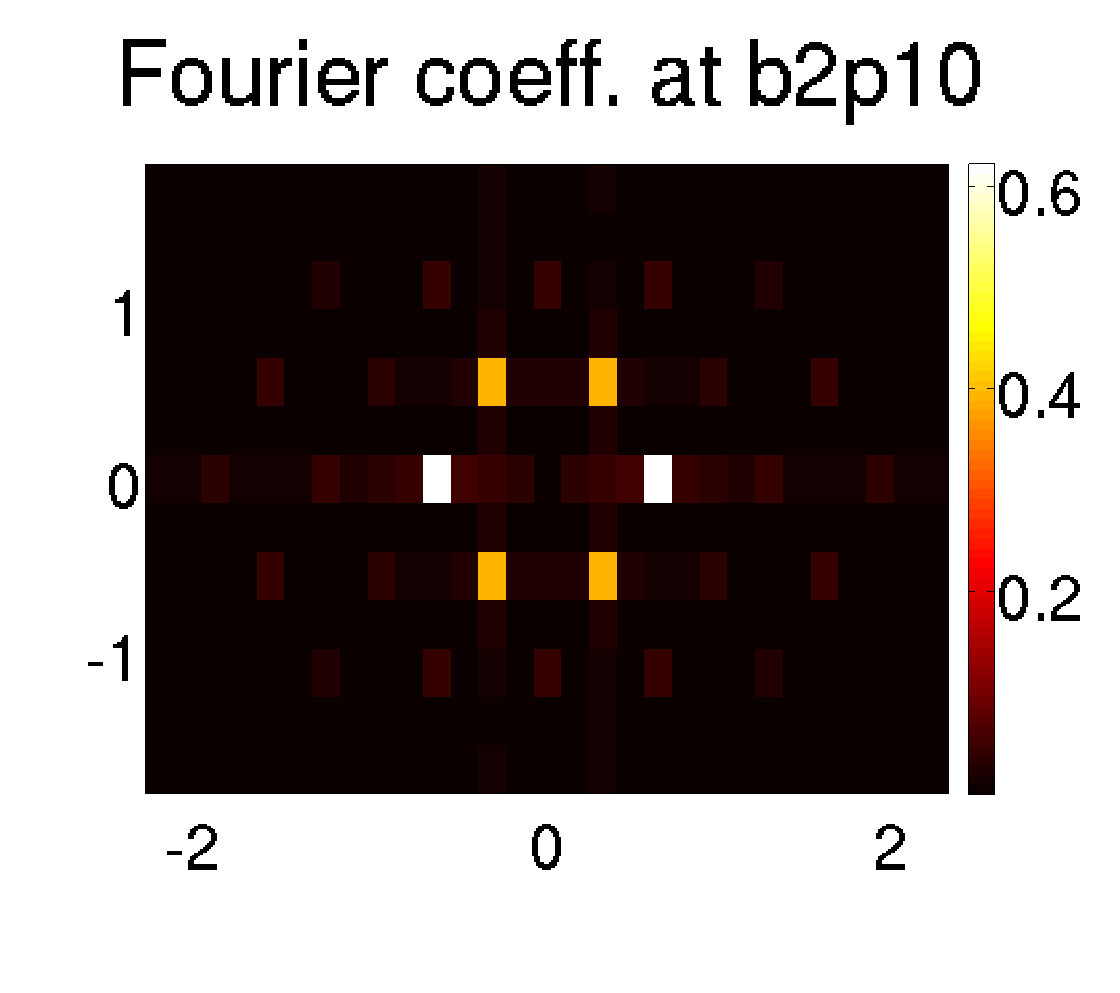}\\
\includegraphics[width=\textwidth]{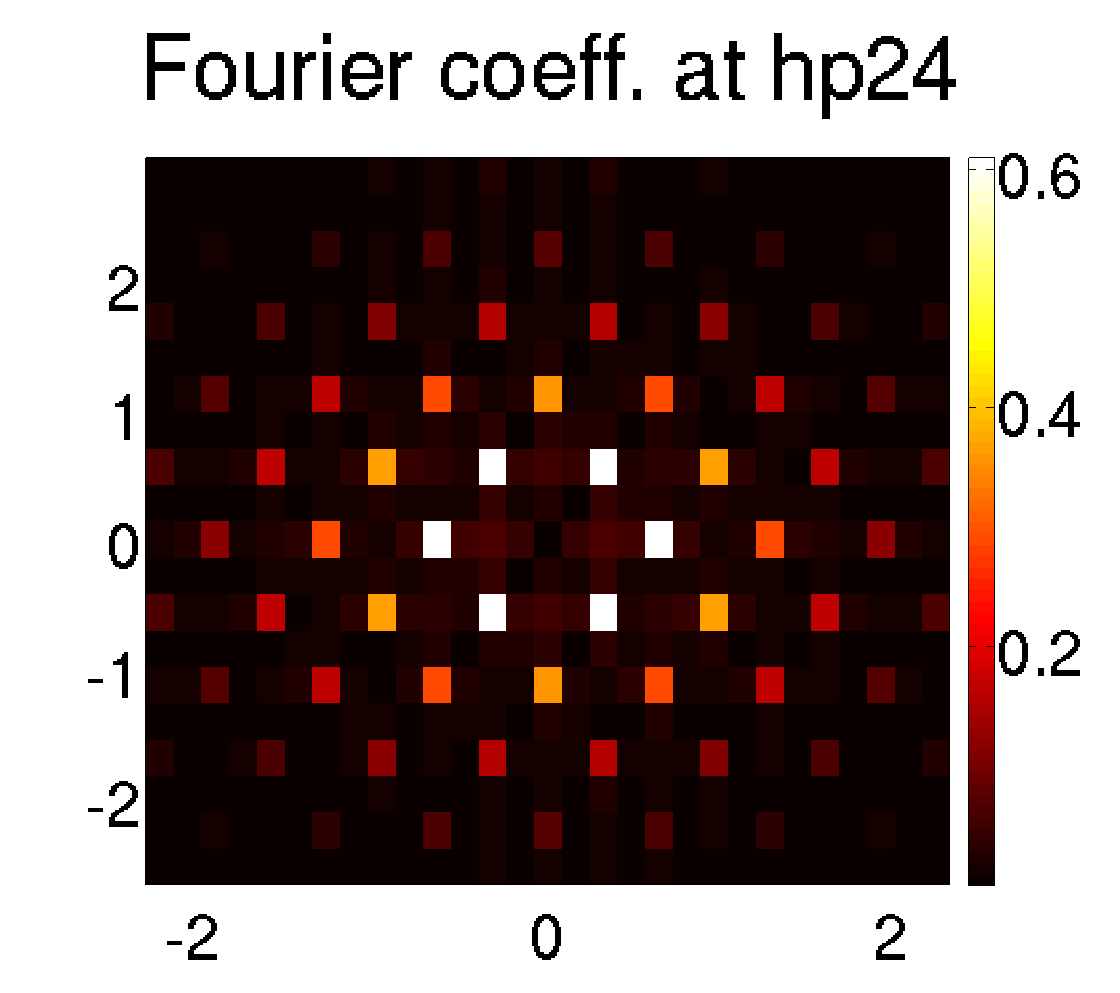}
\end{minipage}

\begin{minipage}{0.16\textwidth}
\includegraphics[width=\textwidth]{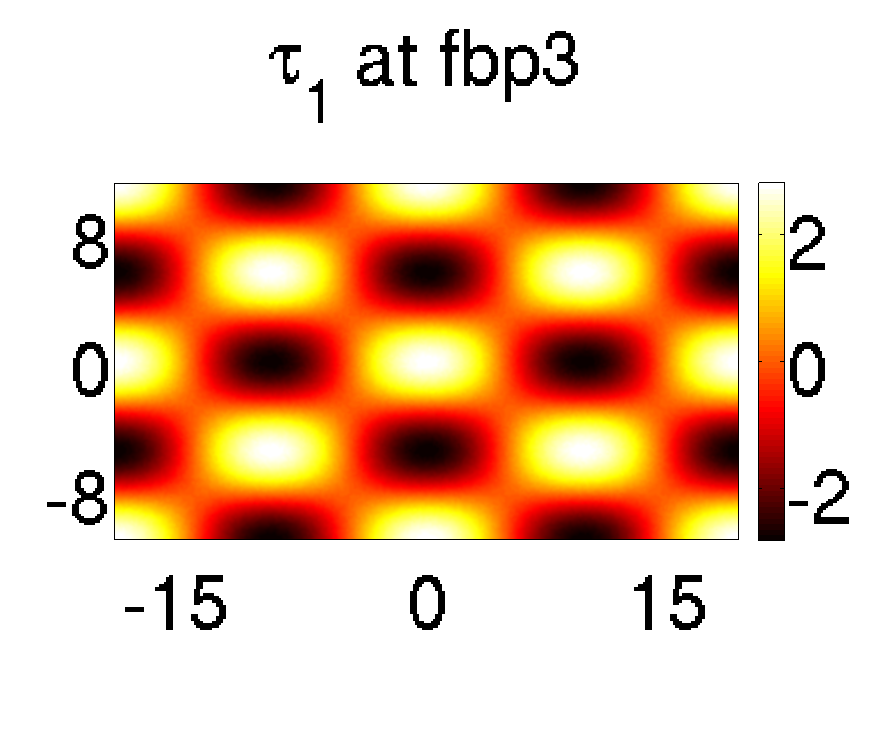}
\end{minipage}
\begin{minipage}{0.16\textwidth}
\includegraphics[width=\textwidth]{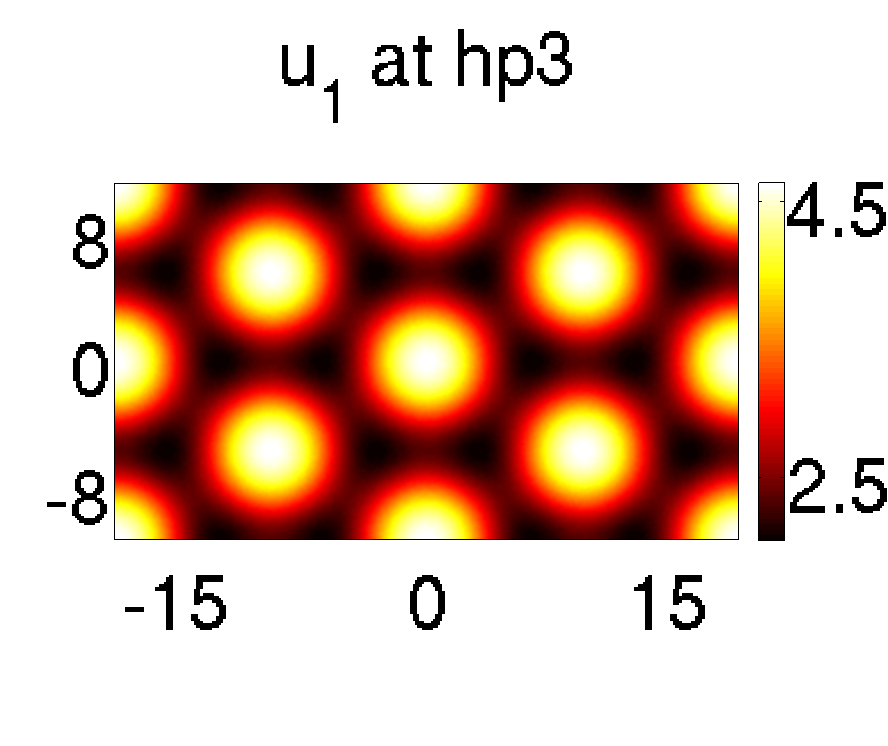}
\end{minipage}
\begin{minipage}{0.16\textwidth}
\includegraphics[width=\textwidth]{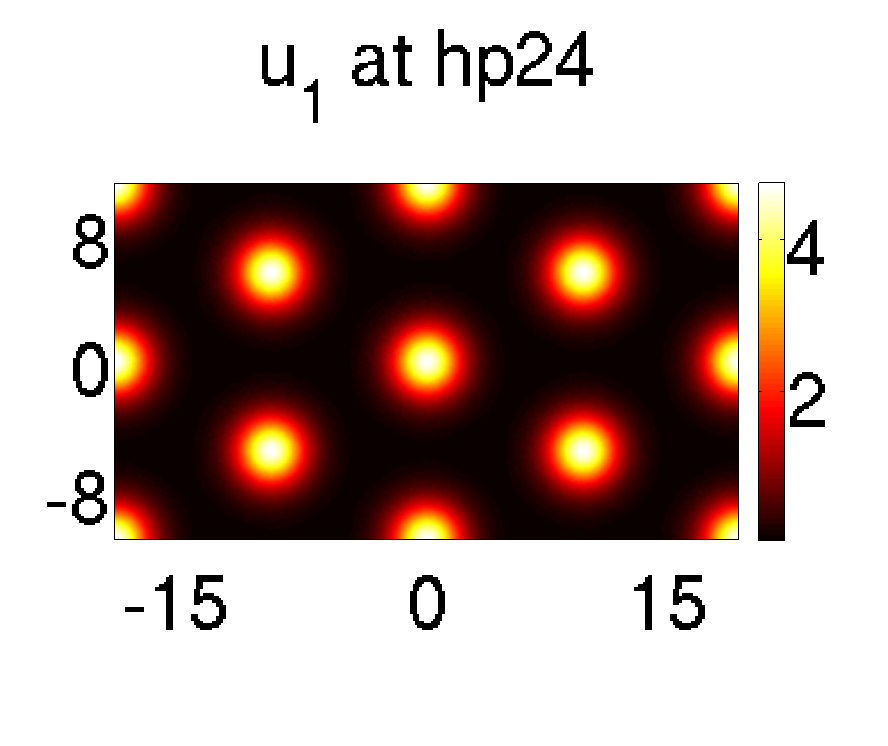}
\end{minipage}
\begin{minipage}{0.16\textwidth}
\includegraphics[width=\textwidth]{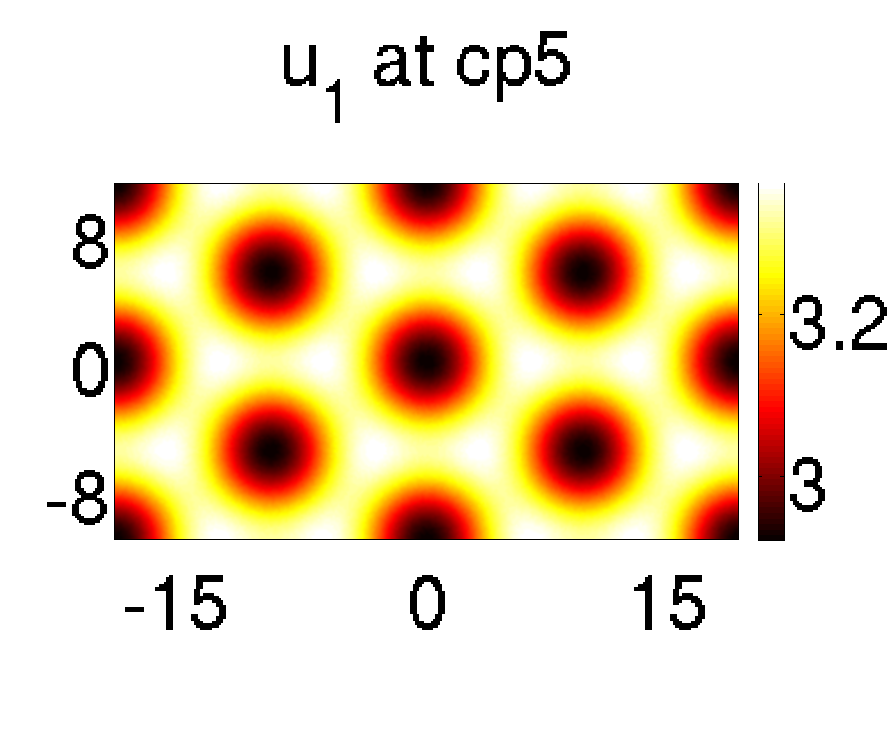}
\end{minipage}
\begin{minipage}{0.16\textwidth}
\includegraphics[width=\textwidth]{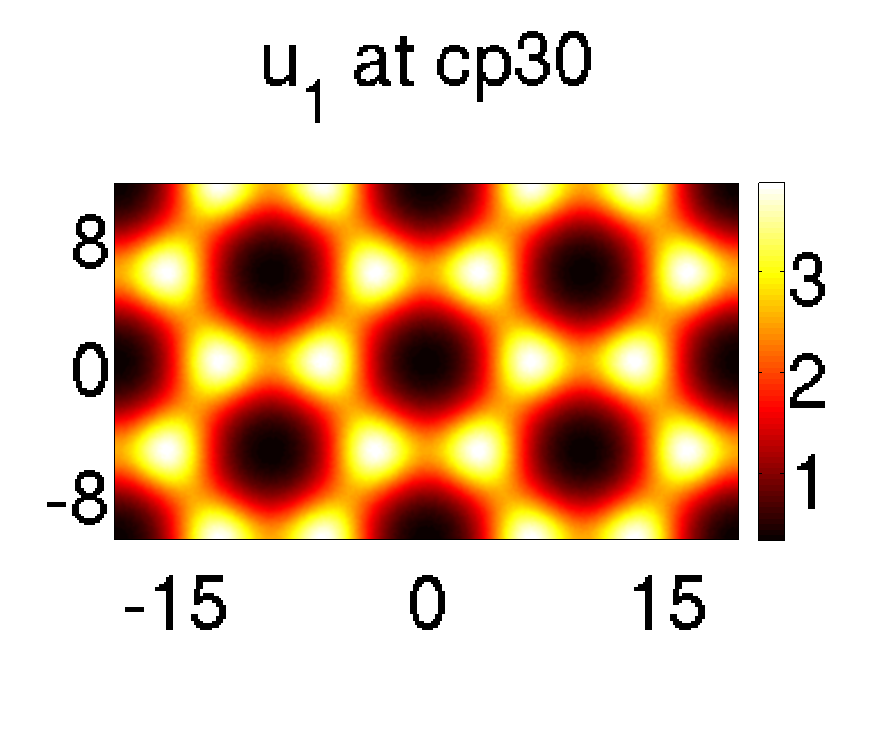}
\end{minipage}
\begin{minipage}{0.16\textwidth}
\includegraphics[width=\textwidth]{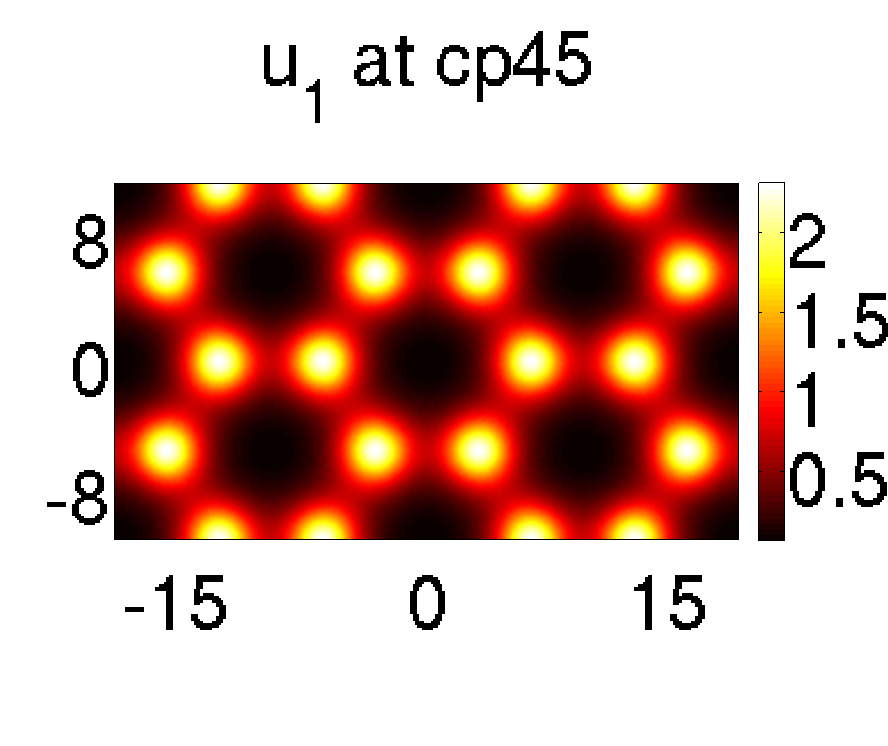}
\end{minipage}
\caption{\small{Bifurcation diagram and a selection of 
patterns for \reff{s1}. 
The branch {\tt f} in the bifurcation diagram
    represents the homogeneous solutions, {\tt s} the stripes, {\tt
      s2} the phase shift of {\tt s}, {\tt h} the hot spots, {\tt c}
    the cold spots. {\tt b} and {\tt b2} are mixed modes, also called 
    beans. $\circ$ are bifurcation points. Thick lines are stable and
    thin unstable. For {\tt s}, {\tt h} and {\tt b} we plot the
    maximum of $u_1$ and for {\tt s2}, {\tt c} and {\tt b2} the
    minimum. As usual, {\tt hp3} stands for the third point of {\tt h}, {\tt
      bbp2} for the second bifurcation point of {\tt b} and $\tau_1$
    at {\tt fbp3} for the tangent in the third bifurcation point of
    {\tt f}. On the right we also plot $|\hat{u}_{kl}|^{1/2}$, where 
$\hat{u}$ is the discrete Fourier transform of $u_1-\spr{u_1}$, see 
{\tt four.m} in directory {\tt schnakenberg}. 
These Fourier plots are often interesting for pattern forming systems: 
for instance {\tt b2p10} shows that the pattern is essentially 
still generated from the basic harmonics $\exp(\ri k_c x)^m$ and 
$\exp(\ri k_c(\frac 1 2 x\pm\frac{\sqrt{3}}{2} y)^n$, $m,n=\pm 1$, 
while at {\tt hp24} a rather large number of 
higher harmonics $\exp(\ri k_c x)^m\exp(\ri k_c(\frac 1 2 x
\pm\frac{\sqrt{3}}{2}y))^n$, $m,n\in\N$, contribute. 
Branch switching at some of the further bifurcation points 
yields a number of further interesting patterns, including some 
``snaking'' between stripes and hexagons, see \cite{uw12}. \label{schnakbifu}}}
\end{figure}

We consider the (stationary) Schnakenberg system in the form 
\huga{\label{s1} 
0=G(u)=\bpm -\Delta u_1+u_1-u_1^2u_2 \\
-d\Delta u_2 -\lambda+u_1^2u_2\epm . 
}
We use $\lambda$ as bifurcation parameter, fix $d=60$ and consider
\eqref{s1} for $(x,y)\in\Omega=[-l_x,l_x]\times [-l_y,l_y]$ with
Neumann BC.  Over $\Om=\R^2$ the spatially homogeneous solution 
$u^*(\lam)=(\lam,1/\lam)$
becomes Turing unstable \cite{Mur} when decreasing $\lam$ 
below $\lam_c\approx 3.2$, with critical wavenumber 
$k_c\approx 0.63$.  In 2D, the most
famous Turing patterns are 
\hualst{
&u(x,y)=u^*+A\cos(k_c x)+\text{h.o.t.}\qquad\text{ (stripes) }, \\
&u(x,y)=u^*+A\cos(k_c x)+
B\cos(\tfrac{k_c}{2} x)\cos(\tfrac{\sqrt{3}}{2}k_c y)+\text{h.o.t.}
\quad\text{ (hexagons, or hexagonal spots)}, 
}
where $A, B \in \R^2$ are suitable amplitudes, 
$\text{h.o.t.}$ stands for higher order terms (in $A,B,\lam-\lam_c$), and we dropped the $\lam$ dependence 
of all terms.  
Over $\Omega=[-l_x,l_x]\times [-l_y,l_y]$, if the domain size and 
BC permit it, both (spots and stripes) bifurcate from the trivial
branch at $\lam=\lam_c$. Here, to make $\lam_c$ a simple bifurcation
point (for vertical stripes), we choose $l_x=2m\pi/k_c$ and 
$l_y=2n \del \pi/(\sqrt{3}k_c)$, $m,n\in N$, where 
$\del\approx 1$ is a deformation parameter, such that for 
$\del\ne 1$ the double bifurcation point splits into 
two simple bifurcation points. To locate these, we start 
on the homogeneous branch and use a bisection type routine  
based on the number of negative eigenvalues of $G_u$, 
see {\tt findbif.m}. 
Figure \ref{schnakbifu} shows a bifurcation diagram and some 
solution plots obtained for $m=n=2$ and $\delta=0.99$. 

Here the problem is that using the standard 
{\tt cont} algorithm we quickly obtain some undesired branch switching. 
For instance when continuing the stripe branch {\tt s} with 
standard settings we switch to the beans branch {\tt b} when 
approaching its bifurcation point. This particular 
branch switching can be avoided by decreasing $\xi$ to 
$\xi=0.1/{\tt p.np}$, say, but only to the effect that 
we get branch--switching at some later point on the {\tt s} branch. 
Such undesired branch--switching is a serious problem 
in all continuation algorithms, see, e.g., \cite[\S3]{seydel}, 
and we use a modification {\tt pmcont} of {\tt cont} explained 
in the next section, which also incorporates some parallel computing 
for speedup.
\footnote{Still, the bifurcation diagram in Fig.~\ref{schnakbifu} is 
computationally quite expensive (about 40 minutes on a quad-core desktop PC, 
with about 60.000 triangles on average), 
and therefore the init-function {\tt p=schnakinit(p,m,n,nx,del)} 
takes the domain sizes $m,n$, the deformation parameter $\delta$ 
and the startup spatial discretization {\tt nx} as parameters. 
{\tt schnak11demo.m} (or schnak11cmds.m)
then uses $m=n=1$ and $\del=0.97$ and only 
takes a few minutes for a bifurcation diagram similar to 
Fig.~\ref{schnakbifu} over the smaller domain, while 
{\tt schnak22cmds.m} treats the $m=n=2$ domain in Fig.~\ref{schnakbifu}.}

\subsection{{\tt pmcont}}\label{pmcont-sec}
\begin{figure}[htbp]
\bce 
\begin{minipage}{0.25\textwidth}
\ig[width=42mm,height=37mm]{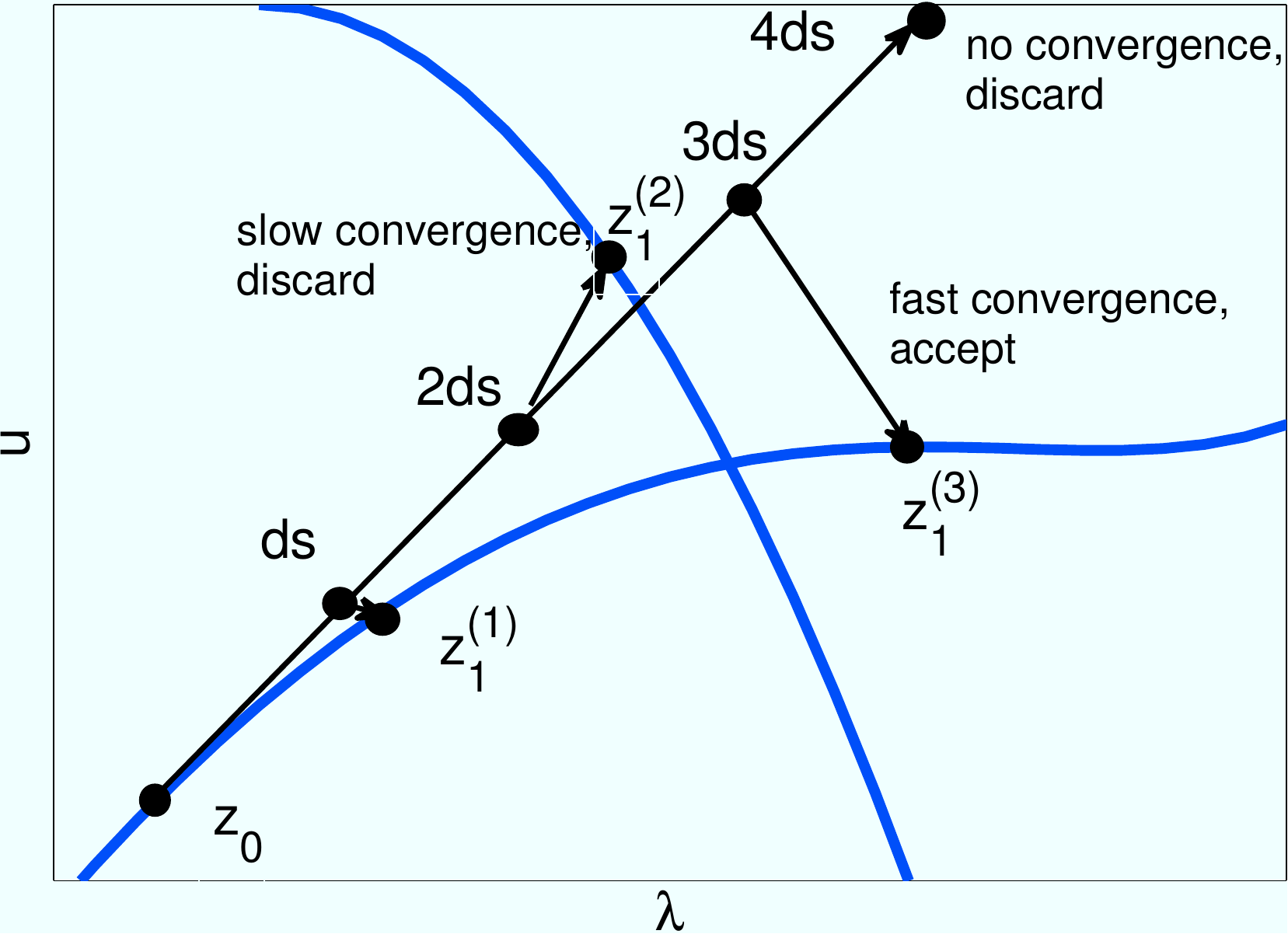}
\end{minipage}\ \ 
\begin{minipage}{0.69\textwidth}
\fbox{\parbox{\textwidth}{\small 
{\bf Algorithm {\tt pmcont}}\\[-7mm]
\ben
\item Multipredictors. $(u^i,\lambda^i)=(u_0,\lambda_0)+i\ {\tt p.ds}\,\tau$, 
$i=1,\dots,{\tt p.mst}$   
\item Newton--loops (parallel). Use \reff{resdec} to identify ``good points''.
\item Tangents (sequentially). Calculate new tangents $\tau_1,\ldots,\tau_m$ at good points, using \reff{tau1}.
\item Bifurcation detection and localization (parallel). 
\item Postprocess (sequentially). Call {\tt ufu}, save, plot, return to 1.
\een 
}}
\end{minipage}
\ece

\vs{-2mm}
\caption{{\small  Sketch of the basic idea of multiple predictors 
and convergence monitoring, and pseudocode of {\tt pmcont}. 
The arrows from, e.g. ``2ds'' to $z_1^{(2)}$ 
just illustrate the result of the Newton loops, 
not the hyperplanes $\{(u,\lam)\in\R^{Nn_p+1}: \spr{\tau_0,(u,\lam)}_\xi={\rm ds}\}$ 
as in Fig.\ref{xif}. \label{pmf}}}
\end{figure}

Theorem 4.4 in \cite{keller77} guarantees that the standard continuation 
converges 
to a given branch for ``suffienctly small'' $\rds$, but near 
bifurcation points only 
in cones around the branch. Thus, near a bifurcation point 
it is often not useful to choose very small $\rds$.
To circumvent this and similar
problems we provide the function {\tt pmcont}. 
The basic idea is explained in Fig.~\ref{pmf}. 

Instead of using just one predictor $(u^1,\lambda^1)=(u_0,\lambda_0)+{\tt p.ds}\,\tau$,  
{\tt pmcont} creates in every continuation step the predictors
$$
(u^i,\lambda^i)=(u_0,\lambda_0)+i\,{\tt p.ds}\,\tau, \quad 
i=1,\dots,{\tt p.mst}, 
$$
and starts a Newton loop for each. 
{\tt pmcont} then monitors 
the convergence behaviour of each loop to decide whether it yields 
a ``good'' point, i.e., a point on the present branch. 
The criterion is that in each Newton step the residual 
has to decrease by a factor $0<\al<1$, i.e., 
\huga{\label{resdec} 
\|G(u_{n+1},\lam_{n+1})\|\le \al\|G(u_n,\lam_n)\|, 
}
otherwise the loop is stopped. The heuristic idea is that if 
the Newton loop converges slowly, then probably the 
solution is on a different branch, because the loop has to (slowly) change 
the solution ``shape''. 

For the crucial parameter $\al$, which describes the desired convergence speed, 
 we recommend trying $\al=0.1$. \footnote{However, for instance on the 
``hot hexagon branch'' {\tt h} in Fig.~\ref{schnakbifu} we need to use $\al=10^{-6}$ 
to avoid branch--switching, see {\tt schnak22cmds.m}.} 
Of course, these heuristics 
in no way guarantee that no branch switching occurs, or anyway that we get 
convergence for long predictors 
$(u^i,\lambda^i)=(u_0,\lambda_0)+i \cdot{\tt ds}\tau$ 
with $i>1$. But in practice 
we find the idea to work remarkably well. See also 
\S\ref{s:rbconv} for an example of the ``unreasonable effectiveness'' 
of {\tt pmcont}, together with an example that long predictors 
in {\tt pmcont} tend to branch--switching in imperfect bifurcations. 

Concerning data structures, we need three additional parameters (with {\tt p} the problem structure): the number of predictors {\tt p.mst}, 
$\al=${\tt p.resfac}, and {\tt p.pmimax}. These are used to gain 
some flexibility for \reff{resdec} via stepsize control. 
If {\tt p.mst} equals the number of
solutions found and {\tt p.ds} is smaller than {\tt
  p.dsmax/p.dsincfac}, then {\tt p.ds} will be increased by the factor
{\tt p.dsincfac}, essentially as in {\tt cont}. 
If no solution is found and {\tt p.ds} is greater than {\tt
  (1+p.mst)}$\cdot${\tt p.dsmin}, then the step size {\tt p.ds} will
be divided by {\tt 1+p.mst} in the next continuation step. 
Finally, if {\tt
  p.ds} is less than {\tt (1+p.mst)}$\cdot${\tt p.dsmin} and
{\tt p.pmimax} is less than {\tt p.imax}, then {\tt p.pmimax} will be
increased to {\tt p.pmimax+1}, and 
$\|G(u_{n+{\tt p.pmimax}},\lam_{n+{\tt p.pmimax}})\|\le {\tt p.resfac}\|G(u_n,\lam_n)\|$ 
is required instead of \reff{resdec}.   
The only new functions used in {\tt pmcont} are {\tt pmnewtonloop.m} and 
{\tt pmbifdetec.m}. 

Another advantage of the {\tt p.mst} predictors is that the 
Newton loops can be calculated in parallel, which 
on suitable machines gives substantial 
speedups.\footnote{Here we use the Matlab Parallel Computing Toolbox 
in an elementary setup 
with basic monitoring of open kernel threads. 
} 
All final Newton iterates with a residual smaller than {\tt p.res} are
taken as solutions, and plotted and saved as in {\tt cont}. Next, the tangents 
$\tau_1,\ldots,\tau_m$ are calculated sequentially, because $\tau_{l+1}$ needs 
$\tau_{l}$, and afterwards the bifurcation detection and localization is again 
in parallel. 
The last solution and its tangent will be used for the next 
continuation step. Mesh adaption/refinement is inquired at the 
start of {\tt pmcont}, i.e., before generating the predictors, 
but not on the individual correctors. 

In summary, for {\tt p.resfac=1} and {\tt p.mst=1} we have that {\tt pmcont} is 
roughly equivalent to {\tt cont}, except for slightly less versatile 
mesh adaption and error estimates. For {\tt p.mst}$>1$, {\tt pmcont} 
 takes advantage of parallel computing, and  
is often useful to avoid convergence problems and undesired branch 
switching close to bifurcation points. 
The main reasons why we 
(currently) keep the two version and do not combine them into one is 
that {\tt cont} is simpler to hack and 
implements Keller's basic, well tested algorithm.

\section{Three classical examples from physics}\label{s:phys}
In this section we consider models for Bose-Einstein (vector) solitons, 
Rayleigh-B\'enard convection, and the \vKar\ system for buckling of an 
elastic plate, as examples for systems with more than two
components, and with \bcs\ implemented via {\tt gnbc} as described in \S\ref{geosec}. 
The largest and most complicated system (in the sense of number of
components and implementation of \bcs) here is the \vKar\ system. 
Hence, for this we also explain the coding in \pdep\ in most detail, while for 
 Bose-Einstein solitons and Rayleigh-B\'enard convection we mostly refer
to the m-files for comments.

\subsection{Bose--Einstein (vector) solitons {\tt (gpsol)}} 
\label{gpsec}
As an example with $x,y$ dependent coefficients, nontrivial advection, 
and interesting localized solutions we consider 
(systems of) Gross--Pitaevskii (GP) equations 
with a parabolic 
potential that arise for instance as amplitude equations in  Bose--Einstein 
condensates.
\subsubsection{The scalar case} 
First, following \cite{lash08} we consider the scalar equation 
\huga{\label{nls1} 
\ri\pa_t\psi=-\Delta\psi+r^2\psi-\sigma|\psi|^2\psi, 
}
where $\psi=\psi(x,y,t)\in\C$, $r^2=x^2+y^2$, and $\sigma=1$ 
(focussing case). This has a huge number of families of localized 
solutions, aka solitons, which may be time periodic, standing or 
rotating in space. 
Going into a frame rotating with speed $\om$ and splitting off 
harmonic oscillations with frequency $\mu$, i.e., 
\huga{\label{nls1a}
\psi(x,y,t)=\Phi(r,\phi-\om t)\er^{-\ri \mu t}, 
}
we obtain 
\huga{\label{nls2}
\bigl[\pa_r^2+\frac 1 r\pa_r+\frac 1 {r^2}\pa_\theta^2
-\ri\om\pa_\theta+\mu-r^2\bigr]\Phi
+\sigma|\Phi|^2\Phi=0.
}
A typical ansatz for (approximate) solutions has the form 
\huga{\label{nlsa} 
\Phi(r,\theta)=A\phi(r/a)(\cos(m\th)+\ri p\sin(m\th)), 
\quad \phi\in\R, \text{ e.g. }
\phi(\rho)=\rho^mL_n^{(m)}(\rho^2)\er^{-\rho^2/2}, 
}
with $L_n^{(m)}$ the $n^{th}$ Laguerre polynomial. 
Plugging this into \reff{nls2} yields expressions for $A,a,p,\om$ 
for approximate solutions. 
The case $n{=}0$ and $p{=}0$ corresponds to so called 
(nonrotating, since $\om=0$) real $m$--poles, 
and $|p|=1$ to a so called radially symmetric vortex of charge $m$, 
while the intermediate cases $0<|p|<1$ give to so called 
rotating azimuthons with interesting angular modulations 
of $|\Phi|$. 

\begin{figure}[ht!]
\bce 
{\small 
\begin{tabular}[t]{p{36mm}p{38mm}p{37mm}p{38mm}} 
(a) \ig[width=36mm,height=37mm]{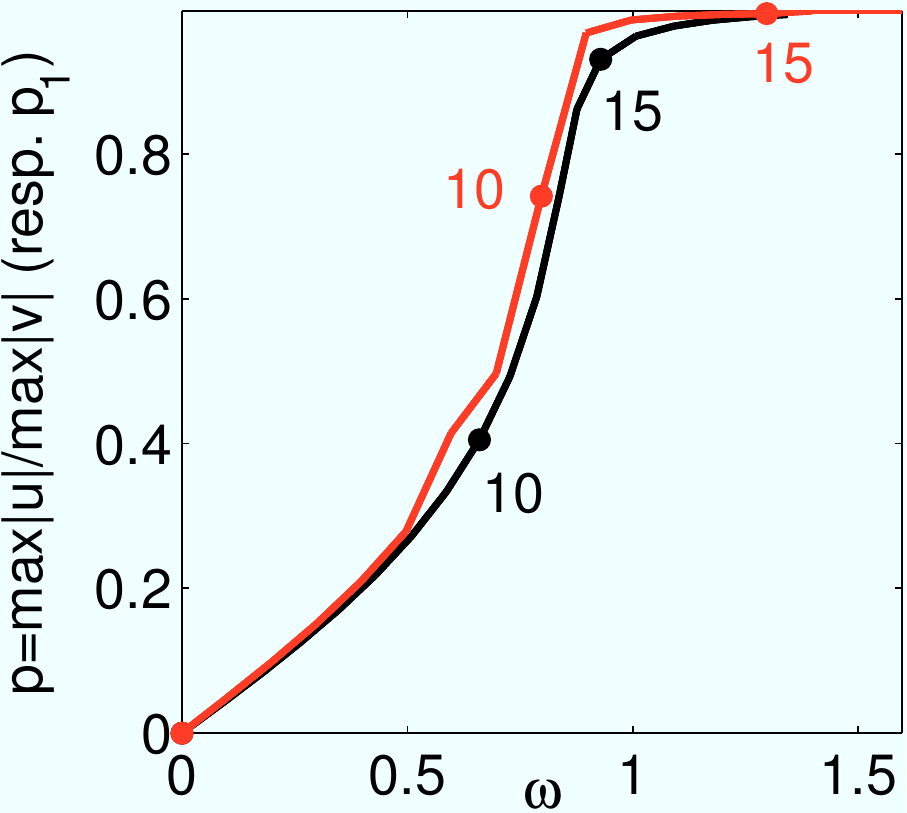}&
(b) \ig[width=35mm,height=39mm]{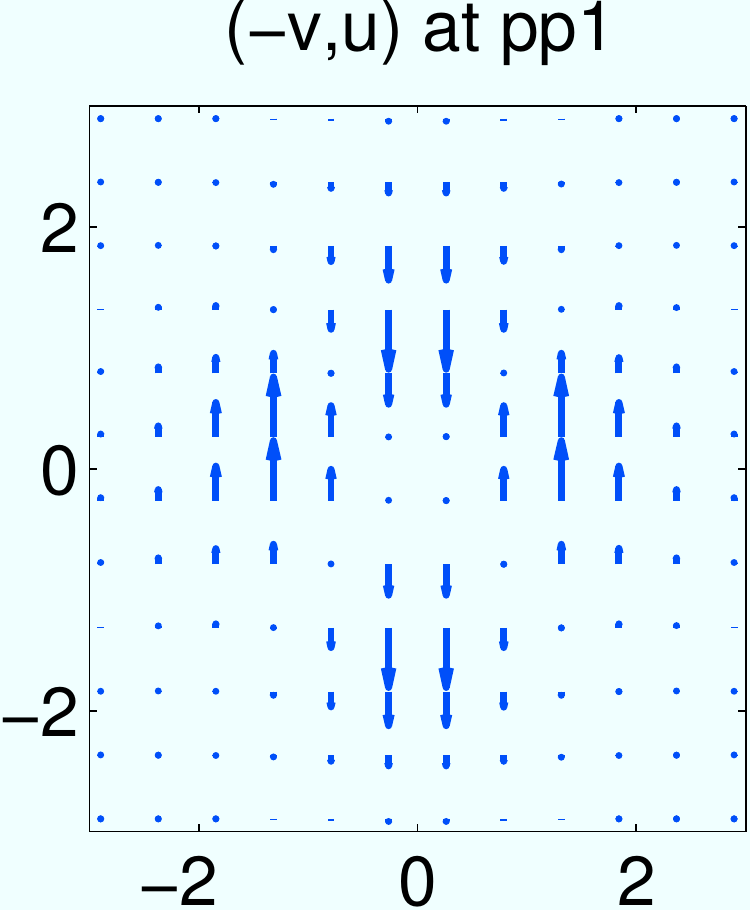}& 
(c) \ig[width=35mm,height=39mm]{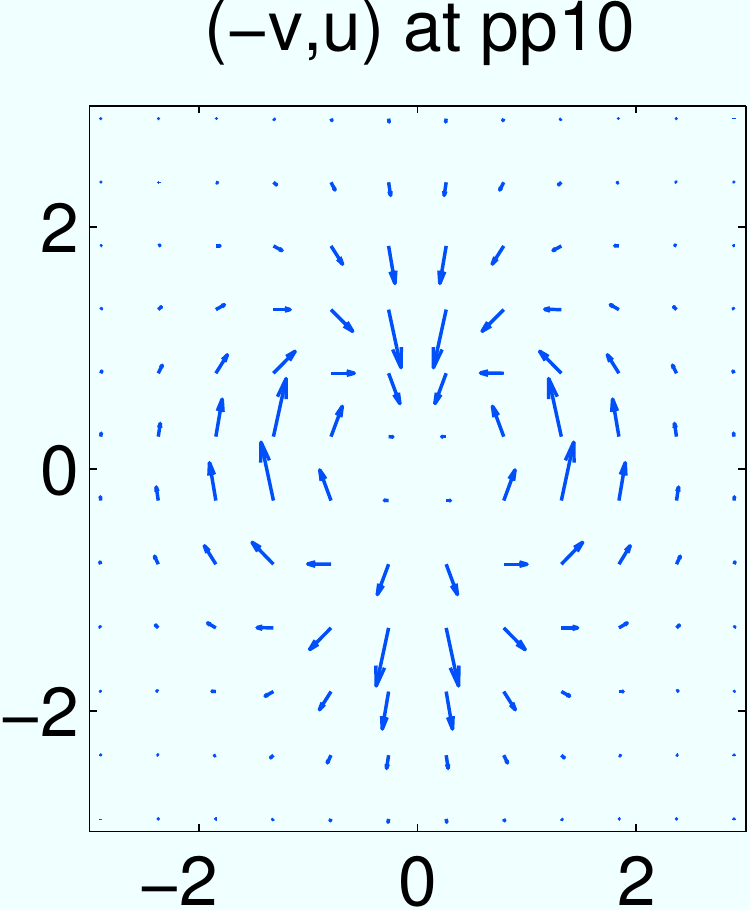}&
(d) \ig[width=38mm,height=42mm]{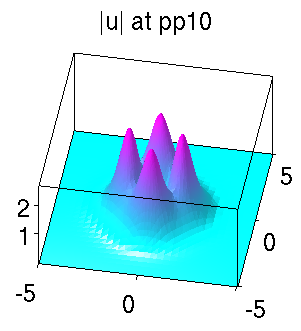}\\
(e) \ig[width=37mm,height=42mm]{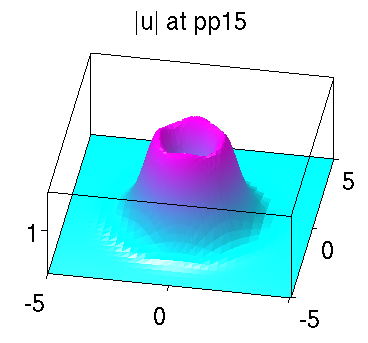}&
(f) \ig[width=38mm,height=35mm]{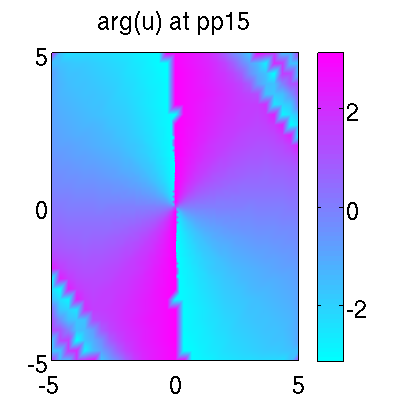}&
(g) \ig[width=34mm,height=40mm]{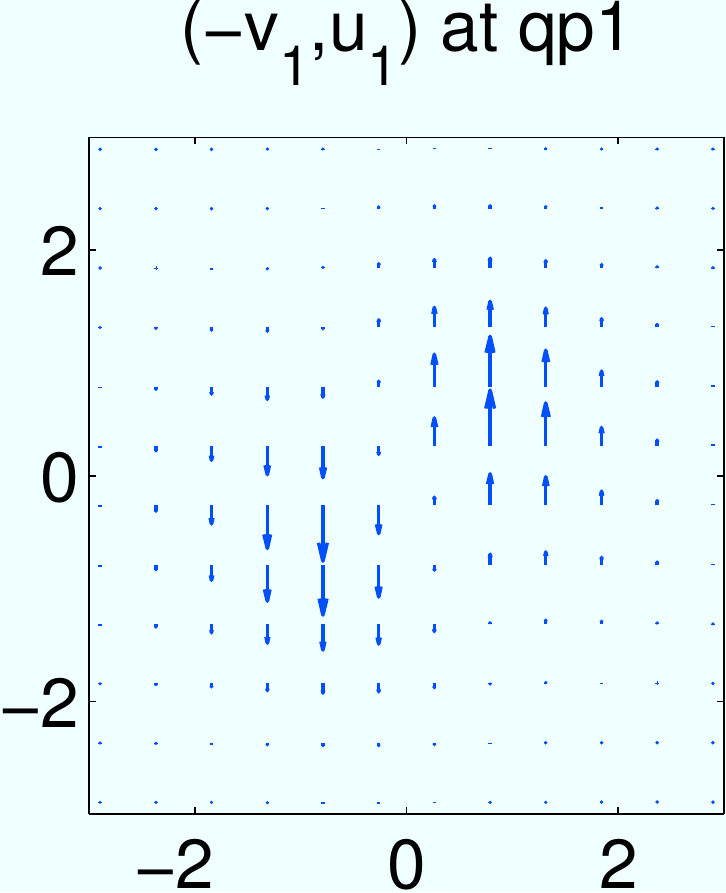}&
(h) \ig[width=34mm,height=40mm]{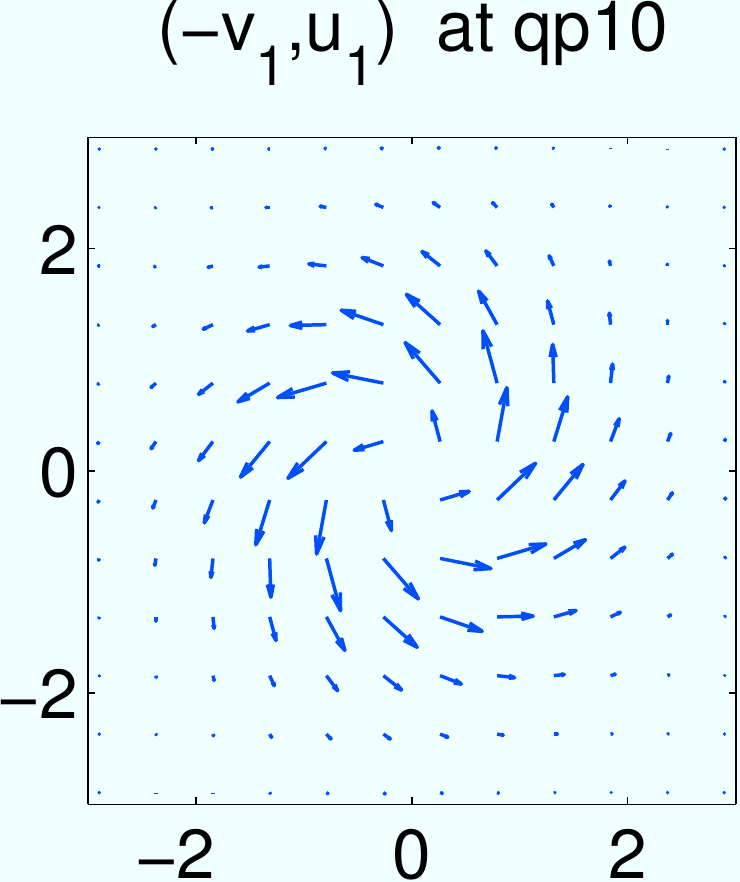}
\end{tabular}}
\ece

\vs{-4mm}
\caption{\label{bef1}{\small (a) Continuing a quadrupole to an azimuthon to 
a 2-vortex, $p=\max|u|/\max|v|$ over $\lam=\om$, ({\tt p} branch, black); 
and continuing a vector dipole  to a vector--azimuthon to a vector--vortex, 
$p_1=\max|v_1|/\max|u_1|$ over $\lam=\om$ ({\tt q} branch, red), $\mu=2$ 
resp.~$\mu_1=2, \mu_2=2.2$. 
(b),(c) vectorfield plots 
at $\lam=0$ (quadrupole) resp.~$\lam\approx 0.66$ (azimuthon). (d)--(f) 
$|U|$ and arg$(u+\ri v)$ as indicated. (g),(h) vectorfield plots for the 
first condensate $\psi_1$ as indicated, the second condensate is similar.  
Also see the customized {\tt plotsol.m} in {\tt gp}. The phase plot in 
(f) is somewhat ragged due to the coarse mesh away from the center. 
For the {\tt p} branch we used fixed {\tt p.nt=5208} from mesh-refinement 
during init. In the {\tt q} branch we used {\tt q.amod=8} with {\tt q.maxt=9000} 
which lead to {\tt q.nt=9500} at, e.g., point 10.  }} 
\end{figure}

Our goal is to calculate these solutions numerically with \pdep. 
Returning to cartesian coordinates,  i.e., setting 
$\Phi(r,\theta)=u(x,y)+\ri v(x,y)$ 
we obtain the 2-component real elliptic system 
\begin{subequations}\label{nls3}
\hual{
&-\Delta u+(r^2-\mu)u-|U|^2 u-\om(x\pa_y v-y\pa_x v)=0, \\
&-\Delta v+(r^2-\mu)v-|U|^2 v-\om(y\pa_x u-x\pa_y u)=0, 
}
\end{subequations} 
where $|U|^2=u^2+v^2$. Our strategy is to use \reff{nlsa} 
for $\om=0$ and to continue in $\lam:=\om$. 
A measure for the deformation of multipoles into vortices
for the numerical solutions is the 
``modulation depth'' $p$ of the soliton intensity 
\huga{\label{nls5}
p=\max|\im\Phi|/\max|\re\Phi|=\max|v|/\max|u|.
}
The result of typical continuation of a quadrupole using a stiff--spring 
approximation of DBC for $u,v$ on domain $\Om=[-5,5]^2$ 
is shown in Fig.\ref{bef1}, (a)-(f), 
see {\tt gpf.m}, {\tt gpjac.m}, {\tt gpcmds.m} and {\tt gpinit.m}, 
and also {\tt plotsol.m} in directory {\tt gpsol} for the customized
of {\tt plotsol}. 

The following remarks are in order. First, we switch off 
stability or bifurcation tracking ({\tt spcalcsw=0}, {\tt bifchecksw=0}), 
see Remark \ref{nlsrem1} below. 
Second, the linearization of \reff{nls3} is given by 
(recall that $\lam=\om$) 
\hualst{
G_u(u,v)=&\bpm -\Delta&0\\0&-\Delta\epm+\bpm-\mu+r^2-3u^2-v^2&-2uv\\
2uv&-\mu+r^2-3v^2-uv\epm \\&
-\lam\bpm 0& x\pa_y-y\pa_x\\-x\pa_y+y\pa_x&0\epm. 
}
Thus, 
the last term is a good example how to use {\tt assemadv} with a 
relatively complicated $b$. 
Third, some problems should be expected from the large number of solutions
of \reff{nls2}, in particular the phase--invariance: 
if $\Phi$ is a solution, so is $\Phi^{\ri \al}$ 
for any $\al$, or equivalently, \reff{nls3} is invariant under 
multiplication with ${\small \bpm \cos \al&-\sin\al\\\sin\al&\cos\al\epm}$. 
Thus, even for all parameters fixed, solutions of \reff{nls3} always 
come in continuous families, and thus $G_u$ as a linear operator 
on $[L^2(\R^2)]^2$, say, always 
has a zero eigenvalue. See also \cite{lash08} and the references 
therein for some tricks for the numerical solution of \reff{nls2}. 
Rather remarkably, we need none of these tricks, presumably since 
numerically in $G_u$ 
the zero eigenvalue is perturbed sufficiently far away from $0$. 
One trick we do use is to start with a coarse mesh of $30\times 30$ points, 
first take some (rather arbitrary) 
monopole as dummy--starting guess, use {\tt meshref} to generate 
a rather fine mesh in the center, define 
a quadrupole initial guess using \reff{nlsa} on that first refined mesh, 
and then refine again, yielding the (still small) number of 
5208 triangles for this continuation, with an error--estimate 
less than 0.01. On a small laptop computer 
the whole continuation takes about a minute.

\subsubsection{A two--component condensate}
The above can be generalized to multi--component condensates \cite{lash-kiv09}. 
For two components, we then have coupled GP equations of the form, e.g., 
\huga{\label{nlss} 
\ri\pa_t\psi_1=[-\Delta+r^2-\sigma|\psi_1|^2-g_{12}|\psi_2|^2]\psi_1, 
\quad  \ri\pa_t\psi_2=[-\Delta+r^2-\sigma|\psi_2|^2-g_{21}|\psi_1|^2]\psi_2, 
}
where $g_{12},g_{21}$ are called interspecies interaction coefficients. 
Physically, it makes sense to use ans\"atze of the form \reff{nls1a} 
with different $\mu$ but equal $\om$, i.e., 
$\psi_j(x,y,t)=\Phi_j(r,\phi-\om t)\er^{-\ri \mu_j t}$.  
Next we can use the form \reff{nlsa} for each component 
$\Phi_j$ and classify the thus obtained approximate solutions 
as soliton-soliton, soliton--vortex, soliton--azimuthon etc pairs. 
To calculate such solutions numerically we set 
$
\Phi_j(r,\th)=u_j(x,y)+v_j(x,y)
$
and obtain an elliptic system of the form \reff{nls3} but 
with four real equations. This has been implemented 
in {\tt vgpf}, with Jacobian {\tt vgpjac}. Figure 
\ref{bef1} (a),(g),(h) shows the continuation 
of a two--dipole obtained from {\tt vgpcmds}. Similar remarks as for the scalar 
case apply, see the comments in {\tt vgpcmds.m}.

\brem\label{nlsrem1} 
This clearly was just a very introductory demo of continuation 
of solutions of \reff{nls1} resp.~\reff{nlss}; there are {\em many} 
more and interesting branches, and further questions, again see, e.g., 
\cite{lash08,lash-kiv09}. 
Interesting questions concern, e.g., the dependence of 
vector solitons on $|\mu_1-\mu_2|$ which can for instance be studied 
by fixing $\om$ and continuing in $\lam=\mu_2$, or the effect of 
including a periodic potential, leading to gap solitons 
\cite{DU09}. 
Some of these questions will be considered elsewhere. 
\eex
\erem 

\subsection{Rayleigh-B\'enard convection 
({\tt rbconv})}\label{s:rbconv}
As an example from fluid dynamics we consider two-dimensional
Rayleigh-B\'enard convection in the Boussinesq approximation in the
domain $\Omega=[-2,2]\times [-0.5,0.5]$. In the
streamfunction formulation the stationary system reads 
\begin{align}\label{e:rbconv}
  -\Delta\psi  + \omega&= 0, \nonumber\\
  -\sigma\Delta \omega -\sigma R \partial_x\theta + \partial_x\psi\partial_z\omega - \partial_z\psi\partial_x\omega &=0, \\
  -\Delta \theta - \partial_x\psi + \partial_x\psi\partial_z\theta
  - \partial_z\psi\partial_x\theta &=0,\nonumber
\end{align}
with streamfunction $\psi$, temperature $\theta$, and the auxiliary
$\omega=\Delta \psi$. Moreover, $\sigma$ is the Prandtl number, set to
$1$ here, and $R$ the Rayleigh number, which will be the continuation
parameter. The implementation of \reff{e:rbconv} in \pdep\ is 
relatively straightforward, including analytical Jacobians, see 
{\tt rbconvf.m} and {\tt rbconvjac.m} 

The boundary conditions at the top and bottom plates are taken at
constant temperature and with zero tangential stress
\[
\psi=\partial_{zz}\psi=\theta=0, \mbox{ at } z=\pm 0.5;
\]
note that this means $\Delta\psi=0$ at $z=\pm 0.5$. 
Motivated by the analysis in \cite{HirschKnob97}, laterally we
consider on the one hand ``no-slip'' (and perfectly insulating) \bcs 
\begin{align}\label{e:bcnoslip}
  \psi=\partial_x\psi=\partial_x\theta=0, \mbox{ at } x=\pm L,
\end{align}
and on the other hand ``stress free'' \bcs 
\begin{align}\label{e:bcstressfree}
  \psi=\partial_{xx}\psi=\partial_x\theta=0, \mbox{ at } x=\pm L.
\end{align}
See the comments in {\tt rbconvbc\_noslip.m} resp.~{\tt
  rbconvbc\_stressfree.m} for the implementation (approximation) of
these \bcs\ based on \reff{e:gnbc}.

In both cases it is known that continuation of the trivial zero state
for increasing $R$ gives a sequence of bifurcations alternating
between even and odd modes. The stability thresholds are plotted in
\cite{HirschKnob97}, Figure 1(a) for \eqref{e:bcstressfree} and 1(b)
for \eqref{e:bcnoslip}. We use these to choose initial values of
$\lambda=R$ for the first two bifurcations, respectively.

\begin{figure}[ht!]
\bce 
{\small 
\begin{minipage}{40mm}
(a) Bifurcation diagram\\ 
\ig[width=40mm,height=34mm]{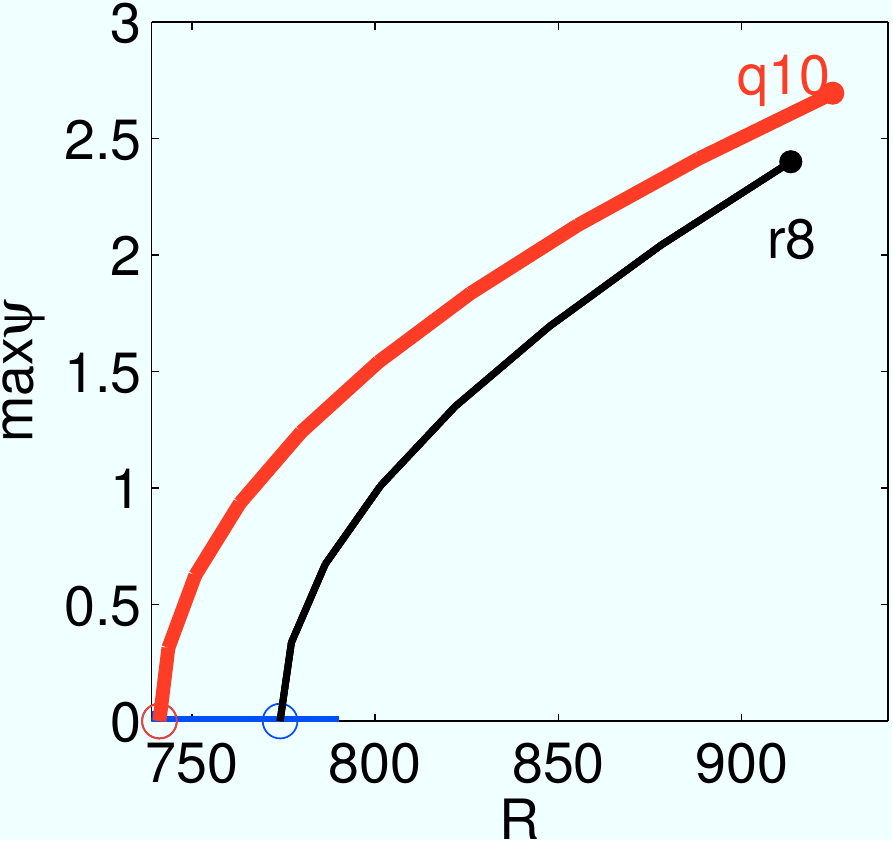}
\end{minipage}
\begin{minipage}{60mm}
(b) q10 (top) and r8 (bottom)\\ 
\ig[width=60mm,height=17mm]{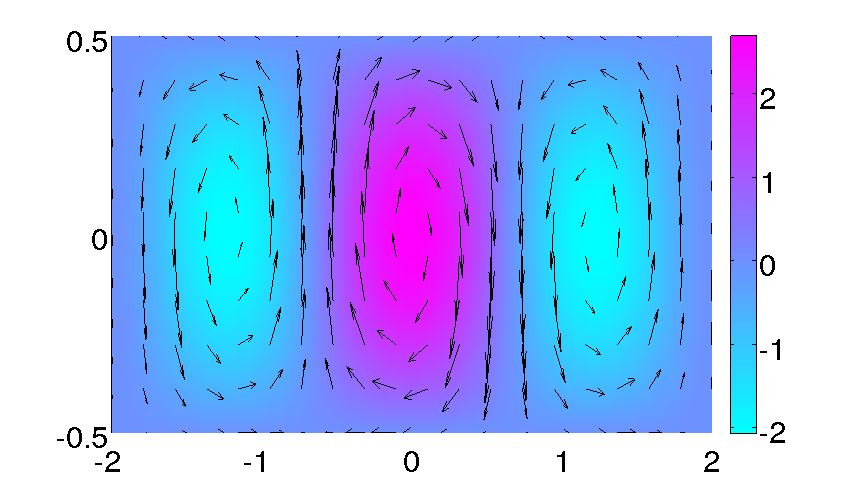}
\ig[width=60mm,height=17mm]{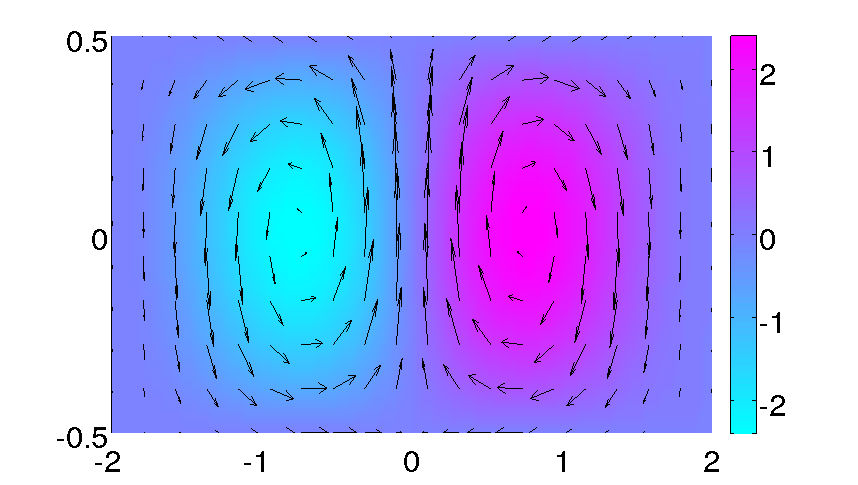}
\end{minipage}
}
\ece

\caption{{\small (a) Partial bifurcation diagram of \eqref{e:rbconv}
    with \reff{e:bcnoslip}. (b) sample solutions ($\psi$, and 
arrows indicating the fluid flow) from (a). See also {\tt arrowplot.m} for how to 
produce the quiver plots. 
\label{rbf1}}}
\end{figure}

For the no-slip case \eqref{e:bcnoslip}, the resulting bifurcation
diagram is plotted in Fig.~\ref{rbf1}, which corresponds
to the sketch Fig.~2 in \cite{HirschKnob97}. No secondary 
bifurcations are found up to $R=900$. 

\begin{figure}[ht!]
\bce 
{\small 
\begin{minipage}{40mm}
(a)  bifurcation diagram\\ 
\ig[width=40mm,height=34mm]{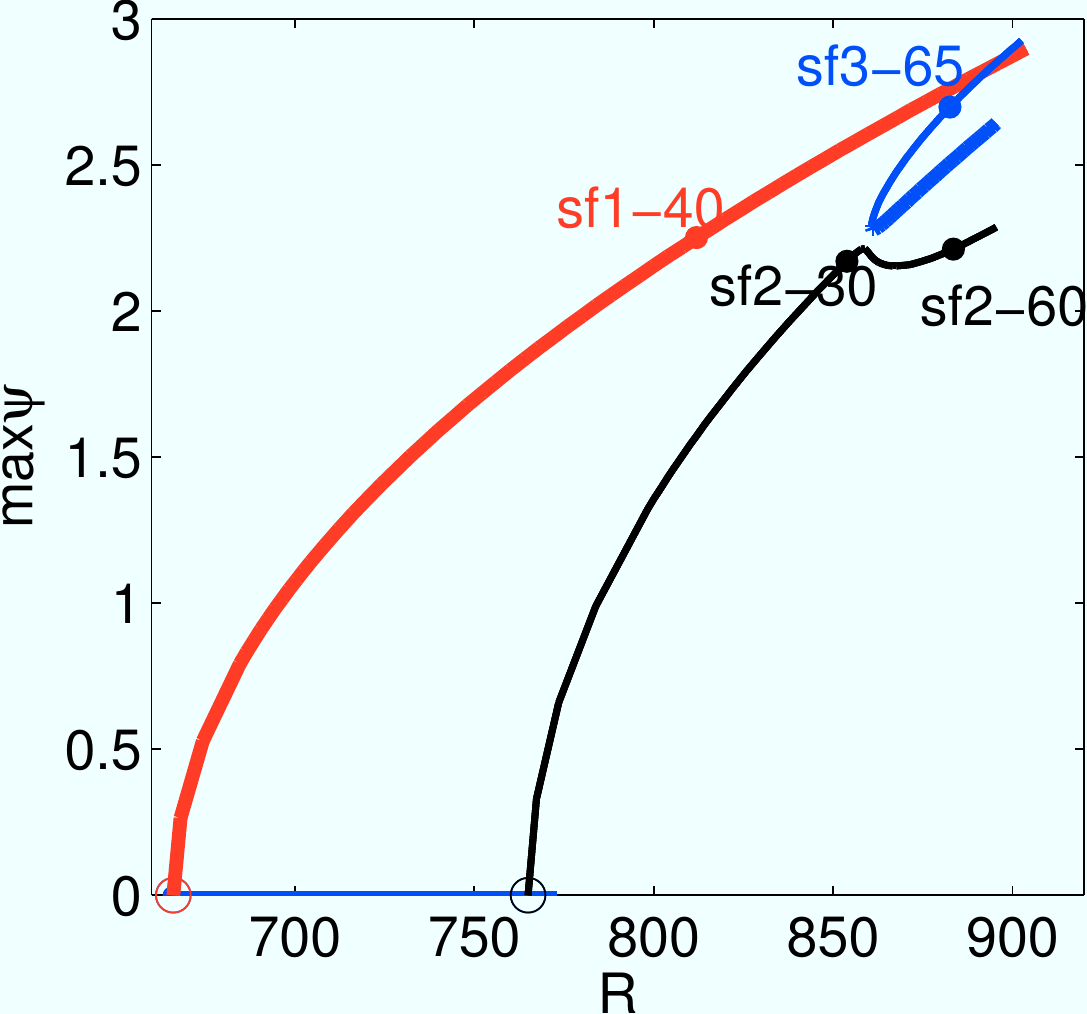}
\end{minipage}
\begin{minipage}{60mm}
(b) sf1-40 (top) and sf2-30 (bottom)\\ 
\ig[width=60mm,height=17mm]{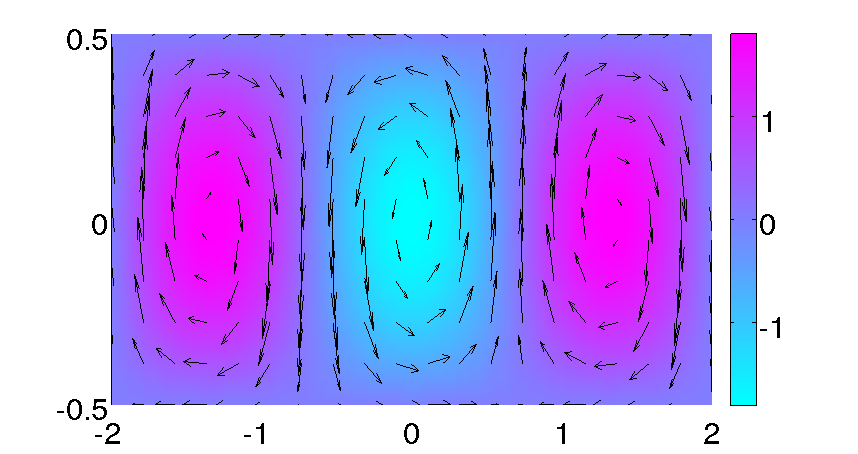}
\ig[width=60mm,height=17mm]{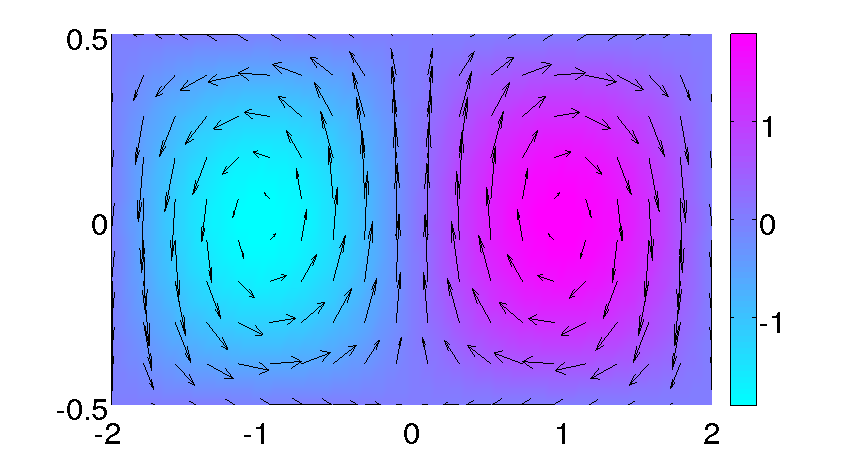}
\end{minipage}
\begin{minipage}{60mm}
(c) sf2-60 (top) and sf3-65 (bottom)\\ 
\ig[width=60mm,height=17mm]{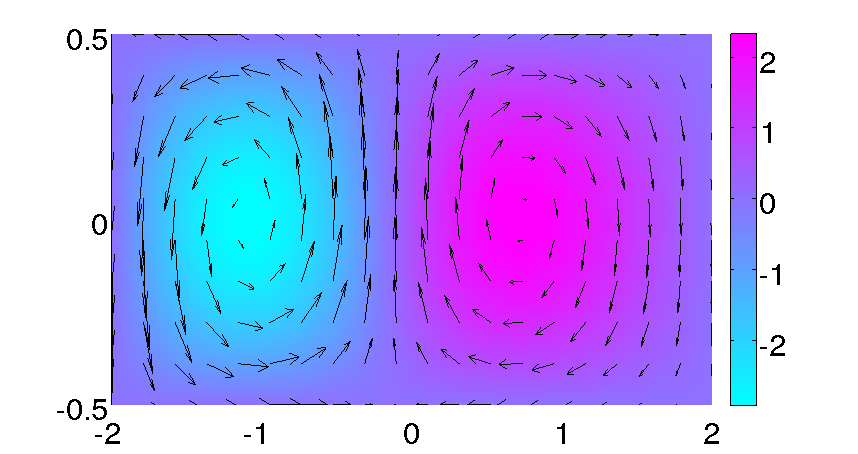}
\ig[width=60mm,height=17mm]{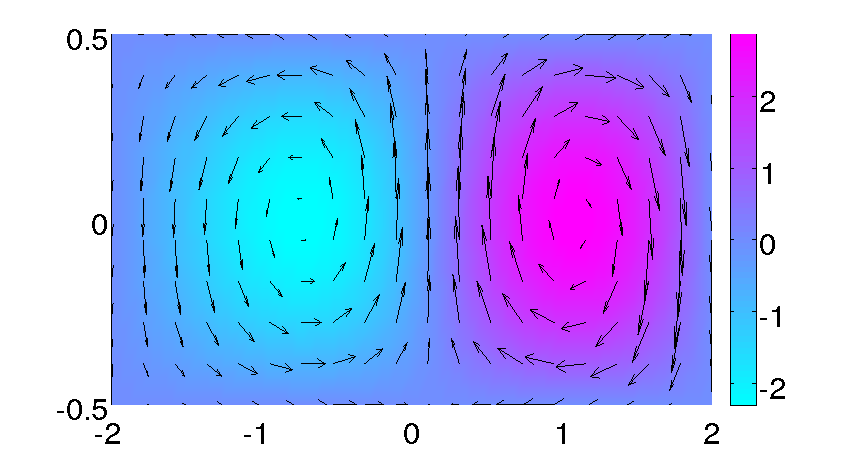}
\end{minipage}
}
\ece

\caption{{\small (a) Partial bifurcation diagram of \eqref{e:rbconv}
    with stress-free b.c.~and sample
    solutions ($\psi$) from (a). Here {\tt sf1-40} and {\tt sf2-30} are
    approximately symmetric, while {\tt sf2-60} and {\tt sf3-65}, 
 generated in a (numerically) imperfect pitchfork around $R=860$, are not.}
\label{rbf2}}
\end{figure}

For stress-free \bcs\ we obtain the
bifurcation diagram in Fig.~\ref{rbf2}(a), which
corresponds to the case $b'^2>a'^2, b'>0$ in Fig.~3 of
\cite{HirschKnob97}. Here the secondary symmetry breaking pitchfork
from \cite{HirschKnob97} is turned into an imperfect pitchfork. 
The $x\to-x$ reflection
symmetry is broken by the triangle data of the mesh (here we use {\tt
  poimesh}) and the stiff-spring approximation of the boundary
conditions. We have located the stable branch {\tt s}  
of the imperfect pitchfork by 
time-integrating
\footnote{which is not equivalent to the time integration of the time-dependent
  Boussinesq equations} with {\tt tint} from the unstable branch in the suitably
chosen unstable direction. 

The demos run on a rather coarse mesh 
of $100\times 25$ gridpoints, because even with assembled Jacobians 
the calculations are rather slow due to a non--simple structure 
of the Jacobians. Thus, we use {\tt pmcont} for the 
bifurcating {\tt q} and {\tt r} branches, which gives a huge speed 
advantage and works remarkably 
well even directly after bifurcation from the trivial branch, 
where long predictors are 
far off the actual branches.\footnote{E.g., in Fig.~\ref{rbf1} 
{\em all} points are calculated from the essentially ``vertical'' 
predictors at bifurcation! To check that this does not miss 
any (possibly imperfect) bifurcation we compared with {\tt cont} 
with small {\tt ds} 
and obtained the same branches but much slower.} In Fig.~\ref{rbf2}, however, 
we switch back to {\tt cont} when approaching the imperfect pitchfork since 
{\tt pmcont} tends to switch from the {\tt r} to the {\tt s} branch via 
long predictors. 

Increasing the number of meshpoints brings the
diagram closer to a symmetry breaking pitchfork. 
In fact, for both boundary conditions,
qualitatively the same bifurcation diagram can be found for rather
coarse meshes, but the location of the branches can be off by order
100 in $R$.

\subsection{Von K\'arm\'an description of the buckling of plates ({\tt vkplate})}
The \vKar\ equations
\huga{
\begin{split}
-\Delta^2 v-\lam\pa_x^2 v+[v,w]&=0, \quad 
-\Delta^2 w-\frac 1 2 [v,v]=0, 
\end{split}
\label{vk1} 
}
can be derived to describe the deformation of an 
elastic (rectangular) 
plate $\Om=[-l_x,l_x]\times [-l_y,l_y]\subset\R^2$ under compression. 
Here $v:\Om\ra\R$ is the out of plane deformation, $w:\Om\ra\R$ is 
the Airy stress function, $\Delta^2=(\pa_x^2+\pa_y^2)^2$ is the 
squared Laplacian, $\lam$ is the compression parameter, and the 
bilinear form $[\cdot,\cdot]$ is given by 
$$
[v,w]:=v_{xx}w_{yy}-2v_{xy}w_{xy}+v_{yy}w_{xx}. 
$$
There are a number of choices for the boundary conditions for 
\reff{vk1}. For $v$ one can choose for instance between 
(in the notation from \cite{gop97})
\hualst{
\text{I(v)}:\quad & v=\Delta v=0 \text{ on }\pa \Om, \quad 
\text{(simply supported)}, \\
\text{II(v)}:\quad & v=\Delta v=0 \text{ on } y=\pm l_y, 
v=\pa_n v=0\text{ on } x=\pm l_x, \\
& \text{(simply supported on the sides, clamped at the ends)}\\
\text{III(v)}:\quad &v=\pa_n v=0 \text{ on } \pa\Om, \quad 
\text{(clamped on whole boundary)}.
}
Similarly, for $w$ we may consider, on $\pa\Om$,  
\hualst{
\text{I(w)}:\ w=\Delta w=0,\quad 
\text{II(w)}:\ \pa_n w=\pa_n(\Delta v)=0,\quad 
\text{III(w)}:\ w=\pa_n w=0. 
}

Clearly, for all \bcs-combinations and all $\lam$ 
the trivial state $v=w=0$ is a solution. 
Mathematically, the combination a) (I(v),I(w)) (sometimes as a whole 
called simply supported) is most simple because it allows 
an easy explicit calculation of bifurcation points from the trivial branch. 
However, \cite{gs79} argues that physically 
the combinations b) (II(v),I(w)) or c) (II(v),II(w)) are more reasonable, 
and various combinations and modifications have been studied since, see 
\cite{cgm00} and the references therein for an overview. 

Here we focus on case b) since this yields secondary bifurcations, called 
``mode jumping'' in this field.  The other cases can be handled quite similarly 
and, e.g., a) is in fact slightly simpler. The aim is to show how 
\reff{vk1} can be put into \pdep\ 
and thus recover a number of interesting bifurcations. 

Clearly, the first idea to set up 
\reff{vk1} would be to introduce auxiliary variables $\Delta v, \Delta w$ 
and set 
$$
u=(u_1,u_2,u_3,u_4)=(v,\Delta v,w,\Delta w)
$$
to obtain the (quasilinear elliptic) system  
$$
\bpm -\Delta& 1& 0&0\\
-\lam\pa_x^2&-\Delta&0&0\\
0&0&-\Delta&1\\
0&0&0&-\Delta\epm u
-\bpm 0\\ -[u_1,u_3]\\0\\\frac 1 2 [u_1,u_1]\epm =0, 
$$
for instance in case a) with homogeneous Dirichlet BC 
$u_1=u_2=u_3=u_4=0$. The problem with this formulation in \pdep\ are 
the derivatives $\pa_x^2 u_1,\ldots,\pa_x\pa_y u_3$ in the nonlinearity. 
In principle, these can be obtained from calling 
{\tt pdegrad}, {\tt pdeprtni}, and {\tt pdegrad} 
again\footnote{e.g., 
{\tt [u1xt,u1yt]=pdegrad(p.points,p.tria,u(1:p.np));  
u1x=pdeprtni(p.points,p.tria,u1xt); 
[u1xx,u1xy]= pdegrad(p.points,p.tria,u1x);} could be used to 
calculate (approximate) $\pa_x^2 u_1$}. 
However, the first problem is that this introduces some averaging 
into the second derivatives, in particular at the boundaries. 
The second problem is that with this approach we have no easy way to 
generate the Jacobian of $G$ since {\tt pdegrad/ pdeprtni} 
neithers fit to matrix assembling nor to numerical 
differentiation.\footnote{For the latter the next-next-neighbor 
effect of {\tt pdegrad/ pdeprtni} does not comply with the Jacobian 
stencil, as explained in Remark \ref{jacrem}.}

Thus, here we choose to introduce additional auxiliary variables, i.e., set 
$$
u=(v,\Delta v, w,\Delta w,\pa_x^2 v,\pa_y^2 v,\pa_{x}\pa_y v, 
\pa_x^2 w,\pa_y^2 w,\pa_{x}\pa_y w)\in \R^{10}.
$$
For instance, $u_5=\pa_x^2 u_1$ can then be simply added 
as a linear equation $-\pa_x^2 u_1 +u_5=0$ in the \pdep\ 
formulation\footnote{see below for the BC for $u_5,\ldots,u_{10}$}.
However, since this way we get a number of indefinite equations, 
in particular the mixed derivatives $-\pa_{x}\pa_y u_1+u_7=0$ and 
$-\pa_{x}\pa_y u_3+u_{10}=0$, here we use an ad hoc regularization 
and set $-\pa_x^2 u_1+(1-\del\Del) u_5=0$ with small $\del>0$
(i.e., $\del=0.05$ numerically) 
and similarly for $u_6,\ldots,u_{10}$. Thus, instead of 
\reff{vk1} we now really treat the problem 
\huga{\begin{split}
-\Delta^2 v-\lam\pa_x^2 v+\bigl(Sv_{xx}Sw_{yy}-2Sv_{xy}Sw_{xy}
+Sv_{yy}Sw_{xx}\bigr)&=0, \\
\quad 
-\Delta^2 w-\bigl(Sv_{xx}Sv_{yy}-Sv_{xy}Sw_{xy}\bigr)&=0 
\end{split}\label{vkreg}
}
with the smoothing operator $S=(1-\del\Delta)^{-1}$. 
However, for small $\del$, comparison of our results with the 
literature shows that the regularization plays no qualitative or even 
quantitative role (in the parameter regimes we consider). 

Thus, we now have a 10 component system, and to illustrate its 
implementation in \pdep\ we write it in the form $(-C+A)u-f=0$ with 
$$
f=(0, -(u_5u_9-2u_7u_{10}+u_6u_8), 0, u_5u_6-u_7^2, 0, 0, 0,0,0,0)^T, \quad\text{and} 
$$
\arraycolsep4pt
\def\Delti{{\tilde{\cal D}}}
$$
{\footnotesize
-C+A=\bpm -\Delta_{1}&1^{11}& & & & & & & &  \\
-\lam\pa_x^2{}_5&-\Delta_{45}& & & & & & &  &\\
&&-\Delta_{89}&1^{33}& & & & & &\\
&&&-\Delta_{133}& & & & & & \\
-\pa_{x}^2{}_{17}&&&&-\Delti_{177}^{45}& & & &  &\\
-\pa_{y}^2{}_{21}&&&&&-\Delti_{221}^{56}&&&  &\\
-\pa_{x}\pa_y {}_{25}&&&&&&-\Delti_{265}^{67}&&&\\
&&-\pa_{x}^2{}_{109}&&&&&-\Delti_{309}^{78}&&\\
&&-\pa_{y}^2{}_{113}&&&&&&-\Delti_{353}^{89}&\\
&&-\pa_{x}\pa_y{}_{117}&&&&&&&-\Delti_{397}^{100}
\epm. 
}
$$
Here, (for layout reasons) $\Delti=\del\Delta+1$, and 
the subscripts $1,5,17,\ldots$ denote the 
starting positions of the respective $2\times 2$ tensor 
stored in the ``400 rows vector'' {\tt c}.\footnote{I.e., $-\Delta_{1}$ 
yields $c_1=[1;0;0;1]$ stored in positions 
1 to 4 in {\tt c}, ${\pa_x^2}_5$ yields $c_2=[1;0;0;0]$ stored in positions 
5 to 8 in {\tt c}, and so on.}
The superscripts $11,33,\ldots$ denote the 
positions in the ``100 rows vector'' {\tt a}, and 
for $\Delti$ subscripts refer to $\del\Delta$ and superscripts to $+1$. 
See {\tt vkf.m}. 
Similarly, it is now rather easy to put the linearization $f_u$ 
into \pdep, i.e., the second and fourth row of $f_u$ as a $10\times 10$ 
matrix read
$$
\begin{array}{llcccccccccc}
f_u, \text{ 2nd row:}\quad
&(0&0&0&0&-u_9{}^{42}&-u_8{}^{52}&2u_{10}{}^{62}&-u_6{}^{72}&-u_5{}^{82}&
2u_7{}^{92}&)\\
f_u, \text{ 4th row:}\quad
&(0&0&0&0&u_6{}^{44}&u_5{}^{54}&-2u_7{}^{64}&0&0&0&). 
\end{array}
$$
Here again the superscripts give the positions in {\tt fu}.  
Of course, the full ${\tt fu}=f_u-a$ also contains the constant coefficient 
terms at positions 11, 33, 45 etc from $A$; see {\tt vkjac.m}. 

It remains to encode the boundary conditions. First note that 
$v=0$ and $w=\Delta w=0$ imply 
\hualst{
&u_5=v_{xx}=0\text{ and } u_6=v_{yy}=0\text{ on horizontal edges, and}\\ 
&u_6=v_{yy}=0 \text{ on vertical edges, but no condition for $u_5=v_{xx}$, and}\\
&u_8=w_{xx}=0 \text{ and }u_9=w_{yy}=0\text { on all edges.}
}
For $u_5$ on the vertical edges and $u_7, u_{10}$ on all edges 
we take homogeneous Neumann boundary conditions. 
To put this into \pdep\ via \reff{e:gnbc} we 
thus need two boundary matrices $q^h$ and $q^v$. For the horizontal 
boundaries ($y=\pm l_y$), 
$q^h$ has diagonal $q^h_d=\bpm s & s &s & s& s& s& 0& s& s& 0\epm $. 
For the vertical boundaries ($x=\pm l_x$) 
$q^v$ has diagonal $q^v_d=\bpm 0 & 0& s& s& 0& s& 0& s& s& 0 \epm$ and 
additionally $q^v_{2,1}=s$, where $s=10^3$ stands for the stiff spring constant. 
Positions 7 and 10 in $q^h_d$ and 
$q^v_d$ give the Neumann BC for $u_7,u_{10}$, while 
the top left $2\times 2$ block $\bpm 0&0\\ s&0\epm$ in $q_v$ 
gives $\pa_n u_1=0$ via the first row and $u_2=0$ via the second row. 

The (analytical) calculation of bifurcation points from $(v,w)=0$ in case b) 
is rather tedious, see \cite{gs79}. There, motivated 
by mode--jumping, the particular interest 
is in (the lowest) double bifurcation points, which yields $l=\sqrt{k(k+2)}$ 
with eigenfunctions  
$w_1(x,y)=\left(\frac{k+2}k\sin(k\frac x l)-\sin((k+2)\frac x l)\right)\sin(y)$ 
resp.~$w_2(x,y)=\left(\cos(k\frac x l)-\cos((k+2)\frac x l)\right)\sin(y)$ 
(over the domain $[0,l\pi]\times [0,\pi]$). 
The first bifurcation is then obtained for $k=1$, hence $l=\sqrt{3}$. 
The idea is to perturb $l$ slightly which may lead to secondary 
bifurcations between branches coming originally from the same $\lam$. 

\begin{figure}[htbp]
\bce 
{\small 
\ig[width=40mm,height=35mm]{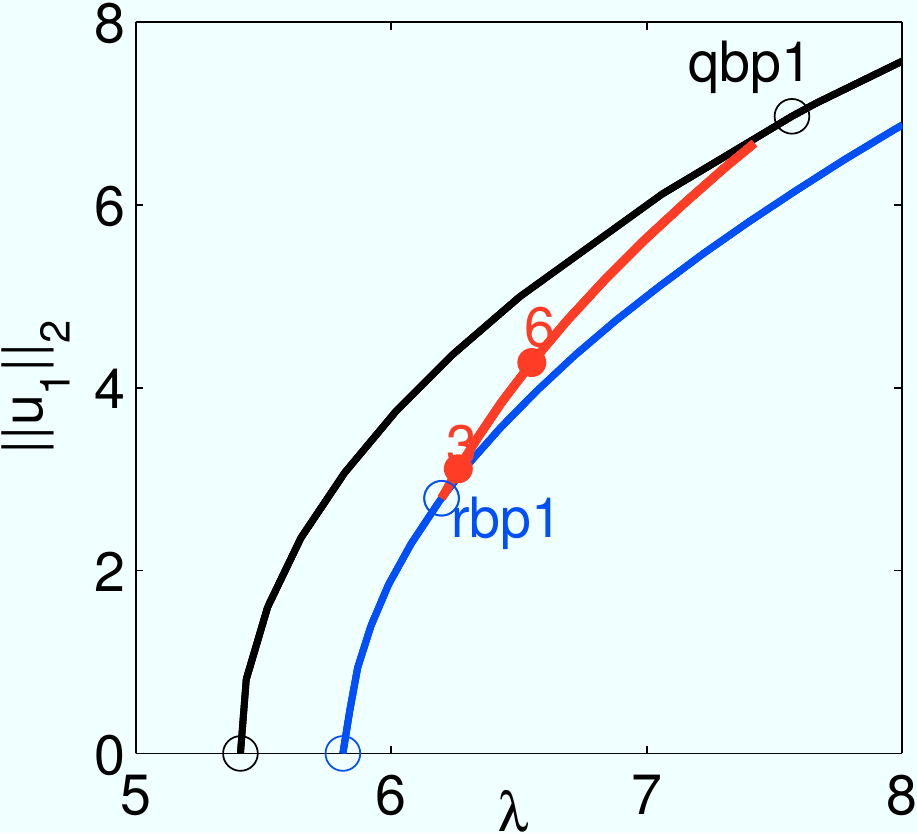}
\ig[width=38mm,height=35mm]{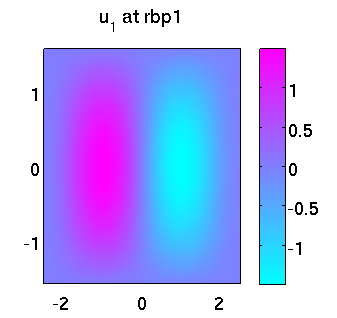}
\ig[width=38mm,height=35mm]{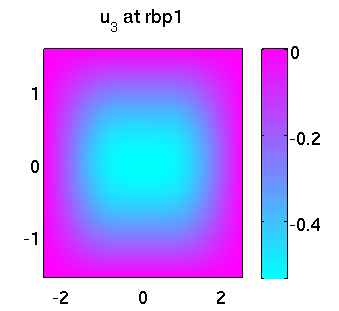}
\ig[width=38mm,height=35mm]{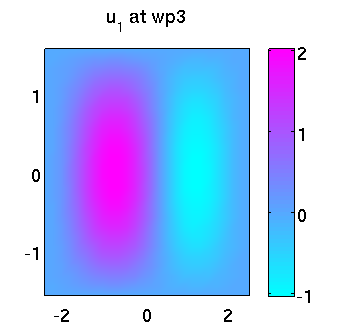}

\ig[width=38mm,height=35mm]{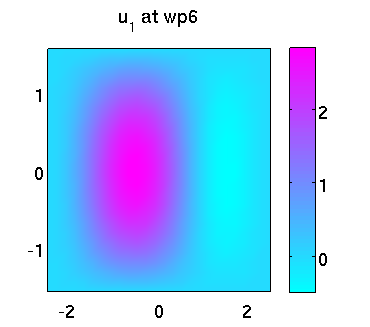}
\ig[width=38mm,height=35mm]{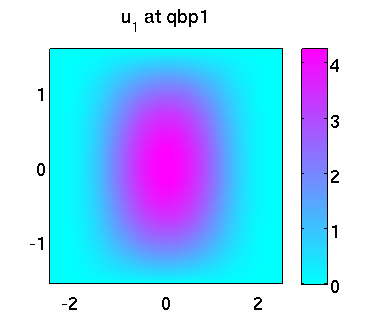}
\ig[width=38mm,height=35mm]{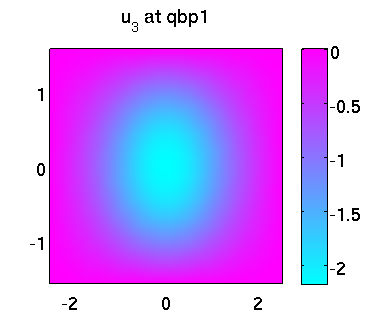}
\ig[width=38mm,height=35mm]{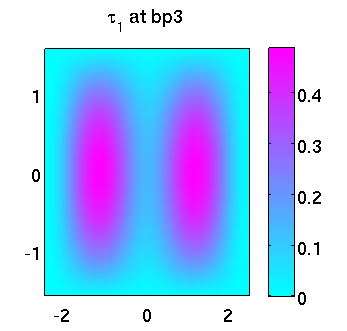}
}
\ece

\vs{-4mm}
\caption{{\small Secondary (``mode jumping'') 
bifurcation (w-branch, red) in the (regularized) partially clamped plate 
\reff{vkreg}: 
$\max u_1$ over $\lam$, selected plots of $u_1$ and $u_3$  from the bifurcation diagram, 
and $\tau_1$ at the third bifurcation ($\lam\approx 9.1$) from the trivial branch. 
$l_y=\pi/2$, $l_x=4\pi/5$, regular mesh with $25\times 25$ points (1250 triangles). 
Error estimate 
$0.3$ at, e.g., {\tt rbp1}. By mesh refinement we can obtain an 
error-estimate$\approx 0.045$ with nt$=14916$. Then, however 
a typical step takes a couple of minutes, where about 80\% of the time is 
spent in {\tt blss} or {\tt lss} (standard setting). We expect that this 
can be optimized considerably, 
but here we content ourselves with the ``proof of principle'' setup 
for the 10 components system for \reff{vkreg}. 
 }\label{vkpcpic}} 
\end{figure}

Putting all these ideas together we  indeed get a secondary bifurcation 
between the first two primary branches, see Fig.~\ref{vkpcpic}. A number 
of further bifurcations from the trivial branch is also detected and can 
be followed. However, in the tutorial run {\tt vkcmds} we use a rather 
coarse mesh with 1250 triangles, which should be refined before following 
higher bifurcations. 

\section{Discussion}\label{discsec}
Clearly, 
numerical continuation and bifurcation analysis for 2D elliptic systems 
poses  additional challenges compared to algebraic equations or 1D BVP, 
partly of course due to the more demanding numerics, but in particular also due 
to the typically very rich solution and bifurcation structure. With \pdep\ we believe 
to provide  a {\em general} tool that works essentially out-of-the-box 
also for non-expert users and allows to start exploring such systems and the 
rich zoo of their solutions. 
Of course, in many respects this is just a first step, and probably the 
main entries on our {\bf to--do--list} are: 
\ben 
\item Implement some more general bifurcation handling, in particular 
bifurcation via nonsimple eigenvalues as these are quite ubiquious 
in 2D systems due to various symmetries. 
\item Implement some (genuine) multi--parameter continuation. For instance, 
the bifurcation to travelling waves generically requires 
a second parameter $\gamma$ (the wave speed) to adapt, and consequently 
we need to further extend the ``extended system'' \reff{esys} by 
one more equation, the ``phase condition''.\footnote{\label{eesfn}
More generally, adding some constraints to \reff{esys} will also be useful 
to remove some degeneracies as, e.g., the phase invariance in \S\ref{gpsec}.}
We believe that our set--up of \pdep\ is sufficiently modular and 
transparent such that this and similar adaptions 
will pose no implementation 
problems, but for now we confine ourselves to the 
basic one--parameter continuation and simple bifurcations. 
\een

\renewcommand{\refname}{References}
\renewcommand{\arraystretch}{1.05}\renewcommand{\baselinestretch}{1}
\small
\bibliographystyle{plain}\bibliography{pdepath.bib}

\end{document}